\documentclass[11pt, leqno,amscd, psamsfonts]{amsart}
\usepackage{amsthm, amsmath}
\usepackage{amssymb} 
\usepackage{amscd}
\usepackage[matrix,arrow,frame,curve]{xy}
\newdir{ (}{{}*!/-5pt/@^{(}}

\newtheorem{thm}{Theorem}[section]
\newtheorem{prop}[thm]{Proposition}
\newtheorem{lemma}[thm]{Lemma}

\newtheorem{Definition}[thm]{Definition}
\newtheorem{Remarknumb}[thm]{Remark}

\newtheorem{Remark}[thm]{Remark}

\newtheorem{conjecture}[thm]{Conjecture}
\newtheorem{cor}[thm]{Corollary}

\newcounter{ex}[section]

\numberwithin{equation}{section}
\numberwithin{thm}{section}

\newcommand{\cal}{\mathcal}

\newcommand{\Ext}{{\rm Ext}}

\newcommand{\Lotimes}{\overset{{\rm L}}{\otimes}}

\newcommand{\mi}{{\hbox{\small -}}}
\newcommand{\nExt}{n{\hbox{\small -}}{\rm Ext}^1}
\newcommand{\CC}{{\mathcal C}}
\newcommand{\Cl}{{\rm Cl}}
\newcommand{\Gal}{{\rm Gal}}

\newcommand{\R}{{\bf R}}

\newcommand{\Q}{{\bf Q}}

\newcommand{\ep}{\epsilon}

\newcommand{\Lder}{\buildrel{\rm L}\over{\otimes}}
\newcommand{\Hom}{{\rm Hom}}

\newcommand{\Gg}{{\mathcal G}}
\newcommand{\Hh}{{\mathcal H}}
\newcommand{\Gm}{{{\bf G}_m}}
\newcommand{\Z}{{\bf Z}}
\newcommand{\F}{{\mathcal F}}
\newcommand{\uno}{\underline {O}}
\newcommand{\ti}{\tilde}
\newcommand{\Spec}{{\rm Spec }\, }
 \renewcommand{\O}{{\mathcal O}}

\newcommand{\Sfr}{{\mathfrak S}}
\newcommand{\M}{{\mathcal M}}
\renewcommand{\L}{{\mathcal L}}

\newcommand{\Gr}{{\rm G}}

\newcommand{\Kr}{{\rm K}}

\newcommand{\lon}{{\longrightarrow}}
\newcommand{\lo}{{\rightarrow}}

\newcommand{\fp}{{{\bf F}_p}}

\newcommand{\Pic}{{\rm Pic}}
\newcommand{\PP}{{\mathcal P}}

\newcommand{\QQ}{{\mathcal Q}}

\newcommand{\loniso}{\buildrel \sim \over \lon}
\newcommand{\noi}{\noindent}

\newcommand{\ord}{{\rm ord}}

\def\thfill{\null\nobreak\hfill}

\def\endproof{\thfill\vbox{\hrule
  \hbox{\vrule\hbox to 5pt{\vbox to 5pt{\vfil}\hfil}\vrule}\hrule}}

\renewcommand{\P}{{\cal P}}

\begin{document}

\title[Galois modules]{Galois modules, ideal class groups\\ and cubic structures}
\author[G. Pappas]{Georgios Pappas*}
\date{\today}
\thanks{*Partially supported by NSF
grant DMS02-01140 and by a Sloan Research Fellowship.}
\address{Dept. of
Mathematics\\
Michigan State
University\\
E. Lansing\\
MI 48824-1027\\
USA}
\email{pappas@math.msu.edu}

\begin{abstract}
We establish a connection between 
the theory of cyclotomic ideal class groups and 
the theory of ``geometric" Galois modules and obtain 
results on the Galois module structure of coherent cohomology groups
of Galois covers of varieties over $\Z$. In particular, we show that an
invariant that measures the obstruction to the existence of 
a virtual normal integral basis for the coherent cohomology of such covers 
is annihilated by a product of  certain Bernoulli numbers with
orders of even ${\rm K}$-groups of $\Z$. We also show that the existence of 
such a normal integral basis is
closely connected to the truth of the Kummer-Vandiver
conjecture for the prime divisors of the degree of the cover.
Our main tool is a theory of ``hypercubic structures" for
line bundles over group schemes.

\end{abstract}
\maketitle

 \centerline{\sc Contents}
\noindent  \S 1. Introduction\hfill\\
 \S 2. Torsors\hfill\\
 \S 3. Hypercubic structures\hfill\\
 \S 4. Multiextensions\hfill\\
 \S 5. Differences and Taylor expansions\hfill\\
 \S 6. Multiextensions and abelian sheaves\hfill\\
 \S 7. Multiextensions of finite multiplicative group schemes\hfill\\
 \S 8. Reflection homomorphisms\hfill\\
 \S 9. The hypercubic structure on the determinant of cohomology\hfill\\
 \S 10. The multiextension associated to a $\Z/p\Z$-torsor\hfill\\
 \S 11. Galois module structure\hfill\\
  \ \ Appendix

\bigskip
\vfill\eject
\section{Introduction}

For a given finite group $G$ and an arithmetic variety
$Y$, we associate to each $G$-torsor $X$ over $Y$ an equivariant
Euler characteristic $\chi^P(\O_X)$ in $\Kr_0(\textbf{Z}[G])$ (the
Grothendieck group of finitely generated projective $\Z[G]$-modules;
see both below and also Section 11). In this paper we develop
powerful new methods for determining Euler characteristics of such
torsors by extending work of Grothendieck, Breen and other authors on cubic
structures. Our main goal   is to establish a connection
between the problem of determining $\chi^P(\O_X)$ and the classical theory
of cyclotomic ideal class groups. 
We then exhibit annihilators
of such Euler characteristics $\chi^P(\O_X)$ constructed from the
divisors of certain Bernoulli numbers and the orders of even
Quillen ${\Kr}$-groups of $\Z$ and describe a relation with the Kummer-Vandiver conjecture. 
To illustrate the power of these
techniques, we mention that our main theorems show in particular that
if $G$ is abelian and if either $\dim(Y)\leq 4$ or the cover $X\to Y$ 
is of ``Albanese type" (see \S\ref{paralb}) then $2\cdot \chi^P(\O_X)$
is  the class of a free module in $\Kr_0(\Z[G])$; furthermore if in
addition we suppose that $G$ has odd order then  we show
that $\chi^P(O_X)$ is in fact   the class of a free module.
Note that these results may be seen as providing higher dimensional counterparts
of a result of M. Taylor, which shows if $Y$ is the spectrum of a ring
of algebraic integers of a number field $K$, then the ring of
integers of a non-ramified abelian extension of $K$ always has a
normal integral basis. It is also important to note that the
techniques developed here, together with the Grothendieck-Riemann-Roch 
theorem, can be used to determine such Euler characteristics even when the 
$G$-cover $X\to Y$ has some ramification  - see [CPT] for further developments in this
direction. Although a number of results in this paper are valid
for non-abelian groups, our approach provides the strongest
results when the group is abelian; this article therefore
naturally raises some fundamental questions for the actions of
non-abelian groups.

To better understand the context of our work, and explain the relation with ideal class groups,
we first review
 some standard cyclotomic theory.  The Stickelberger relations in the ideal class groups
of cyclotomic fields are undoubtedly 
among the cornerstones of classical 
algebraic number theory.  Hilbert gave the first ``Galois module theoretic" proof of 
these relations for cyclotomic fields generated by a root is
 of unity of prime order by showing that, in this case, they can deduced from the fact that,
 for $n$ square free, 
the primitive $n$-th root of unity $\zeta_n$ gives a normal integral basis
for the ring of integers $\Z[\zeta_n]$ in the cyclotomic field $\Q(\zeta_n)$
(Zahlbericht, Satz 136).
Recall that a normal integral basis for the ring of integers
$\O_K$ in the Galois extension $K/\Q$ is an algebraic integer  $a\in \O_K$
such that $\{\sigma(a)\}_{\sigma\in {\rm Gal}(K/\Q)}$ give a $\Z$-basis of $\O_K$;
alternatively, we can say that $a$ is a normal integral basis when
it is a generator of $\O_K$ as a
$\Z[{\rm Gal}(K/\Q)]$-module. This proof of the 
Stickelberger relations was later extended by Fr\"ohlich to composite
$f$ [F1]. Conversely, the Stickelberger relations  play a central role
 in the study of the problem of
existence of a normal integral basis: Indeed, Fr\"ohlich uses the prime factorization
of the Gauss sum which underlies the Stickelberger relation, to establish his fundamental relation 
between Galois Gauss sums and tame local resolvents (see [F2] Theorem 27 and VI \S 4). 
In the present paper,
we  extend this classical connection between the structure of ideal 
class groups of cyclotomic fields and the theory of ``additive" Galois modules
in a new direction 
by considering Galois covers of higher dimensional algebraic varieties over the integers. 

Let $Y$ be a projective algebraic variety over $\Z$
(i.e an integral scheme
which is projective and flat over $\Spec(\Z)$). We will
consider finite Galois covers $\pi: X\to Y$ with group $G$.
By definition, these covers are everywhere unramified, i.e $\pi$ is a $G$-torsor 
(see \S \ref{torsors}). The Galois modules we are considering are the 
finitely generated $G$-modules ${\rm H}^i(X, \O_X)$; more generally we will also
consider ${\rm H}^i(X, \F)$ where $\F$ is a $G$-equivariant coherent sheaf on $X$
(see \S \ref{torsors2}). These (coherent) cohomology groups can be calculated as the
cohomology groups of a complex ${\bf R}\Gamma(X,\F)$ in the derived category
of complexes of  modules over the group ring $\Z[G]$. 
The appropriate generalization of the question of the existence
of a normal integral basis is now the question: 

Is the complex ${\bf R}\Gamma(X,\F)$ in the derived category
of complexes of $\Z[G]$-modules isomorphic to a bounded complex of finitely 
generated {\sl free} $\Z[G]$-modules?

This question has first been considered in a geometric context
by Chinburg [C] for more general (tamely ramified) covers. He showed that the obstruction to 
a positive answer is an element $\bar\chi^P(\F)$ (``the projective equivariant Euler characteristic")
in the class group $\Cl(\Z[G])=\Kr_0(\Z[G])/\pm\{\hbox{\rm free classes}\}$ of finitely generated projective $\Z[G]$-modules. 
In [P] we stated the expectation that, for unramified covers, this obstruction
should vanish ``most of the time" and we obtained results which point in this direction when $Y$ is 
an arithmetic surface (i.e of dimension $2$). In this paper, we extend                                                         
the results of [P] to varieties of arbitrary dimension. In particular, we 
show that class field theory combined with results on  
ideal class groups of cyclotomic fields and their relation with the Quillen $\rm K$-groups ${\rm K}_m(\Z)$
(Herbrand's theorem, work of Soul\'e, Kurihara
 and others) implies general vanishing theorems for the classes $\bar\chi^P(\F)$. 
In fact, as we will explain below, our results indicate that it is reasonable to expect that 
a positive answer to the above question for all unramified Galois covers $X\to Y$ 
of prime order $p>{\rm dim}(Y)$ is equivalent to
the truth of the Vandiver conjecture for $p$.
Vandiver's conjecture (sometimes also attributed to Kummer) for the prime number $p$ 
is the statement that $p$ does not divide the class number $h^+_p=\#{\rm Cl}(\Q(\zeta_p+\zeta_p^{-1}))$.
This has been verified numerically for all $p<12\cdot 10^6$ [BCEM]. 
However, there is widespread doubt about its truth in general (see the discussion in Washington's book [W], p. 158).
Up to date, this conjecture remains one of the most intractable problems in classical 
algebraic number theory. Note that a simpler connection between Vandiver's conjecture and 
the ``relative" Galois module structure of the ring of integers in the $\Z_p$-extension of $\Q(\zeta_p)$ 
is already known (see [KMi], [Gr]). Our results lie deeper and in a different direction.

In order to state our main theorem, we now introduce various numbers which feature in the annihilators
of our equivariant Euler characteristics. For a given finite group $G$ we set $\ep(G)={\rm gcd}(2,\# G)$.
Next recall that the $k$-th Bernoulli number $B_k$ is defined by the powers series identity:
${t}/{(e^t-1)}=\sum^{\infty}_{k=0}B_k {t^k}/{k!}$. 
 For a prime $p$, we denote by $\ord_p(a) $
the highest power of $p$ that divides $a\in\Z_{>0}$.  
For $n\geq 2$, let us set 
\begin{equation*}
e(n)=\begin{cases}
\ \ {\rm numerator\,}({B_{n}}/{n})&,\ \text{if $n$ is even},\\
\displaystyle{\prod_{\ p,\, p|h^+_p}}\ord_p(\#\Kr_{2n-2}(\Z))&,\ \text{if $n$ is odd},
\end{cases}
\end{equation*}
where $\Kr_{2n-2}(\Z)$ is the Quillen $\Kr$-group 
(by work of Borel [Bo], this is a finite group for $n>1$). 
We set 
\begin{equation*}
M_{n}(G)= \prod_{k=2}^n\bigl(\prod_{p |e(k)}\ord_p(\# G)\bigr)\ , \quad M'_n(G)= \prod_{s=1}^{[n/2]} 
\bigl(\prod_{p|e(2s)}\ord_p(\# G)\bigr).
\end{equation*}
Thus $M'_n(G)$ is determined solely by the order of the
group $G$ and the divisibility properties of Bernoulli numbers;
furthermore it follows at once from the definitions that of course
$M'_n(G)$ divides $M_n(G)$. 

Now let ${\mathcal M}_G$ be a maximal $\Z$-order in $\Q[G]$ which contains $\Z[G]$ 
and denote by $\Cl({\mathcal M}_G)$ the class group of finitely generated 
projective ${\mathcal M}_G$-modules. 
Tensoring   with ${\mathcal M}_G$ over $\Z[G]$ induces a 
group homomorphism $\Cl(\Z[G])\to \Cl( {\mathcal M}_G)$. 
Denote its kernel by ${\rm D}(\Z[G])$; this subgroup 
of $\Cl(\Z[G])$ is 
independent of the choice of maximal order ${\mathcal M}_G$ 
([F2] I \S 2). 

Let $\pi: X\to Y$ be a $G$-torsor with $h: Y\to \Spec(\Z)$ projective and flat of 
relative dimension $d$.  Suppose that 
$\F$ is a $G$-equivariant coherent sheaf on $X$.

\begin{thm} \label{main1}
  
a) The multiple $M_{d+1}(G)\cdot \bar\chi^P(\F)$ belongs to the kernel subgroup ${\rm D}(\Z[G])$.

b) Assume that all the Sylow subgroups of the group $G$ are abelian. Then 
 the multiple \hbox{${\ep(G)M_{d+1}(G)\cdot \bar\chi^P(\F)}$} is trivial in $\Cl(\Z[G])$.
\end{thm}
 
\begin{thm} \label{main2}
Assume that all the prime divisors of the order $\#G$ satisfy  Vandiver's conjecture.  
Let us set $C_{d+1}(G)={\rm gcd}(M'_{d+1}(G), 2((d+1)!!))$ with $(d+1)!!=(d+1)!\, d!\cdots 2!$. Then 
$C_{d+1}(G)\cdot \bar\chi^P(\F)$ belongs to the kernel subgroup ${\rm D}(\Z[G])$.
In particular, if all the prime divisors of $\#G$ satisfy 
Vandiver's conjecture and are larger than $d+1$
then $\bar\chi^P(\F)$ belongs to the kernel subgroup ${\rm D}(\Z[G])$. 
\end{thm}

\begin{thm} \label{main3}
Suppose that $G$ is abelian.

a) If $\pi: X\to Y$ is of Albanese type
(see \S\ref{paralb}) then $\bar\chi^P(\F)$ is trivial in $\Cl(\Z[G])$.

b) If $\pi: X_\Q\to Y_\Q$ is of Albanese type  
then the multiple $2((d+1)!!)\cdot \bar\chi^P(\F)$ is trivial in $\Cl(\Z[G])$.
\end{thm}

If $Y_\Q$ is smooth, has a rational point and $X_\Q$ is geometrically connected,
the condition that 
$X_\Q\to Y_\Q$ is of Albanese type means that the cover
is obtained by specializing an isogeny of the Albanese variety of $Y_\Q$
(see \S\ref{paralb}).

Note that for $d=0$, Fr\"ohlich's conjecture (shown by M. Taylor in [Ta])
implies that, for all finite Galois groups $G$, the class $\bar\chi^P(\F)$ 
is $2$-torsion and  belongs 
to the kernel subgroup ${\rm D}(\Z[G])$. When $d=1$ (the case of arithmetic surfaces), 
$M_{d+1}(G)=1$ for all $G$; in this case, Theorem \ref{main1}
was shown in [P]. In fact, we also have $M_{d+1}(G)=1$ for $d=2, 3$ 
and all $G$ (see Remark \ref{low}).
In general, if $\#G$ is relatively prime to 
$N_{d+1}:=\prod_{k=2}^{d+1}e(k)$ then $M_{d+1}(G)=1$. 
Note that by a result of Soul\'e on the size of $\Kr_{2m}(\Z)$ [So]
and standard estimates on Bernoulli numbers,
there is an effective (although doubly exponential) bound on $N_{d+1}$.

\begin{cor}\label{1.2}
Assume that either all the Sylow subgroups of $G$ are abelian and   we have 
$\ep(G)M_{d+1}(G)=1$ or that $G$ is abelian  
and $\pi: X\to Y$ is of Albanese type.
Then the complex ${\bf R}\Gamma(X,\F)$ 
is isomorphic in the derived category to a bounded complex of finitely 
generated free $\Z[G]$-modules. 
\end{cor}

Let $\ep(G,Y):={\rm gcd}(2,\#G, g(Y_\Q))$ where $g(Y_\Q)$
is the arithmetic genus of the generic fiber of $Y\to \Spec(\Z)$.
When we restrict our attention to $\F=\O_X$ the conclusions of 
Theorem \ref{main1} (b) and of Corollary \ref{1.2} can be slightly improved: 
in these statements, we can replace  $\ep(G)$ by
$\ep(G,Y)$. 

In a slightly different direction, we can use the above work to obtain
the following very down to earth result on  Galois modules:

\begin{cor}{\bf  (Projective integral normal basis)}\label{1.3}
Assume that either all the Sylow subgroups of $G$ are abelian and   we have 
$\ep(G)M_{d+1}(G)=1$ or that $G$ is abelian  
and $\pi: X\to Y$ is of Albanese type.  Also suppose in addition that $X$ is normal. Then
there is a graded 
commutative ring $\oplus_{n\geq 0}A_n$
with (degree $0$) $G$-action such that for each $n>0$, $A_n$ is a free 
$\Z[G]$-module and $X\simeq {\bf Proj}(\oplus_{n\geq 0}A_n)$
as schemes with $G$-action.
\end{cor}

To better explain the relation with Vandiver's conjecture we now assume
that $G=\Z/p\Z$, $p$ an odd prime. Then the kernel subgroup ${\rm D}(\Z[G])$ is trivial
and the class group $\Cl(\Z[G])$ coincides with the ideal class group $\Cl(\Q(\zeta_p))$.
For simplicity, we will restrict our discussion to $\F=\O_X$.
From Theorem \ref{main2}   above we obtain:
\begin{cor}
Assume in addition that $G=\Z/p\Z$, $p>d+1$, and that $p$ satisfies
Vandiver's conjecture. Then $\bar\chi^P(\O_X)=0$.
\end{cor}
Denote by $_pC$ (resp. $C/p$)
the kernel (resp. the cokernel) of multiplication
by $p$ on $C=\Cl(\Q(\zeta_p))$.
We use the superscript 
$(k)$ to denote the eigenspace of $_pC$ or $C/p$
on which $\sigma_a\in \Gal(\Q(\zeta_p)/\Q)$,
$\sigma_a(\zeta_p)=\zeta_p^a$, acts via multiplication by
$a^k$. Recall the classical 
``reflection" homomorphisms ([W] \S 10.2)
$$
R^{(i)}: \Hom((C/p)^{(1-i)}, \Z/p\Z)\to {}_pC^{(i)}.
$$
We continue to assume that $p>\dim(Y)$. Our results  actually 
produce elements 
$t_i(X/Y)\in\Hom((C/p)^{(1-i)}, \Z/p\Z)$, $1\leq i\leq d+1$,
such that
\begin{equation*}
\bar\chi^P(\O_X)=\sum_{i=1}^{d+1}R^{(i)}(t_i(X/Y)) 
\end{equation*}
in $\Cl(\Q(\zeta_p))=\Cl(\Z[G])$. By class field theory,
each   element $t_i(X/Y)$ corresponds to an unramified
$\Z/p$-extension of $\Q(\zeta_p)$. We may think of them 
as characteristic classes of the cover $X\to Y$; at least $t_{d+1}(X/Y)$
can be obtained from the cover of generic fibers $X_\Q\to Y_\Q$ 
in a very explicit manner (see \S \ref{multidet} where a relation
with a cohomological Abel-Jacobi map considered by
Bloch and Colliot-Th\'el\`ene--Sansuc is also explained). 
When $p$ satisfies Vandiver's conjecture, we have 
$R^{(i)}=0$. Hence, $\bar\chi^P(\O_X)=0$. We conjecture that
the elements $\{t_{d+1}(X/Y)\}$, $X\to Y$ ranging over all $\Z/p\Z$-torsors
with $Y\to \Spec(\Z)$ projective and flat of relative dimension $d$,
generate the group $\Hom((C/p)^{(-d)}, \Z/p\Z)$ (this has the 
 flavor 
of a refined version of a higher dimensional inverse Galois problem).
Since $R^{(i)}$ is injective for $i$ odd ([W]), this conjecture has the following 
(possibly vacuous) implication:
\begin{conjecture}
Suppose $p$ does not satisfy 
Vandiver's conjecture.
Then there is a $\Z/p\Z$-torsor
$\pi: X\to Y$ with $Y\to \Spec(\Z)$ projective and flat 
of relative dimension $d<p-1$ such that $\bar\chi^P(\O_X)\neq 0$.  
\end{conjecture}

The problem of determining the classes $\bar\chi^P(\F)$ is very subtle.
Indeed, these lie in the finite group ${\rm Cl}(\Z[G])$ and cannot
be calculated using the known ``Riemann-Roch type" theorems that usually
neglect torsion information. Instead, our basic tool  
is the notion of ``$n$-cubic structure" on line bundles 
over commutative group schemes. This is a generalization
of the notion of cubic structure (for $n=3$) which was introduced by Breen [Br].
Breen's motivation was to explain the properties of the trivializations
of line bundles on abelian varieties which are given by the theorem of the cube.
The starting point for us is the fact, shown by F.~Ducrot [Du], that the 
determinant of cohomology along a projective and flat morphism 
$Y\to S$ of relative dimension $d$ essentially supports a 
$d+2$-cubic structure (see \S \ref{hypercubedet}). When $G$ is abelian, we  deduce  that 
the square of the determinant $\det_{\Z[G]}{\bf R}\Gamma(X,\O_X)^{\otimes 2}$ 
gives a line bundle
over $G^D:=\Spec(\Z[G])$ which supports a $d+2$-cubic structure. 
We proceed to analyze line bundles with $n$-cubic structures over
the group scheme $G^D$. We show that they can be understood using 
``multiextensions". 
This is a notion which generalizes that of a biextension
and already appears in the work of Grothendieck [SGA6]. 
One of our main technical results is the ``Taylor expansion" of \S \ref{taylor}
that expresses line bundles with $n$-cubic structures in terms of
multiextensions. The Taylor expansion, when applied to the line 
bundle $\det_{\Z[G]}{\bf R}\Gamma(X,\O_X)^{\otimes 2}$,  substitutes for 
a (functorial) Riemann-Roch theorem ``without denominators":
it provides a formula for the $2((d+1)!!)$-th power of the determinant of the cohomology. 
(Notice that in the classical Riemann-Roch theorem 
multiplying both sides by $(d+1)!!$ clears all the denominators; see Remark \ref{RRexpl}.)
When $G=\Z/p\Z$, we see that
a $k$-multiextension over $G^D=\mu_p$ is given 
by an unramified $\Z/p\Z$-extension $L$ of $\Q(\zeta_p)$ 
such that  $\sigma_a\in {\rm Gal}(\Q(\zeta_p)/\Q)$ 
acts by conjugation on $\Z/p\Z={\rm Gal}(L/\Q(\zeta_p))$ via multiplication by $a^{1-k}$. 
In fact, the elements $t_i(X/Y)$, $1\leq i\leq d+1$,
mentioned above,
correspond to $i$-multiextensions associated naturally to
the $d+2$-cubic structure on $\det_{\Z[G]}{\bf R}\Gamma(X,\O_X)^{\otimes 2}$.
In general, class field theory implies 
that the group of multiextensions over $G^D$ can 
always be described using eigenspaces of
cyclotomic ideal class groups. We are able to bound 
the support of the group of multiextensions
using results on the corresponding eigenspaces of 
these class groups. This technique applies directly when $G$ is abelian
and to $\F=\O_X$;  the case of more general finite groups 
and general $\F$ follows from these results using Noetherian induction
and the localization methods of [P]. 

The theory developed in this paper also
applies to various other situations.
Most of the basic set-up easily 
extends to torsors for a general
finite and flat commutative group scheme over 
an arbitrary Dedekind ring;  
this leads to 
results on ``relative" Galois module structure.
In fact, the technique of cubic structures
provides the most precise general method for determining the Galois module
structure of abelian covers even when we allow some (tame) ramification.
Then one can not expect $\bar\chi^P(\F)$ to vanish. However, by 
combining the theory of the present paper
with the Grothendieck-Riemann-Roch theorem we can obtain 
``fixed point formulas" which calculate the classes $\bar\chi^P(\F)$ 
via localization on the ramification locus. This is carried out
in joint work with T. Chinburg and M. Taylor [CPT2]. 
Let us mention here that in the present paper we deal with 
the Galois modules given by the cohomology of general equivariant sheaves. 
These are significantly more difficult to determine than the virtual Galois module
of the Euler characteristic of the de Rham complex 
which was considered
in [CEPT] and [CPT1].
In these references, an important part of the calculation 
was quickly reduced to the case $d=0$ which was 
then treated using the work of Fr\"ohlich and Taylor. Here 
we are discussing a genuinely higher dimensional problem. In addition, we 
cannot in general expect a relation with the $\ep$-constants of Artin-Hasse-Weil 
$\rm L$-functions as in loc. cit. (Although, for curves, one can 
anticipate a relation of $\bar\chi^P(\O_X)$
with the leading terms of Artin-Hasse-Weil $\rm L$-functions at $s=1$.) 
On the negative side, an extension of Theorem \ref{main1} (b) to the case that $G$ has
a non-abelian Sylow appears to be outside the reach of 
our techniques; this then raises interesting and fundamental
questions concerning torsors of non-abelian groups. 
In particular our results lead us to formulate the following conjecture
for arbitrary finite groups:

\begin{conjecture}
For any $d\geq 1$, there exists an integer $Q(d)$
with the following property: For any finite group 
$G$ such that ${\rm gcd}(\#G, Q(d))=1$ and any $G$-torsor $X\to Y$ with $Y$ projective and flat
of relative dimension $d$,  $\bar\chi^P(\F)$ is trivial in ${\rm Cl}(\Z[G])$ for all
$G$-equivariant coherent sheaves $\F$ on $X$.
\end{conjecture}

Let us now describe briefly the contents of the various sections of the paper.
In \S \ref{torsors} we give some background on commutative Picard categories and torsors.
In \S \ref{hypercubic} we define the notion of an $n$-cubic structure on a line bundle
over a commutative group scheme.  In \S \ref{multiext} and \S \ref{taylor}
we show how we can analyze line bundles with hypercubic structures
using multiextensions; these can be thought of as giving derivatives
of the hypercubic structure. The ``Taylor expansion" of \S \ref{taylor} 
expresses line bundles with hypercubic structures in terms of
multiextensions. In \S \ref{sheaves} we recall the 
interpretation of multiextensions as extensions in the derived category of 
abelian sheaves  due to Grothendieck and deduce some corollaries.
In \S \ref{overZ} and \S \ref{reflsect} we analyze multiextensions of multiplicative group schemes
over $\Z$ and relate them to cyclotomic ideal class groups
and to the reflection homomorphisms; we deduce results on line bundles
with hypercubic structures over $\Spec(\Z[G])$. In \S \ref{hypercubedet} we explain the results of [Du]
on the hypercubic structure on the determinant of cohomology. 
In \S \ref{multidet} we 
explain how we can calculate the multiextensions derived from
the determinant of the cohomology and, in the case of cyclic 
groups of prime order, give the explicit description
of the invariant $t_{d+1}(X/Y)$ that was referred to above. The proofs of the results
stated in this introduction are completed in \S \ref{Galoismod}. 
Finally, in the Appendix we show how an argument 
due to Godeaux and Serre combined with
an ``integral" version of Bertini's theorem
(based on a theorem of Rumely on the existence 
of integral points) allows us to construct ``geometric" $G$-torsors $\pi: X\to Y$
with $Y$ regular and $Y\to \Spec(\Z)$ projective and flat of relative dimension $d$
for any finite group $G$ and any integer $d\geq 1$. This is a result 
of independent interest.
\smallskip

The author would like to thank T. Chinburg, M. Taylor, R. Kulkarni and L. Washington
for   useful conversations and comments,
and F. Ducrot for making available
a preliminary version of the preprint [Du].

\bigskip
\bigskip

\bigskip

\section{Torsors}\label{torsors}

\subsection{} \label{picard}

Recall ([SGA4] XVIII, 1.4) that a (commutative) {\sl Picard category} is a non-empty category
$\P$
in which all morphisms are isomorphisms and which is equipped with an ``addition" functor
\begin{equation*}
+ \ : \P\times \P\to \P,\ \ (p_1,p_2)\mapsto p_1+ p_2,
\end{equation*}
 associativity isomorphisms
\begin{equation*}
\sigma_{p_1,p_2,p_3}: (p_1+ p_2)+ p_3\xrightarrow{\sim} p_1+ (p_2+ p_3),
\end{equation*}
functorial in $p_1$, $p_2$, $p_3$, and commutativity isomorphisms
\begin{equation*}
\tau_{p_1,p_2}: p_1+p_1\xrightarrow{\sim}p_2+p_1
\end{equation*}
functorial in $p_1$, $p_2$ which satisfy the following axioms:
\smallskip

(1) For every object $p$ of $\P$ the functor $\P\to \P$; $q\mapsto p+q$ is an equivalence of categories,

(2) $\tau_{p_2,p_1}\cdot \tau_{p_1,p_2}={\rm Id}_{p_1+p_2}$ for all objects $p_1$, $p_2$ of $\P$,

(3) (Pentagon Axiom) 
\begin{multline*}
({\rm Id}_{p_1}+\sigma_{p_2,p_3,p_4})\cdot \sigma_{p_1,p_2+p_3,p_4}\cdot (\sigma_{p_1,p_2,p_3}+{\rm Id}_{p_4})=
\sigma_{p_1, p_2, p_3+p_4}\cdot \sigma_{p_1+p_2,p_3,p_4},
\end{multline*}
for all objects $p_1$, $p_2$, $p_3$, $p_4$ of $\P$,
\smallskip

(4) (Hexagon Axiom) 
\begin{equation*}
\sigma_{p_1,p_2, p_3}\cdot \tau_{p_3,p_1+p_2}\cdot \sigma_{p_3,p_1,p_2}=({\rm Id}_{p_1}+\tau_{p_3,p_2})\cdot
\sigma_{p_1,p_3,p_2}\cdot (\tau_{p_3,p_1}+{\rm Id}_{p_2})
\end{equation*}
for all objects $p_1$, $p_2$, $p_3$ of $\P$.

A Picard category is always equipped with an ``identity object";
this is a pair $(\uno,\epsilon)$ of an object $\uno$ with an isomorphism
$\epsilon: \uno+\uno\xrightarrow{\sim} \uno$ which is unique up to unique isomorphism.
For every object $p$ of $\P$ there are isomorphisms $p+\uno\xrightarrow {\sim}p\xleftarrow{\sim} \uno+p$.
Fix an identity object $(\uno,\epsilon)$ for $\P$. 
Then for every object $p$ of $\P$ there is an ``inverse", i.e a pair $(-p, i_p)$
of an object $-p$ of $\P$ with a ``contraction" isomorphism $i_p: p+(-p)\xrightarrow{\sim} \uno$ which
is unique up to unique isomorphism (see [SGA4] XVIII for more details). We associate to 
$\P$ two commutative groups: The group $\pi_0(\P)$ of isomorphism classes of objects of $\P$ and 
the group $\pi_1(\P)={\rm Aut}_{\P}(\uno)$.
For every object $q$ of $\P$ the translation functor provides a canonical isomorphism between $\pi_1(\P)$
and ${\rm Aut}_{\P}(q)$. The symmetry isomorphism $\tau_{p,p}$ of $p+p$ defines an element 
$\varepsilon(p)$ of $\pi_1(\P)$ such that $\varepsilon(p)^2=1$; this gives a group 
homomorphism $\varepsilon: \pi_0(\P)\to \pi_1(\P)$.
If in addition to the above axioms we have
\smallskip

(s.c.) \ \ $\tau_{p, p}={\rm Id}_{p+p}$, for all objects $p$ of $\P$ \ (i.e if $\varepsilon= 1$),
\smallskip

\noindent then we will say that the Picard category is ``strictly commutative".
\smallskip

Note that a commutative group defines a ``discrete" s.c. Picard category:
The objects are the elements of the group, the only morphisms are the
identity morphisms and the addition is given by the group law.
A non-trivial example of a s.c. Picard category is provided by
the category ${\rm PIC}(S)$ of invertible $\O_S$-sheaves 
over a scheme $S$ with morphisms isomorphisms of $\O_S$-modules and with ``addition" given by tensor product 
([De]). We will also occasionally use the Picard category ${\rm PIC}_*(S)$ of ``graded"
invertible $\O_S$-sheaves on $S$: The objects here are pairs $(\L, d)$
of an invertible $\O_S$-sheaf together with a (Zariski) locally constant function $d: S\to \Z$.
There exist morphisms from $(\L,d)$ to $(\M, e)$ only if $d=e$ and in this case a
morphism is an $\O_S$-module isomorphism $\L\xrightarrow{\sim} \M$. The addition is defined by
\begin{equation*}
(\L, d)+(\M, e)=(\L\otimes_{\O_S}\M, d+e).
\end{equation*}
This is endowed with the usual associativity morphisms of the tensor product. The commutativity
morphism is given by
\begin{equation*}
\tau:  (\L\otimes_{\O_S}\M, d+e)\to  (\M\otimes_{\O_S}\L, e+d)\ ;\quad l\otimes m\mapsto (-1)^{de} m\otimes l,\ \ \ \ \ \ 
\end{equation*}
if $l$ and $m$ are  sections of $\L$ and $\M$. Notice that ${\rm PIC}_*(S)$ {\sl is not 
a strictly commutative} Picard category; we have $\ep((\L, d))=(-1)^{d^2}=(-1)^d$.
\medskip

We will  use  the following:

\begin{lemma}\label{construct}
([SGA4] VVIII 1.4.3) Let $(P_i)_{i\in I}$ be a family of objects 
of the s.c. Picard category $\PP$. Denote by $e: I\to \Z^{(I)}={\rm Maps}(I, \Z)$ 
the canonical map $e(i)(j)=\delta_{ij}$. Then there exists a map
$
\Sigma: \Z^{(I)}\to {\rm Ob}(\PP),
$
isomorphisms $a_i: \Sigma(e(i))\xrightarrow{\sim}P_i$
and $a_{\underline {n}, \underline {m}}: \Sigma(\underline {n}+\underline {m})\xrightarrow{\sim} \Sigma(\underline {n})+
\Sigma(\underline {m})$, $\underline {n}$, $\underline {m}\in \Z^{(I)}$,
such that the diagrams 
\begin{equation*}
 \xymatrix
{\Sigma(\underline {n}+\underline {m}+\underline {k})\ar@{=}[d]\ar[r]^-{a_{\underline {n},\underline {m}+\underline {k}}}&
\Sigma(\underline {n})+\Sigma(\underline {m}+\underline {k})\ar[r]^-{\text{\rm id}+a_{\underline {m}, \underline {k}}}&
\Sigma(\underline {n})+(\Sigma(\underline {m})+\Sigma(\underline {k}))\\
\Sigma(\underline {n}+\underline {m}+\underline {k})\ar[r]^-{a_{\underline {n}+\underline {m},\underline {k}}}&
\Sigma(\underline {n}+\underline {m})+\Sigma(\underline {k})\ar[r]^-{a_{\underline {n}, \underline {n}}+\text{\rm id}}
&(\Sigma(\underline {n})+\Sigma(\underline {m}))+\Sigma(\underline {k})\ar[u]^{\sigma}\\
}
\end{equation*}
and,
\begin{equation*}
\begin{CD}
\Sigma(\underline {n}+\underline {m})@>a_{\underline {n},\underline {m}}>>\Sigma(\underline {n})+\Sigma(\underline {m})\\
@| @V\tau VV\\
\Sigma(\underline {m}+\underline {n})@>a_{\underline {m},\underline {n}}>>\Sigma(\underline {m})+\Sigma(\underline {n})
\end{CD}
\end{equation*}
are commutative. The system $(\Sigma, (a_i), (a_{\underline {n}, \underline {m}}))$
is unique up to unique isomorphism and is functorial in $(P_i)_{i\in I}$.\endproof
\end{lemma}

If $\P$ and $\QQ$ are two Picard categories an {\sl additive} functor $F: \P\to \QQ$ 
is a functor equipped with isomorphisms
\begin{equation*}
f_{p_1,p_2}: F(p_1)+F(p_2)\xrightarrow {\sim} F(p_1+p_2)
\end{equation*}
for all objects $p_1$, $p_2$ of $\P$ which are functorial in $p_1$, $p_2$ and 
which are compatible with the associativity and commutativity isomorphisms
of $\P$ and $\QQ$ in the sense that we have
\begin{multline*}
\ \ f_{p_1,p_2+p_3}\cdot ({\rm Id}_{F(p_1)}+f_{p_2,p_3})\cdot \sigma_{F(p_1), F(p_2), F(p_3)}=\\=
F(\sigma_{p_1,p_2,p_3})\cdot f_{p_1+p_2,p_3}\cdot (f_{p_1,p_2}+{\rm Id}_{F(p_3)}), \ \ 
\end{multline*}
\begin{equation*}
f_{p_2,p_1}\cdot \tau_{F(p_1), F(p_2)}=F(\tau_{p_1,p_2})\cdot f_{p_1,p_2},
\end{equation*}
for all objects $p_1$, $p_2$, $p_3$ of $\P$.

We refer the reader to [SGA4] XVIII 1.4
for the definitions of a (s.c.) Picard stack
over a site $\mathfrak S$ and of an additive ($\Sfr$-)functor between two
Picard $\mathfrak S$-stacks (these are modeled on the definitions above). 
There is a construction corresponding to Lemma \ref{construct} 
in this context.
In our applications, the site $\mathfrak S$ will be given by a subcategory 
of the category of $S$-schemes $(\text{Sch}/S)$ equipped with a Grothendieck topology TOP which is either
the fpqc or the fppf topology. In this case, we may think of a Picard $\mathfrak S$-stack
as a fibered category  such that the fiber over each $S'\to S$
is a Picard category, with additive pull-back functors and which satisfies
a certain descent condition for both the objects and the morphisms.
As before, a commutative group scheme over $S$ gives a ``discrete" fpqc or fppf s.c. Picard 
stack.  Very often we will work 
with the  site $S_{\rm fppf}$ of $S$-schemes which are locally of finite presentation
with the fppf topology.  
If $\phi: T\to S$ is an $S$-scheme we will denote by 
${\mathcal {PIC}}(T)$ (resp. ${\mathcal {PIC}}_*(T)$) the Picard $S_{\rm fppf}$-stack 
given by (resp. graded) invertible $\O_{T\times_SS'}$-sheaves on the schemes
$T\times_SS'$, $S'\to S$ a  morphism in $S_{\rm fppf}$. 

\subsection{} \label{torsors1} Let $H\to S$ be a group scheme {\sl flat and affine over $S$}.
A $H$-torsor (or ``a torsor for $H$")
is a sheaf $T$ for the fpqc topology on $(\text{Sch}/S)$ 
with morphisms 
\begin{equation*}
m: T\times_S H\to T\ ;\quad m(t,h)=t\cdot h\  \end{equation*}
and $p: T\to S$ such that
\smallskip

a) $H$ operates on $T$ (via $m$) over $S$.
\smallskip

b) The map $T\times_S H\to T\times_S T $ given by $(t, h)\mapsto (m(t,h), t)$
is an isomorphism.
\smallskip

c) $p$ is a sheaf epimorphism.
\smallskip

By [DG] III \S 4 Prop. 1.9, under our assumption on $H$,
every $H$-torsor $T$ is representable
by an $S$-scheme which we will also denote by $T$; 
furthermore, the resulting morphism of schemes $p: T\to S$ is affine and faithfully 
flat and identifies $S$ 
with the (categorical) quotient $T/H$.
In what follows,
we will think of torsors as either 
fpqc sheaves or schemes without making
the distinction.

For more details on the following the 
reader can refer to [DG] III \S 4. 
A morphism between two torsors $T_1\to S$, $T_2\to S$,
is an $S$-morphism $f:T_1\to T_2$ which commutes 
with the $H$-action; by descent such a morphism
is necessarily an isomorphism.
If $Y$ is an $S$-scheme we will occasionally
use the expression ``$X\to Y$ is an $H$-torsor" to 
mean that $X\to Y$ is a torsor for the group scheme 
$H_Y:=H\times_SY\to Y$. By the above,  
if $\pi: X\to Y$ is an $H$-torsor,
then $\pi$ is affine and flat
and identifies $Y$ 
with the (categorical) quotient $X/H$.
In fact, this quotient is {\sl universal} in the 
sense that for every base change $S'\to S$,
the natural morphism
\begin{equation}\label{basechange}
(X\times_SS')/H\to (X/H)\times_SS'
\end{equation}
is an isomorphism. If $H\to S$ is in addition
of finite presentation then so is $\pi: X\to Y$.

Assume now in addition that $H\to S$ is commutative; let
$T_1\to  S$, $T_2\to S$ be two $H$-torsors. We let the group scheme
$H$ act on the fiber product $T_1\times_S T_2$ by 
$(t_1,t_2)\cdot h=(t_1\cdot h, t_2\cdot h^{-1})$.
The quotient $(T_1\times_S T_2)/H$
then gives an $H$-torsor over $S$ (the action is
via $(t_1,t_2)\cdot h=(t_1\cdot h, t_2)$)
which we will denote by $T_1\cdot T_2$. 
We can see that there are canonical isomorphisms of $H$-torsors
\begin{equation}\label{asscom}
T_1\cdot (T_2\cdot T_3)\simeq(T_1\cdot T_2)\cdot T_3, \quad T_1\cdot T_2\simeq T_2\cdot T_1.
\end{equation}
The above satisfy the appropriate axioms so that
we obtain on the category of $H$-torsors
the structure of a s.c. Picard
category.
The identity object is given by the {\sl trivial torsor}. This is $H\to S$ with 
action given by multiplication; there are canonical isomorphisms
\begin{equation}\label{torid}
H\cdot T\simeq T\cdot H\simeq T.
\end{equation}
The {\sl inverse} $T^{-1}\to S$ of an $H$-torsor $T\to S$
is the same scheme with new action $m^{-1}(t,h):=m(t,h^{-1})$;
there are canonical isomorphisms
\begin{equation}\label{torinv}
T^{-1}\cdot T\simeq T\cdot T^{-1}\simeq H.
\end{equation}
 
Denote by $\Gm$ the multiplicative group scheme over $S$. 
There is a natural additive equivalence between the category ${\rm PIC}(T)$ of invertible $\O_T$-sheaves over an $S$-scheme
$T$ and the category of $\Gm_T$-torsors
on $T$ given by 
\begin{equation*}
\L\to \underline{\rm Isom}_{\O_T}(\O_T,\L).
\end{equation*}
(In what follows, for simplicity, we will denote
the $\Gm_T$-torsor associated to the invertible sheaf $\L$ again by $\L$).

\subsection{}\label{torsors2}
Assume now that $G$ is a finite group; for a scheme $S$ we will
denote by $G_S$ the constant group scheme 
$\sqcup_{g\in G} S$ given by $G$. When $S=\Spec(\Z)$ 
we will often abuse notation and simply write 
$G$ instead of $G_S$. Let $T$ be an $S$-scheme 
with a right $G$-action (this is the same as a right
$G_S$-action). We will say that {\sl $G$ acts freely on $T$}
if for every $S$-scheme $S'$ the action of $G$ on the
set of $T(S')$ of $S'$-points of the $S$-scheme $T$ is free
(i.e all the stabilizers are trivial).

Assume now that $T\to S$ is quasi-projective. Then the quotient $T/G$ exists as a scheme;
when in addition $G$ acts freely on $T$ then $\pi: T\to T/G$ is a $G$-torsor
and the morphism $\pi$ is finite \'etale ([D-G] III, \S 2, ${\rm n}^{\rm o}$ 6). 
We will continue with these assumptions in the rest of this section.

A coherent (resp. locally free coherent) sheaf of $\O_T$-$G$--modules $\F$ on $T$ is a coherent (resp. locally free
coherent) sheaf of $\O_T$-modules with an action of $G$
compatible with the action of $G$ on $(T, \O_T)$ in the 
appropriate sense. We will often use the term $G$-equivariant coherent 
(resp. $G$-equivariant locally free coherent) sheaf on $T$ instead of coherent 
(resp. locally free coherent) $\O_T$-$G$--sheaf.  
Let $\F$ be a $G$-equivariant coherent sheaf on $T$
and suppose that $V=\Spec(C)$ is an open affine subscheme of $T/G$.
Since $\pi$ is finite,  $U=\pi^{-1}(V)$ is an
affine $G$-equivariant open subscheme of $T$.
The sections $\F(U)$ form a left $G$-module and since $\F(U)$ is also a $C$-module,
we obtain on it the structure of a left $C[G]$-module. Hence, $\F$ provides us with
a coherent sheaf of $\O_{T/G}[G]$-modules on $T/G$, which we will
denote by $\pi_*(\F)$. On the other hand, if $\Gg$ is a coherent
(locally free coherent) sheaf of $\O_{T/G}$-modules on $T/G$,
the pull back $\pi^*(\Gg)$ is a $G$-equivariant coherent 
(resp. locally free coherent) sheaf on $T$. The pull-back functor $\pi^*$ 
gives an equivalence between the category of coherent (resp. locally free coherent) sheaves on the quotient $T/G$
and the category of $G$-equivariant coherent 
(resp. locally free coherent) sheaves on $T$.
Its inverse functor is obtained by taking invariants
under $G$ of the direct image $\pi_*$. In particular, there
are adjunction functor isomorphisms
\begin{equation}\label{adjunction}
{\rm id}\simeq \pi^*\cdot (\pi_*)^G,\quad {\rm id}\simeq (\pi_*)^G\cdot \pi^*.
\end{equation}

\subsection{}\label{2d} Now suppose that in addition $G$ is commutative. Denote by $G^D_S$
the Cartier dual group scheme of $G_S$; by definition, this represents 
the sheaf of characters $\underline{\Hom}(G_S,\Gm_S)$. Let $\chi: G\to \Gamma(S',\O^*_{S'})=\Gm(S')$
be a character with $S'\to S$ an $S$-scheme; then
$\chi$ corresponds to an $S'$-point of $G^D_S$.  If $\F$ is a $G$-equivariant coherent 
(resp. locally free coherent) sheaf on $T$, consider the  coherent
sheaf $\F_{S'}:=\F\otimes_{\O_S}\O_{S'}$ on $T':=T\times_SS'$ obtained by pulling back $\F$ 
via ${\rm pr}_1: T'=T\times_SS'\to T$. We can use the character $\chi$ to define the 
structure  of an $\O_{T'}$-$G$-module $\F_{S'}(\chi)$ on $\F_{S'}$
by setting:
\begin{equation}
g\cdot (m\otimes a')=g\cdot m\otimes \chi(g)^{-1}a', 
\end{equation}
where $m$, $a'$ stand for sections of $\F$, resp. $\O_{S'}$. Taking the $G$-invariants
$(\pi'_*(\F_{S'}(\chi)))^G$
of the direct image by $\pi':T'=T\times_SS'\to (T\times_SS')/G=(T/G)\times_SS'$
gives a coherent sheaf on $(T/G)\times_SS'$ which we will denote by $\F_{\chi}$.
Now apply the above to the structure sheaf $\O_T$. The resulting $\chi\mapsto \O_{T,\chi}$ 
provides us with a functor from
the discrete Picard category $G^D(S')$ to the s.c.
Picard category ${\rm PIC}(S')$.

Recall that, by using the $G$-action on $\O_T$, we have given  on $\pi_*(\O_T)$ 
the structure of a coherent sheaf of $\O_S[G]$-modules on $S$; 
we may think of $\pi_*(\O_T)$  as a coherent 
$\O_{G^D_S}$-sheaf on $G^D_S$. 
Suppose now that $\chi_0: G\to \Gm(G^D_S)$ is the ``universal" $G^D_S$-valued character
obtained from the natural pairing $G^D_S\times_S G_S\to \Gm_S$.
We have 
\begin{equation*}
\O_{T,\chi_0}\simeq (\pi_*(\O_T)\otimes_{\O_S}\O_{G^D_S})^G=(\pi_*(\O_T)[G])^G 
\end{equation*}
with the $G$-action defined by $g\cdot ( \sum_h a_hh)=\sum_h g(a_h)g^{-1}h$.
We can see that this gives an isomorphism
\begin{equation}\label{27}
\O_{T,\chi_0}\simeq \pi_*(\O_T)
\end{equation}
of (invertible) $\O_{G^D_S}$-sheaves on $G^D_S$. 

Now if $\chi_1$, $\chi_2$ are $S'$-valued characters
of $G$ as above, the multiplication morphism
$\O_{T'}\otimes_{\O_{S'}}\O_T\to \O_{T'}$ 
induces a $G$-equivariant isomorphism
\begin{equation}
\O_{T'}(\chi_1)\otimes_{\O_{S'}}\O_{T'}(\chi_2)\xrightarrow{\sim} \O_{T'}(\chi_1\chi_2).
\end{equation}
It follows from (\ref{adjunction}) that, by passing to $G$-invariants, 
this gives an isomorphism of invertible $\O_{S'}$-sheaves
\begin{equation}\label{torsorIso}
c_{\chi_1,\chi_2}: \O_{T,\chi_1} \otimes_{\O_{S'}}\O_{T,\chi_2}\xrightarrow{\sim}\O_{T,\chi_1\chi_2} .
\end{equation}
We can see that the diagrams
\begin{equation}\label{assocStr}
\begin{CD}
\O_{T,\chi_1}\otimes_{\O_{S'}}\O_{T,\chi_1}\otimes_{\O_{S'}} \O_{T,\chi_1} @>c_{\chi_1,\chi_2}\cdot{\rm id}>> \O_{T, \chi_1\chi_2}\otimes_{\O_{S'}} \O_{T,\chi_3}\\
@V{\rm id}\cdot c_{\chi_2,\chi_3}VV @Vc_{\chi_1\chi_2,\chi_3}VV\\
\O_{T,\chi_1}\otimes_{\O_{S'}} \O_{T,\chi_2\chi_3}@>c_{\chi_1,\chi_2\chi_3}>> \O_{T,\chi_1\chi_2\chi_3}\\
\end{CD}
\end{equation}
and
\begin{equation}\label{commuStr}
\begin{CD}
\O_{T,\chi_1}\otimes_{\O_{S'}} \O_{T,\chi_2}  @>c_{\chi_1,\chi_2} >> \O_{T,\chi_1\chi_2}\\
@V\wr VV @|\\
 \O_{T,\chi_2}\otimes_{\O_{S'}} \O_{T,\chi_1}@>c_{\chi_2,\chi_1}>> \O_{T,\chi_2\chi_1}\\
\end{CD}
\end{equation}
are commutative.  (The left vertical arrow in (\ref{commuStr}) is the commutativity isomorphism;
in (\ref{assocStr}), for simplicity, we supress in the notation the associativity isomorphism
for the tensor product).
Therefore the isomorphisms $c_{\chi_1,\chi_2}$ are compatible with the 
associativity and commutativity
constraints of the source and target commutative Picard
categories. Hence, they equip the functor $\chi\mapsto \O_{T,\chi}$
with additivity data. Notice also that for any base change $q: S''\to S'$
there are natural isomorphisms
\begin{equation}
q^*\O_{T,\chi} \simeq \O_{T,\chi\cdot q}.
\end{equation}
(on the right hand side, we view $\chi$ as a 
morphism $S'\to G^D_S$). These satisfy the usual cocycle condition
and they are compatible with the additivity isomorphisms (\ref{torsorIso}).
Hence, we can view $\chi\mapsto \O_{T,\chi}$ as providing
an additive functor
\begin{equation*}
G^D_S\to {\mathcal {PIC}}(S)
\end{equation*}
 from the  Picard $S_{\rm fppf}$-stack
given by $G^D_S$ to the  Picard   $S_{\rm fppf}$-stack 
of invertible sheaves.
\smallskip

\begin{Remark}\label{rem22}
{\rm As we shall see in \S \ref{extension}, the above allows us to construct a commutative extension 
\begin{equation}\label{groth}
1\to \Gm_S\to E\to G^D_S\to 0
\end{equation}
such that the fiber of $E$ over the $S'$-point of $G^D_S$ which is given by the character $\chi$
is isomorphic to the $\Gm_{S'}$-torsor that corresponds to the invertible $\O_{S'}$-sheaf $\O_{T,\chi}$.
We will see in \S \ref{sheaves} that this construction  
gives an isomorphism between the group ${\rm H}^1(S,G)$ of isomorphism classes of $G$-torsors
over $S$ and the group of isomorphism classes of commutative extensions of $G^D_S$ by $\Gm_{S}$.}
\end{Remark}

\begin{Remark}\label{rem23}
{\rm Suppose that $G=\Z/n$ and that $S$ is a Noetherian $\Spec(\Z[\zeta_n,1/n])$-scheme
with $\zeta_n$ a primitive $n$-th root of unity. 

a)  Let us consider $\psi_0: G\to \O_S^*$  defined by $\psi_0(a)=\zeta_n^a$. In this case, the association
$T\mapsto \O_{T,\psi_0}$ gives an equivalence between the category 
of $\Z/n$-torsors $T\to S$ and the category of pairs $(\L, t)$ of invertible sheaves $\L$ over $S$
with an isomorphism $t: \L^{\otimes n}\xrightarrow {\sim}\O_S$. (For $\L=\O_{T,\psi_0}$, $t$ is given by 
$\O_{T,\psi_0}^{\otimes n}\simeq \O_{T,\psi_0^n}=\O_{T,1}\simeq \O_S$, where the first isomorphism
is obtained using (\ref{torsorIso})).
If $S$ is 
semilocal, we can choose a trivialization $t': \O_S\xrightarrow{\sim} \L$
and consider $f=t(t'(1)^{\otimes n})\in \O_S^*$. Then we can see that
the isomorphism class of the pair $(\L, t)$ is determined by the ``Kummer element"
$f\in \O_S^*/(\O_S^*)^n$; this gives the isomorphism of Kummer theory
\begin{equation}\label{kummer}
{\rm H}^1(S,\Z/n)\xrightarrow {\sim} \O_S^*/(\O_S^*)^n.
\end{equation}
Note that if we replace in the construction $\psi_0$ by the inverse character $\psi_0^{-1}$, 
we obtain (\ref{kummer}) composed with the homomorphism $f\mapsto f^{-1}$.

b) Suppose that $\phi: S'\to S$ is a finite locally free morphism.   
If $\L$ is an invertible sheaf over $S'$, then we can define its norm by
$$
{\rm Norm}_{S'/S}(\L)={\rm det}(\phi_*(\L))\otimes (\det(\phi_*(\O_{S'}))^{-1}
$$
(here ${\rm det}$ simply means highest exterior power). The norm gives an additive functor from
the Picard category of invertible sheaves on $S'$ to the Picard 
category of invertible sheaves on $S$.
For a finite locally free morphism $\phi: S'\to S$
and a $\Z/n$-torsor $T'\to S'$ we can consider the $\Z/n$-torsor
over $S$ which corresponds to the pair $({\rm Norm}_{S'/S}(\L), {\rm Norm}_{S'/S}(t))$. 
By the construction of [SGA4] XVIII this corresponds to the ``trace" map
\begin{equation}
{\rm Tr}_{\phi}: {\rm H}^1(S',\Z/n)\to {\rm H}^1(S,\Z/n)\ .
\end{equation}
If $S$ is semilocal then so  is $S'$. The above then shows that 
the diagram
\begin{equation}\label{kummerNorm}
\begin{CD}
{\rm H}^1(S',\Z/n)@> {\sim}>> \O_{S'}^*/(\O_{S'}^*)^n\\
@V{\rm Tr}_{\phi}VV @VV{\rm Norm}_{S'/S} V\\
{\rm H}^1(S,\Z/n)@> {\sim}>> \O_S^*/(\O_S^*)^n
\end{CD}
\end{equation}
commutes.
}
\end{Remark}

\bigskip

\section{Hypercubic structures}\label{hypercubic}
\setcounter{equation}{0}

\subsection{} Let $H\to S$ be a commutative $S$-group scheme. 
For $n\geq 1$, we will denote by $H^n:=H\times\cdots \times H$ the $n$-fold 
fiber product over $S$ (for simplicity, we will often omit the subscript $S$ in the 
notation of the product). 
If $I$ is a subset of the index set $\{1,\ldots, n\}$,
we will denote by $m_{I}$ the morphism $H^n\to H$ given 
on points by $(h_1,\ldots, h_n)\mapsto \sum_{i\in I}h_i$ (if $I=\emptyset$, $m_{I}(h_1,\ldots, h_n)=0$).
 When $I=\{i\}$, then $m_{I}$ is the $i$-th projection
$p_{i,H}: H^n\to H$. 
Recall that we identify the s.c.
Picard category  of invertible sheaves over an $S$-scheme
$T$ with  the s.c. Picard category of $\Gm_T$-torsors
on $T$ (see \ref{torsors1}).
Lemma \ref{construct} implies that the ``tensor operations" we use in what follows to define
invertible sheaves or $\Gm$-torsors give results that are well-defined 
up to coherent canonical isomorphism. 

 If $\L$ is an invertible
sheaf on $H$, then we set 
\begin{equation}\label{theta}
\Theta_{n}(\L)=\bigotimes_{I\subset\{1,\ldots, n\}}m_{I}^*(\L)^{(-1)^{n-\#I}}
\end{equation}
(an invertible sheaf on $H^n$). A permutation $\sigma:\{1,\ldots, n\}\to \{1,\ldots, n\}$
induces a corresponding $S$-isomorphism $\sigma: H^n\to H^n$. Since
$ m_{I}\cdot \sigma=m_{\sigma(I)}$ permuting the factors of \ref{theta} gives 
a canonical isomorphism 
\begin{equation}\label{permu} {\mathfrak P}_\sigma: \sigma^*\Theta_{n}(\L)\xrightarrow{\sim} \Theta_{n}(\L)\ .
\end{equation}
Now suppose that $n\geq 2$
and consider the morphisms
$A,B,C,D: H^{n+1}\to H^n$ given by 
\begin{eqnarray}
A(h_0,h_1,h_2,\ldots, h_n)&=& (h_0+h_1,h_2,\ldots, h_n),\\
B(h_0,h_1,h_2,\ldots, h_n)&=&(h_0,h_1,h_3, \ldots, h_n),\\
C(h_0,h_1,h_2,\ldots, h_n)&=&(h_0,h_1+h_2,h_3,\ldots, h_n),\\
D(h_0,h_1,h_2,\ldots, h_n)&=&(h_1,h_2,h_3,\ldots, h_n).
\end{eqnarray}
We can observe that there is a canonical isomorphism
\begin{equation}\label{Q}
{\mathfrak Q}: A^*\Theta_{n}(\L)\otimes B^*\Theta_{n}(\L)\xrightarrow{\sim}  C^*\Theta_{n}(\L)\otimes D^*\Theta_{n}(\L)
\end{equation}
which is obtained by contracting duals and permuting factors
(cf. [Br] \S 2 or [AHS] \S 2). (The order in which these operations are performed 
in the s.c. Picard category is of no consequence;
the isomorphism remains the same. This can be viewed as a consequence
of Lemma \ref{construct}.)

Finally observe that if $(0,\ldots,0): S\to H^n$ is the zero section,
there is a canonical isomorphism
\begin{equation}\label{zero}
{\mathfrak R}: (0,\ldots,0)^*\Theta_n(\L)\xrightarrow{\sim} \O_S\ .
\end{equation}
\begin{Definition} \label{cubicdef}
Let $n\geq 2$. An $n$-cubic structure on the invertible sheaf $\L$ over $H$ is 
an isomorphism of invertible sheaves on $H^n$
\begin{equation}
\xi: \O_{H^n}\xrightarrow{\sim}\Theta_{n}(\L)
\end{equation}
(i.e a choice of a global generator $\xi(1)$ of
$\Theta_{n}(\L)$)
which satisfies the following conditions:
\smallskip

c0) It is ``rigid", i.e if $(0,\ldots,0): S\to H^n$
is the zero section, then
\begin{equation*}
{\mathfrak R}((0,\ldots ,0)^*(\xi(1)))=1\ .
\end{equation*}

c1) It is ``symmetric", i.e 
for all $\sigma\in S_n$,
\begin{equation*}
{\mathfrak P}_\sigma(\sigma^*(\xi(1)))=\xi(1)\ .
\end{equation*}

c2) It satisfies the ``cocycle condition"
\begin{equation*}
{\mathfrak Q}(A^*(\xi(1))\otimes B^*(\xi(1)))=C^*(\xi(1))\otimes D^*(\xi(1))\ .
\end{equation*}
\end{Definition}

\begin{Remark}\label{32re}
{\rm  a) This definition also appears in [AHS] 2.39. 
For $n=3$ it is a slight variant of Breen's definition 
of a cubic structure ([Br] \S 2; see also Moret-Bailly [MB]). To explain this set 
\begin{equation*}
\Theta (\L):=\bigotimes_{\emptyset\neq I\subset \{1,2,3\}}m_I^*(\L)^{(-1)^{3-\#I}}=\Theta_3(\L)\otimes m_\emptyset^*\L.
\end{equation*}
(Recall that $m_\emptyset: H^3\to H$ is the zero homomorphism).
There are isomorphisms analogous to (\ref{permu}) and (\ref{Q}) for $\Theta(\L)$. 
According to Breen, a cubic structure on $\L$ is a trivialization $t:\O_{H^3}\xrightarrow{\sim}\Theta (\L)$
which respects these isomorphisms (i.e satisfies conditions analogous to (c1) and (c2)). On the other hand,  
(\ref{zero}) induces a canonical isomorphism
\begin{equation*}
0^*\L\xrightarrow{\sim}(0,0,0)^*\Theta (\L).
\end{equation*}
Hence, $t$ also induces a ``rigidification" of $\L$, i.e an isomorphism $\O_S\xrightarrow{\sim} 0^*\L$
which we will denote by $r(t)$.
For any invertible sheaf $\L$ on $H$ now set $\L^{\rm rig}:=\L\otimes p^*0^*\L^{-1}$, $p: H\to S$
the structure morphism. The invertible sheaf $\L^{\rm rig}$ is equipped with a canonical rigidification
$r_{\rm can}$ and so there is a canonical isomorphism
\begin{equation*}
\phi_{\rm can}:  \Theta_3(\L)\xrightarrow{\sim}\Theta (\L^{\rm rig}).
\end{equation*}
One can now verify that $\xi\mapsto t(\xi):=\phi_{\rm can}\cdot \xi$ gives
a bijective correspondence between the set of $3$-cubic structures $\xi$ on $\L$ in the sense 
above and the set of Breen's cubic structures $t$ on $\L^{\rm rig}$ which satisfy $r(t)=r_{\rm can}$ 
(cf. [Br] \S 2.8 and [AHS] Remark 2.44).  
\smallskip

b) In what follows, we will often denote various invertible sheaves
by giving their fibers over a ``general" point of the base. For example,
we can denote $\Theta_3(\L)$ as
$$
\L_{x+y+z}\otimes\L_{x+y}^{-1}\otimes\L_{y+z}^{-1}\otimes\L_{z+x}^{-1}\otimes\L_x\otimes\L_y\otimes\L_z\otimes\L_{0}^{-1}.
$$
(This gives the fiber of $\Theta_3(\L)$ over the point $(x,y,z)$ of $H^3$.)}
\end{Remark}
 \subsection{} \label{various}
i) By definition, an isomorphism between the invertible sheaves
with $n$-cubic structures $(\L,\xi)$ and $(\L',\xi')$ is an isomorphism 
$\phi:  \L\xrightarrow{\sim}\L'$ such that
\begin{equation*}
\Theta_n(\phi)\cdot \xi=\xi'\ ,
\end{equation*}
where $\Theta_n(\phi): \Theta_n(\L)\xrightarrow{\sim}\Theta_n(\L')$ is functorially induced
from $\phi$. If $(\L,\xi)$, $(\L',\xi')$ are invertible sheaves
with $n$-cubic structures we define their product 
\begin{equation*}
(\L,\xi)\cdot(\L',\xi')=(\L\otimes\L',\xi*\xi')
\end{equation*}
where $\xi*\xi'$ is the composition
$$
\O_{H^n}=\O_{H^n}\otimes_{\O_{H^n}}\O_{H^n}\xrightarrow{\xi\otimes\xi'}\Theta_n(\L)\otimes_{\O_{H^n}}
\Theta_n(\L')\xrightarrow{\alpha}\Theta_n(\L\otimes\L')
$$
with $\alpha$ the standard natural isomorphism. 
We can see that the pairs $(\L,\xi)$ give the objects of a s.c.
Picard category $n{\hbox{\small -}}{\rm CUB}(H,\Gm)$ with arrows given by isomorphisms
as above and ``addition" given by the above product. (This is similar to the corresponding statement 
for Breen's cubic structures; see [Br] \S 2).
 
\medskip

ii) Suppose that  
the invertible sheaf $\L$ on $H$ is trivial
via $\psi: \O_H\xrightarrow{\sim}\L$. This then
induces a trivialization 
\begin{equation}\label{thetatriv}
\Theta_n(\psi): \O_{H^n}=\Theta_n(\O_H)\xrightarrow{\sim}\Theta_n(\L)\ .
\end{equation}
A second trivialization $\xi: \O_{H^n}\xrightarrow{\sim}\Theta_{n}(\L)$
can now be given via the ratio of the generators
\begin{equation}\label{elem}
c=\xi(1)/\Theta_n(\psi)(1)\in \Gamma(H^n, \O_{H^n}^*)\ .
\end{equation}
 In this case, we can see 
that $\xi: \O_{H^n}\xrightarrow{\sim}\Theta_{n}(\L)$
gives an $n$-cubic structure on $\L$ if and only if the element $c$ satisfies:
\smallskip

\noi c0)  $c(0,\ldots, 0)=1$,
\smallskip

\noi c1)  $c(h_{\sigma(1)},\ldots, h_{\sigma(n)})=c(h_{1},\ldots, h_{n})$, for all $\sigma\in S_n$, 
\smallskip

\noi c2) 
$
c(h_0+h_1,h_2,\ldots, h_n)c(h_0,h_1,h_3, \ldots, h_n)=
c(h_0,h_1+h_2,h_3,\ldots, h_n)c(h_1,h_2,\ldots, h_n).
$
\\

\noi (Here $h_i$, $0\leq i\leq n$, range over all $T$-valued points of $H$, $T$ any $S$-scheme.
Also, for example, $c(h_1,h_2,\ldots, h_n)\in \Gamma(T, \O^*_T)$ is obtained from $c$ by pulling back along 
the morphism $T\to H^n$ given by $(h_1,h_2,\ldots, h_n)$).

An inductive argument shows that if $c\in \Gamma(H^n, \O_{H^n}^*)$ satisfies (c0), (c1), (c2) above,  
then it also satisfies
\smallskip

c0${}^\prime$)\ \ \ \ \  $c(h_1,h_2,\ldots, h_n)=1$, if at least one of the $h_i$ 
is $0$. 

\medskip

iii) Suppose that $S=\Spec(R)$ and $H=G^D_S$, the Cartier dual
of a finite abelian {\sl constant} group scheme $G$. Then, $H=\Spec(R[G])$, $H^n=\Spec(R[G\times\cdots \times G])$.
If $T=\Spec(R')$, then $T$-valued points $h_i: T\to H$ correspond to $R'$-valued characters
$\chi_i: G\to {R'}^*$. Suppose now that $R$ is local; then $R[G]$ is semi-local
and any invertible sheaf $\L$ on $H=\Spec(R[G])$ is trivial. Hence, from (ii), we see that 
$n$-cubic structures on $\L$ are given by units $c\in R[G^n]^*$
which satisfy
\smallskip

c0)  \ \ \ $  (1\otimes\cdots \otimes 1)(c)=1$,
\smallskip

c1)  \ \ \ $(\chi_{\sigma(1)}\otimes\cdots\otimes\chi_{\sigma(n)})(c)=
(\chi_{1}\otimes\cdots\otimes\chi_{n})(c)$,
\begin{eqnarray*}
c2)\ \ \  \ \ (\chi_0\chi_1\otimes \chi_2\otimes\cdots\otimes \chi_n)(c) 
(\chi_0\otimes \chi_1\otimes\chi_3\otimes\cdots\otimes \chi_n)(c)&=&\ \ \ \ \ \ \ \ \ \ \ \ \ \ \ \ \ \ \ \ \ \\
\ \ \ \ \ (\chi_0\otimes\chi_1 \chi_2\otimes\chi_3\otimes\cdots\otimes \chi_n)(c)
(\chi_1\otimes \chi_2\otimes\cdots\otimes \chi_n)(c)\ .\  
\end{eqnarray*}

(In these relations $\chi_1\otimes\cdots \otimes\chi_n$ etc. are characters of $G^n$ which are evaluated
on the element $c$ of $R[G^n]$).
\smallskip

As above if $c\in R[G^n]^*$ satisfies (c0), (c1), (c2) above (for all characters of $G$),  
then it also satisfies
\smallskip

c0${}^\prime$)\ \ \ \ \  $(\chi_{1}\otimes\cdots\otimes\chi_{n})(c)=1$, if at least one of the characters $\chi_i$ 
is trivial.

\begin{Definition} An element $c$ of $R[G^n]$ which
satisfies (c0), (c1), (c2) above (for all characters of $G$) 
is called {\sl $n$-cubic}.
\end{Definition}

\subsection{} \label{augmenG} Suppose that $A$ is an abelian group and $n\geq 1$. Denote by $I[A]$ the 
augmentation ideal of the group ring $\Z[A]$; by definition, this is the kernel of the 
ring homomorphism $\Z[A]\to \Z$; $\sum_{a}n_a[a]\mapsto \sum_an_a$.
Set
\begin{equation}
C_n(A):={\rm Sym}^n_{\Z[A]}I[A]
\end{equation}
(the $n$-th symmetric power of the $\Z[A]$-module $I[A]$; cf. [AHS] 2.3.1
where the reader is also referred to for more details). The abelian group $C_n(A)$ is the quotient
of ${\rm Sym}^n_{\Z}I[A]$ by all the relations of the form
\begin{eqnarray*}
([b+a_1]-[b])\otimes([a_2]-[0])\otimes\cdots\otimes([a_n]-[0])=\ \ \ \ \ \ \ \ \ \\
\ \ \ \ \ \ \ \ \ \ \ \ \ \ \qquad =([a_1]-[0])\otimes([b+a_2]-[b])\otimes\cdots \otimes([a_n]-[0])\ 
\end{eqnarray*}
with $a_1,\ldots, a_n, b\in A$.
After rearranging and reindexing, this relation 
can be expressed in terms of the generators $[a_1,\ldots, a_n]:=([a_1]-[0])\otimes\cdots\otimes([a_n]-[0])$
of ${\rm Sym}^n_{\Z}I[A]$ as
\begin{eqnarray*}
[a_1,a_2,\ldots, a_n]+[a_0,a_1+a_2,a_3,\ldots, a_n]=\ \ \ \ \ \ \ \\
\qquad\quad=[a_0,a_1,a_3, \ldots, a_n]+[a_0+a_1,a_2,\ldots, a_n].
\end{eqnarray*}

Now suppose that $A=H(\Spec(R'))$, the  group of characters of
the finite abelian group $G$ with values in the $R$-algebra $R'$.
We can see by the above that an  
$n$-cubic element $c\in R[G^n]^*$ gives a
group homomorphism 
\begin{equation}
\alpha(c): C_n(A)\to {R'}^*\ ;\ \ [\chi_1,\ldots, \chi_n]\mapsto (\chi_{1}\otimes\cdots\otimes\chi_{n})(c)\ .
\end{equation}
In fact, every element $c\in R[G^n]^*$ which satisfies (c0$'$) gives
a group homomorphism 
$$
\ \ \ \ \bigotimes^n _{\Z}I[A]\to {R'}^*;\qquad ([\chi_1]-[1])\otimes\cdots\otimes([\chi_n]-[1])\mapsto (\chi_{1}\otimes\cdots\otimes\chi_{n})(c)\ .
$$
By the above, if this homomorphism factors through $C_n(A)$ (for all $R$-algebras $R'$) 
then $c$ is $n$-cubic.
\bigskip

\section{Multiextensions}\label{multiext}
\setcounter{equation}{0}

\subsection{} \label{extension} Suppose that $J$ and $H$ are two  
flat commutative group schemes
over the scheme $S$. We will assume that $J\to S$ is affine. 
(Most of the time we will take $J={\bf G}_m$ and $H=G^D_S$, the Cartier dual
of a finite constant abelian group scheme $G_S$.)  According to [SGA7I] 
Expos\'e VII \S 1 giving a commutative group scheme extension 
$E$ of $H$
by $J$ is equivalent to giving, for every $S$-scheme $U\to S$ and each 
$U$-point $a: U\to H$ over $S$, a $J_U$-torsor $E_a$ with the following additional structure: These torsors
should come together with isomorphisms
\begin{equation}\label{51}
c_{a,a'}:  E_{a} \cdot E_{a'}\xrightarrow{\sim}E_{a+a'}
\end{equation}
such that the diagrams
\begin{equation}\label{assoc}
\begin{CD}
E_a\cdot E_{a'}\cdot E_{a''} @>c_{a,a'}\cdot{\rm id}>> E_{a+a'}\cdot E_{a''}\\
@V{\rm id}\cdot c_{a',a''}VV @Vc_{a+a',a''}VV\\
E_a\cdot E_{a'+a''}@>c_{a,a'+a''}>> E_{a+a'+a''}\\
\end{CD}
\end{equation}
and
\begin{equation}\label{commu}
\begin{CD}
E_a\cdot E_{a'}  @>c_{a,a'} >> E_{a+a'} \\
@V\wr VV @|\\
 E_{a'}\cdot E_{a}@>c_{a',a}>> E_{a'+a}\\
\end{CD}
\end{equation}
are commutative. Both the torsors and the isomorphisms should be functorial in the $S$-scheme $U$.
(Of course, $E$ is then $E_{\rm id}$ for ${\rm id}: H\to H$ the ``universal" $H$-point of $H$.)
For simplicity, in the first diagram we write
$E_a\cdot E_{a'}\cdot E_{a''}$ instead of using the canonical associativity isomorphism 
$E_a\cdot (E_{a'}\cdot E_{a''})\simeq (E_a\cdot E_{a'})\cdot E_{a''}$.
In the second, the left vertical isomorphism is the commutativity isomorphism
(see (1.1.4.1) and (1.2.1) in loc. cit.) Another way to describe these
conditions is to say that the above data define an additive 
functor from the discrete   Picard category of the abelian 
group $H(U)$ to the Picard category of $J_U$-torsors.
As a result, when $H\to S$ is fppf,
we can view extensions of $H$ by $J$
as being given by additive functors from the discrete Picard  $S_{\rm fppf}$-stack
of $H$ to the Picard $S_{\rm fppf}$-stack
of  $J$-torsors (this point of view is explained in more detail
in [SGA4] XVIII). 

\subsection{} \label{nexteNew} We refer the reader to [SGA7I] Expos\'e VII
for the definition of $J$-biextensions of commutative group schemes.
There is an obvious generalization of both the notions
of extension and biextension: the notion of an $n$-extension of $(H,\ldots,H)$ 
by $J$. (Often, for simplicity, we will just say ``$n$-extension of $H$ by $J$";
for $n=1$ this gives a usual extension of commutative group
schemes as above and for $n=2$ a $J$-biextension of $(H,H)$).
By definition (see loc. cit. Def. 2.1 for $n=2$
and 2.10.2 in general) such an $n$-extension is a $J$-torsor $E$ over $H^n$
equipped with ``compatible partial composition laws". Giving an $n$-extension of $(H,\ldots, H)$ by $J$
is equivalent to giving, for
each $S$-scheme $U\to S$ and $U$-valued point $(a_1,\ldots, a_n)$ of $H^n$ over $S$,
 a $J$-torsor  
$E_{(a_1, \ldots, a_n)}$ over $U$ with additional structure: These torsors should come together with  
isomorphisms ($i=1,\ldots, n$)
\begin{equation}\label{isomulti}
c^i_{a_1,\ldots, a_i;a'_i,\ldots, a_n}:E_{(a_1,\ldots, a_i,\ldots, a_n)}\cdot E_{(a_1,\ldots, a'_i,\ldots, a_n)} \xrightarrow{\sim} E_{(a_1,\ldots, a_i+ a'_i,\ldots, a_n)}
\end{equation}
which  satisfy conditions as in (\ref{assoc}), (\ref{commu}). In addition, we require the
following  
compatibility requirement between the isomorphisms for various $i$:
For all pairs $i\neq j$, the diagram
\begin{equation}\label{compa}
\begin{CD}
E_{( a_i ,  a_j  )}\cdot E_{( a'_i ,  a_j  )}\cdot 
E_{( a_i , a_j' )}\cdot E_{( a_i',   a_j' )}@>{c^i\cdot c^i}>> 
E_{( a_i+a_i' , a_j )}\cdot E_{( a_i+a_i', a'_j  )}\\
@V{(c^j\cdot c^j)\circ \phi}VV @V{c^j}VV\\
E_{( a_i ,  a_j+a_j' )}\cdot E_{(  a_i', a_j+ a'_j  )}@>{c^i}>>
E_{( a_i+a_i',  a_j+a_j'  )}\\
\end{CD}
\end{equation}
is commutative. (Here, for simplicity, we write $E_{( b_i ,  b_j  )}$ instead of 
$E_{(b_1,\ldots, b_i,\ldots, b_j,\ldots, b_n)}$ and omit the subscripts from 
the notation of the composition laws. Also, we denote by $\phi$ the canonical
isomorphism $E_{( a_i ,  a_j  )}\cdot E_{( a'_i ,  a_j  )}\cdot 
E_{( a_i , a_j' )}\cdot E_{( a_i',   a_j' )}\simeq E_{( a_i ,  a_j  )}\cdot E_{( a_i ,  a_j'  )}\cdot 
E_{( a_i' , a_j )}\cdot E_{( a_i',   a_j' )}$.) Once again both the torsors and the isomorphisms should be 
functorial on the $S$-scheme $U$. 

In fact, these conditions on (\ref{isomulti}) are such that,
for each $i$, these isomorphisms equip $E$
with a structure of an extension of commutative group schemes over 
$H\times\cdots\times\underset{i}{\hat H}\times\cdots\times H$
(this notation means that the $i$-th factor is omitted from the product)
\begin{equation}\label{exti}
0\to J_{H\times\cdots\times\underset{i}{\hat H}\times\cdots\times H} \to E\to H_{H\times\cdots\times\underset{i}{\hat H}\times\cdots\times H}\to 0\ .
\end{equation}
There is an obvious notion of isomorphism
between $n$-extensions of $H$ by $J$
(it is given by an isomorphism of the corresponding torsors that
respects the composition laws (\ref{isomulti})). The $n$-extensions
of $H$ by $J$ form the objects of a strictly commutative 
Picard category $n{\hbox{\small -}}{\rm EXT}(H,J)$ with morphisms given by isomorphisms of $n$-extensions
and product corresponding to product of the extensions (\ref{exti}).
The identity object is given by the 
``trivial" $n$-extension on the trivial torsor $J\times_S H^n$.
These facts are explained in detail in [SGA7I] Exp. VII \S 1 and \S 2 when $n=1$, $2$. The same constructions apply 
to the general case (see loc. cit. Remark 3.6.7). We will denote 
by $n{\hbox{\small -}}{\rm Ext}^1(H,J)$ the commutative group  
of isomorphism classes of $n$-extensions of $H$ by $J$ and by $n{\hbox{\small -}}{\rm Ext}^0(H,J)$
the commutative group of the endomorphisms of the identity object.
 
Note that sending the class of an $n$-extension to 
the class of the underlying $J$-torsor 
over $H^n$
defines a group homomorphism
\begin{equation}\label{4.7}
t: n{\hbox{\small -}}{\rm Ext}^1(H,J)\to {\rm H}^1(H^n,J) .
\end{equation}
When $J=\Gm$, we can view this as a homomorphism
\begin{equation}\label{4.8}
t: n{\hbox{\small -}}{\rm Ext}^1(H,\Gm)\to {\rm Pic}(H^n) .
\end{equation}

\subsection{} \label{symmpar}Suppose that $E$ is an $n$-extension of $H$
by $J$.  If $\sigma\in S_n$ is a permutation, then we also denote by $\sigma: H^n\to H^n$ 
the corresponding automorphism. We can see that the pull-back $J$-torsor $\sigma^*E$ 
 also supports a canonical structure of an $n$-extension 
of $H$ by $J$. Denote by $\Delta_n$   the diagonal homomorphism $H\to H^n$.

We will say that the $n$-extension $E$ of $H$ by $J$ is {\sl symmetric} 
if it comes together with isomorphisms of $n$-extensions
\begin{equation*}
\Psi_\sigma: \sigma^*E\xrightarrow {\sim} E,\quad \text{for each\ }\sigma\in S_n,
\end{equation*}
which satisfy the following properties:

i) $\Delta_n^*\Psi_\sigma=i_\sigma$ where 
$i_\sigma: \Delta_n^*\sigma^*E\xrightarrow{\sim} \Delta_n^*E$ is the 
natural isomorphism of $J$-torsors obtained by $\sigma\cdot\Delta_n=\Delta_n$.

ii) For every pair $\sigma$, $\tau\in S_n$, the following diagram is commutative
\begin{equation*}
\begin{CD}
\sigma^*(\tau^* E) @>\sigma^*\Psi_\tau>> \sigma^*E\\
@V\wr VV @VV\Psi_\sigma V\\
(\tau\sigma)^*E @>\Psi_{\tau\sigma}>>E\\
\end{CD}
\end{equation*}
where the left vertical arrow is the natural isomorphism of $J$-torsors. 

Notice that the trivial $n$-extension is naturally symmetric.
When $n=1$ every extension is symmetric with $\Psi_{\rm id}={\rm id}$.

By definition, an isomorphism between two symmetric $n$-extensions
$(E,\{\Psi_\sigma\})$ and $(E',\{\Psi'_\sigma\})$
is  an isomorphism 
$
f:E\xrightarrow{\sim} E'
$
of $n$-extensions such that for any $\sigma\in S_n$
the following diagram commutes:
\begin{equation}\label{isosymm}
\begin{CD}
\sigma^*E@>\sigma^*f>>\sigma^*E'\\
@V\Psi_\sigma VV @VV\Psi'_\sigma V\\
E@>f>> E'.\\
\end{CD}
\end{equation}

\bigskip

\section{Differences and Taylor expansions}\label{taylor}
\setcounter{equation}{0}

In this section  we assume that $S=\Spec(R)$ and $H=G^D_S$, the Cartier dual
of the finite {\sl constant} group scheme given by the abelian group $G$. 
The constructions which we will describe 
below are certainly valid under less restrictive hypotheses. However,
we are only going to need them under these assumptions and so we choose to explain 
them only in this case since then the presentation simplifies considerably.
With the exception of the ``Taylor expansion" of \S \ref{taylorsect} 
most of them are relatively straightforward generalizations of similar constructions
described in [Br] for $n=3$. We  consider the ``$n-1$-th symmetric difference"
$\Theta_{n-1}(\L)$ of the invertible sheaf with an 
$n$-cubic structure $(\L,\xi)$ on $H$. We first show that $\Theta_{n-1}(\L)$ is naturally equipped with the structure
of an $n-1$-extension and then explain  how we can recover a power of $(\L,\xi)$
from such symmetric differences using a Taylor-like expansion.

\subsection{}\label{cubicext}
Let $n\geq 2$. For $1\leq i\leq n-1$  consider the morphisms
$A_i,B_i, C_i: H^{n}\to H^{n-1}$ given on points by 
\begin{eqnarray*}
A_i(h_1,h_2,\ldots, h_n)&=& (h_1, \ldots, \hat h_i,\ldots, h_n),\\
B_i(h_1,h_2,\ldots, h_n)&=&(h_1,\ldots, \hat h_{i+1},\ldots,h_n),\\
C_i(h_1,h_2,\ldots, h_n)&=&(h_1,\ldots, h_i+ h_{i+1},\ldots,  h_n),
\end{eqnarray*}
where $\hat h_j$ means ``omit $h_j$" and where in the last expression
$h_i+h_{i+1}$ is placed in the $i$-th position.
If $\L$ is an invertible sheaf on $H$, we can see from the 
definitions that there is a canonical isomorphism
\begin{equation}\label{nplus1}
\Theta_{n}(\L)\xrightarrow{\sim}C_i^*\Theta_{n-1}(\L) \otimes A_i^*\Theta_{n-1}(\L)^{-1}\otimes B_i^*\Theta_{n-1}(\L)^{-1}\ .
\end{equation}

Let now $(\L,\xi)$ be an invertible sheaf with an $n$-cubic structure 
over $H$. We will show how we can associate to the pair $(\L,\xi)$ an $n-1$-extension
$E(\L,\xi)$
of $H$ by $\Gm$. The corresponding $\Gm$-torsor on $H^{n-1}$ is 
given by $\Theta_{n-1}(\L)$. For $n=3$ a similar construction is described in [Br] \S 2; 
the general cases follows along the same lines. We sketch the argument below:
 By composing (\ref{nplus1}) with $\xi$ we obtain an isomorphism
\begin{equation}\label{compolaw}
c^i: A_i^*\Theta_{n-1}(\L)\otimes B_i^*\Theta_{n-1}(\L)\xrightarrow{\sim}C_i^*\Theta_{n-1}(\L) 
\end{equation}
of invertible sheaves on $H^n$. We can verify that these isomorphisms provide
the partial composition laws (\ref{isomulti}) of an $n-1$-extension: To check that the
$c^i$ are commutative, associative (as in (\ref{commu}) and (\ref{assoc})) and compatible with each other 
(as in (\ref{compa})) we can reduce to the
case that $R$ is local and $\L$ is the trivial  invertible sheaf on $H$
(\ref{various} (ii)-(iii)).
Then the $n$-cubic structure $\xi$ on $\L$ is given by an
$n$-cubic element $c\in R[G^n]^*$ and, by the above, the 
``composition law" $c^i$ is given via multiplication by the element $c$. Hence, we are reduced to checking
certain identities for $c$. These follow directly from 
properties (c1) and (c2) of \ref{various} (ii). More specifically, the commutativity, resp. associativity,
property for $c^i$ follows  directly from property (c1), resp. (c2), for $c$. 
The compatibility (\ref{compa}) between the partial composition laws 
$c^i$ for various $i$ also follows immediately from (c1) and (c2). As a result, the isomorphisms
$c^i$, $1\leq i\leq n-1$, provide $\Theta_{n-1}(\L)$ with the structure 
of an $n-1$-extension $E(\L,\xi)$. In fact, we can see that the construction $(\L,\xi)\mapsto E(\L,\xi)$
is functorial and gives an additive functor
\begin{equation*}
n{\hbox{\small -}}{\rm CUB}(H,\Gm)\to (n-1){\hbox{\small -}}{\rm EXT}(H,\Gm)\ .
\end{equation*}
Actually, we can see that the $n-1$-extension $E(\L,\xi)$ is symmetric
(in the sense of the previous paragraph) with the symmetry isomorphisms $\Psi_\sigma$ given by 
the isomorphisms
${\mathfrak P}_\sigma$ of (\ref{permu}) ($\sigma\in S_{n-1}$).
Indeed, it is easy to see that the isomorphisms ${\mathfrak P}_\sigma$ of $\Gm$-torsors satisfy the conditions (i) and (ii) of \ref{symmpar} and it remains to show that
they actually give (iso)morphisms of $n-1$-extensions. An argument as above now
shows (by reducing to the case $R$ local and $\L$ trivial) that this follows
from the definitions and property (c1).
\bigskip

\subsection{} \label{4a} 
In this paragraph, we assume that $n\geq 3$. Let $\L$ be an invertible sheaf over $H$ equipped with
an isomorphism 
\begin{equation*}
\xi':\O_{H^{n-1}}\xrightarrow{\sim} \Theta_{n-1}(\L)
\end{equation*}
For simplicity, we set $A'=A_1$, $B'=B_1$, $C'=C_1$
for the morphisms $H^n\to H^{n-1}$ of \ref{cubicext}.
Recall the canonical isomorphism (\ref{nplus1})
\begin{equation*}
\Theta_n(\L)\xrightarrow{\sim}{C'}^*\Theta_{n-1}(\L)\otimes{A'}^*\Theta_{n-1}(\L)^{-1}\otimes{B'}^*\Theta_{n-1}(\L)^{-1}.
\end{equation*} Define an isomorphism
\begin{equation}\label{induce}
\xi: \O_{H^{n}}\xrightarrow{\sim} \Theta_{n}(\L)
\end{equation}
by composing the inverse of (\ref{nplus1}) with the trivialization of 
$$
{C'}^*\Theta_{n-1}(\L)\otimes{A'}^*\Theta_{n-1}(\L)^{-1}\otimes{B'}^*\Theta_{n-1}(\L)^{-1}
$$
induced by $\xi':\O_{H^{n-1}}\xrightarrow{\sim} \Theta_{n-1}(\L)$. 
 
\begin{lemma}\label{61}
If $(\L, \xi')$ is a line bundle with an $n-1$-cubic structure on $H$, 
then $\xi$ given by (\ref{induce}) gives an $n$-cubic structure on $\L$.
In fact, the construction $(\L,\xi')\mapsto (\L, \xi)$ gives an additive functor
$$
{(n-1)}{\hbox{\small -}}{\rm CUB}(H,\Gm)\to {n}{\hbox{\small -}}{\rm CUB}(H,\Gm)\ .
$$
\end{lemma}

\begin{Proof}
To show the first statement we have to show that the isomorphism (\ref{induce}) above satisfies the
conditions (c0), (c1), (c2) of Definition \ref{cubicdef} for an $n$-cubic structure. For this purpose,
we may assume that $R$ is local and that $\L$ is the trivial invertible sheaf on $H$
(see \ref{various} (ii)).   Then the $n-1$-cubic 
structure $\xi'$ is given by
an $n-1$-cubic element $c'\in R[G^{n-1}]^*=\Gamma(H^{n-1},\O^*_{H^{n-1}})$ and we can see
that $\xi$ is given by $c\in R[G^n]^*=\Gamma(H^{n},\O^*_{H^{n}})$ which is defined
by
\begin{equation}\label{relation}
c={C'}^*({c'}){A'}^*(c')^{-1}{B'}^*(c')^{-1}\ .
\end{equation}
In other words, we have
\begin{multline}\label{relation'}
c(h_1,h_2,\ldots, h_n)=\\
=c'(h_1+h_2,h_3,\ldots, h_n)c'(h_1,h_3,\ldots, h_n)^{-1}c'(h_2,h_3,\ldots, h_n)^{-1}
\end{multline}
for all points $h_i$, $1\leq i\leq n$, of $H$. We now have to show that
if $c'$ satisfies (c0), (c1), (c2) of \ref{various} (ii)
with $n$ replaced by $n-1$, then $c$ satisfies (c0), (c1), (c2) for $n$: 
It is clear that $c$ satisfies (c0) and that $c$ is symmetric in the ``variables" $h_1$, $h_2$
and in $h_3,\ldots ,h_n$ separately. 
To show that $c$ satisfies (c1) in general, it is enough to show that, in addition, we have
\begin{equation}\label{symSpe}
c(h_1,h_2,h_3,\ldots, h_n)=c(h_1,h_3,h_2,\ldots, h_n).
\end{equation}
To explain this we may assume that $n=3$ (the argument for $n>3$ is essentially the same). 
By the cocycle condition (c2) for $c'$ we obtain: 
$c'(h_2+h_1,h_3)c'(h_1,h_3)^{-1}=c'(h_2, h_1+h_3)c'(h_2,h_1)^{-1}$.
By multiplying both sides with $c'(h_2,h_3)^{-1}$ and using the symmetry condition for $c'$ we obtain
(\ref{symSpe}) and this shows condition (c1) for $c$. The cocycle condition (c2) for $c$ now follows directly from
(\ref{relation'}). This proves the first statement of the Lemma.
To show the second statement  we first observe that our construction
is functorial. The rest follows from the definition of the product of multiextensions.
\endproof
\end{Proof}
\smallskip

\begin{lemma}\label{62}
Suppose that $(\L,\xi)$ is an invertible sheaf with an $n$-cubic 
structure over $H$ which is such that the corresponding $n-1$-extension 
$E(\L,\xi)$ of \ref{cubicext} is trivial as a {\sl symmetric} multiextension. 
Then there is an $n-1$-cubic structure $\xi'$ on $\L$ which
induces the $n$-cubic structure $\xi$ by the procedure 
of Lemma \ref{61}. Conversely, if the $n$-cubic structure $\xi$ is induced 
from an $n-1$-cubic structure $\xi'$ by the procedure of Lemma \ref{61} then
$E(\L, \xi)$ is trivial as a symmetric multiextension.
\end{lemma}
\begin{Proof}
(For $n=3$ and general $H$ this is essentially
[Br] Prop. 2.11.) Suppose that $E(\L,\xi)$ is trivial
as a symmetric $n-1$-extension. By definition, this means that
there is an isomorphism
\begin{equation*}
\xi': \O_{H^{n-1}}\xrightarrow{\sim} E(\L,\xi):=\Theta_{n-1}(\L)
\end{equation*}
which is compatible with the partial composition laws (\ref{isomulti})
and the symmetry isomorphisms 
\begin{equation*}
{\mathfrak P}_{\tau}: \tau^*\Theta_{n-1}(\L)\xrightarrow{\sim}\Theta_{n-1}(\L)
\end{equation*}
for all $\tau\in S_{n-1}$. Recall that the composition laws  
on $E(\L,\xi)$ are given by (\ref{compolaw}). We can now see that the isomorphism $\xi'$ 
is compatible with the composition law for $i=1$ if and only
if $\xi: \O_{H^n}\xrightarrow{\sim} \Theta_n(\L)$ is obtained from 
the isomorphism $\xi'$
by the procedure described in the beginning of \ref{4a}.
We just have to show that the isomorphism $\xi'$ defines an $n-1$-cubic 
structure. For this purpose,
we may assume that $R$ is local and that $\L$ is the trivial invertible sheaf on $H$
(see \ref{various} (ii)). As in the proof of the previous lemma, we see
that the isomorphisms $\xi$, $\xi'$ are given by elements $c\in R[G^n]^*$, $c'\in R[G^{n-1}]^*$
respectively which are related by (\ref{relation'}). Since $\xi$ is an $n$-cubic structure,
$c$ satisfies (c0), (c1) and (c2) of \ref{various} (ii). We would like to show
that $c'$ satisfies (c0), (c1) and (c2) with $n$ replaced by $n-1$. Property (c0)
follows immediately from (\ref{relation'}). 
Since $\xi'$ is compatible with the symmetry isomorphisms $c'$ satisfies
(c1). It remains to show property (c2); the relevant equation can be written
\begin{multline}
c'(h_1+h_2,h_3,\ldots, h_n)c'(h_2,h_3,\ldots, h_n)^{-1}=\\=c'(h_1,h_2+h_3,\ldots, h_n)c'(h_1,h_2,\ldots, h_n)^{-1}\ . 
\end{multline}
This now follows from Property (c1) for $c$ and (\ref{relation'}).
We will leave the converse to the reader.\endproof
\end{Proof}
\smallskip
\begin{Remark} \label{extra}
{\rm Note that in the paragraph above we assumed that $n\geq 3$. 
Suppose now that $n=2$. Then we have $\Theta_{n-1}(\L)=\Theta_1(\L)=\L\otimes 0^*\L^{-1}$.
Hence, if $E(\L,\xi)$ is a trivial ($1$-)extension and $0^*\L$ a trivial invertible $\O_H$-sheaf
then $\L$ is also a trivial invertible $\O_H$-sheaf.}
\end{Remark}

\subsection{} \label{taylorsect} Let $n\geq 1$.
Suppose that $(\L,\xi)$ is an invertible sheaf with an $n+1$-cubic 
structure over $H$. If $\Delta_n: H\to H^n$ is the diagonal morphism, then we can
consider the invertible sheaf $\delta(\L, \xi):=\Delta^*_n\Theta_n(\L)=\Delta_n^*E(\L,\xi)$ on
$H$ and set
\begin{equation}
\L^\flat:=\delta(\L, \xi)\otimes\L^{\otimes -n!}\ .
\end{equation}

\begin{prop}\label{flat}
Suppose that $n\geq 2$. Then the invertible sheaf $\L^\flat$
defined above is equipped with a canonical $n$-cubic structure
$\xi^\flat$.
\end{prop}

\begin{Remark}\label{flat2a} {\rm    
Notice that the invertible sheaf $\delta(\L,\xi)$ is always ``rigid", i.e equipped
with an isomorphism $0^*\delta(\L,\xi)\simeq \O_H$. Hence, there is an isomorphism
$$
0^*\L^\flat\simeq 0^*\L^{-\otimes n!}.
$$
Also notice that when $n=1$, we have $\L^{\flat}=\delta(\L,\xi)\otimes\L^{-1}=\Theta_1(\L)\otimes\L^{-1}=0^*\L^{-1}$.}
\end{Remark}

Notice that successive application of Proposition \ref{flat}, 
combined with the above remark, gives the following. (See Remark \ref{RRexpl} for an interpretation 
of a special case of this formula in the context of a Riemann-Roch theorem.)

\begin{cor}\label{flat2} { {\bf (Taylor expansion)}\ 
There is an isomorphism of invertible sheaves
\begin{equation}
\L^{\otimes n!!}\simeq  (0^*\L)^{\otimes n!!}
\otimes\bigotimes_{i=0}^{n-1} \delta(\L^{(i)},\xi^{(i)})^{\otimes (-1)^i{(n-i-1)}!!}
\end{equation}
where $(\L^{(0)},\xi^{(0)}):=(\L, \xi)$, $(\L^{(i)}, \xi^{(i)}):=((\L^{(i-1)})^\flat, (\xi^{(i-1)})^\flat)$ 
and $n!!=n!(n-1)!\cdots 2!$. }
\end{cor}

\noindent{\sc Proof of \ref{flat}.} 
For simplicity, we set $E=E(\L,\xi)$, $\delta=\delta(\L,\xi)$. If $R'$ is an $R$-algebra
we consider $H(R')=H(\Spec(R'))$; this is the group of  characters of $G$ with values in 
$R'$. Let $\chi_0, \chi_1,\ldots,\chi_n$ be $R'$-valued characters of $G$.
If $S$ is a subset of $\{0,\ldots, n\}$,
we set $\chi_S=\prod_{i\in S}\chi_i$. (Here and below a product, resp. a tensor product,
over the empty set is $1$, resp. the trivial invertible sheaf.) 
By the definition, we have
\begin{equation}\label{1}
\Theta_{n}(\delta)_{(\chi_1,\ldots,\chi_n)}=
\bigotimes_{S\subset\{1,\ldots, n\}}E_{(\chi_S, \ldots, \chi_S)}{}^{(-1)^{n-\#S}}.
\end{equation}
Repeated application of the composition laws (\ref{isomulti}) now provide 
functorial isomorphisms
\begin{equation} \label{2}
\bigotimes_{p: \{1,\ldots, n\}\to S}E_{(\chi_{p(1)},\ldots, \chi_{p(n)})}\xrightarrow{\sim}E_{(\chi_S, \ldots, \chi_S)} \  \end{equation}
(the tensor product runs over all maps $p:\{1,\ldots, n\}\to S$). 
Observe that if $S\neq \{1,\ldots, n\}$ we have
\begin{equation}
\sum_{S'; S\subset S'\subset \{1,\ldots, n\}}(-1)^{n-\# S'}=0\ .
\end{equation}
This shows that in the tensor product
\begin{equation}
\bigotimes_{S\subset\{1,\ldots, n\}}\bigotimes_{p: \{1,\ldots, n\}\to S}E_{(\chi_{p(1)},\ldots, \chi_{p(n)})}{}^{(-1)^{n-\#S}}
\end{equation}
the terms for which either $S\neq \{1,\ldots, n\}$ or $S=\{1,\ldots, n\}$
and $p$ is not surjective contract (canonically). Therefore, we are left with
\begin{equation}
\bigotimes_{p: \{1,\ldots, n\}\xrightarrow{\sim} \{1,\ldots, n\}}E_{(\chi_{p(1)},\ldots, \chi_{p(n)})}\ .
\end{equation}
Hence, using the symmetry isomorphisms and (\ref{1}) we can see that there is a canonical isomorphism
\begin{equation}\label{5}
(E_{(\chi_1,\ldots, \chi_n)})^{\otimes n!}\xrightarrow{\sim} \Theta_{n}(\delta)_{(\chi_1,\ldots,\chi_n)}\ .
\end{equation}
Since by definition $E=\Theta_n(\L)$ we obtain from (\ref{5}) a canonical
isomorphism
\begin{equation} \label{isom}
\xi^\flat: \O_{H^n}\xrightarrow{\sim}  \Theta_n(\L^\flat)=\Theta_n(\delta\otimes\L^{-\otimes n!})\simeq \Theta_n(\delta) \otimes \Theta_n(\L)^{-\otimes n!}\ .
\end{equation}
We will now show that the isomorphism (\ref{isom}) above satisfies the
conditions (c0), (c1), (c2) of an $n$-cubic structure. For this purpose,
we may assume that $R$ is local and that $\L$ is the trivial invertible sheaf on $H$
(see \ref{various} (ii)-(iii)). Then all the
invertible sheaves in the construction above are also trivial and the hypercubic structure $\xi$ is given by
an $n+1$-cubic element $c\in R[G^{n+1}]^*$. By unraveling the definition
above we can now see that the isomorphism (\ref{isom}) is given as multiplication by 
\begin{equation}
\prod_{S\subset \{1,\ldots, n\}} (\chi_1\otimes\cdots \otimes\chi_n)(d_S)^{(-1)^{n-\#S}}\ 
\end{equation}
where $d_S\in R[G^n]^*$ and the term $(\chi_1\otimes\cdots \otimes\chi_n)(d_S)$
gives the isomorphism (\ref{2}). (The element $d_S$ gives the isomorphism (\ref{2}) 
for $\chi_i$, $i=1,\ldots ,n$, the ``universal" $R[G^n]$-valued characters
$G\to G^n\subset R[G^n]$, given by $\chi_i(g)= (g_j)_j$, with $g_j=g$ if $j=i$, $g_j=1$ 

if $j\neq i$. Notice that if $\#S\leq 1$, then $d_S=1$.)
In fact, it is more convenient to consider the inverse of (\ref{2})
and view that as the composition
of several isomorphisms in which the arguments $\chi_S$
are unraveled one by one. Suppose that 
$S=\{i_1<i_2<\cdots <i_m\}\neq \emptyset$. Then the first of these isomorphisms
is 
\begin{equation} \label{6}
 E_{(\chi_S, \ldots, \chi_S)}\xrightarrow{\sim}\bigotimes_{k=1}^mE_{(\chi_{i_k},\chi_S,\ldots, \chi_{S})}\ .
\end{equation}
By the definition of the composition law
of the $n$-extension $E=E(\L,\xi)$  (see \ref{cubicext}), this is described by the inverse of the element:
\begin{eqnarray*}
(\chi_{i_1}\otimes \prod_{k>1}\chi_{i_k}\otimes\chi_S\otimes&\cdots&\otimes\chi_S)(c)\cdot\\
\cdot(\chi_{i_2}\otimes \prod_{k>2}\chi_{i_k}\otimes\chi_S\otimes&\cdots&\otimes\chi_S)(c)\cdot\\
&\vdots&\\
\cdot(\chi_{i_{m-1}}\otimes\chi_{i_m}\otimes\chi_S\otimes&\cdots&\otimes\chi_S)(c) .
\end{eqnarray*}
Using   \ref{augmenG} we see that we can write this 
as the value of the element
\begin{eqnarray*}
 \biggl(\sum_{p=1}^{m-1}\bigl\{([\chi_{i_p}]-[1])
\otimes([\prod_{k>p}\chi_{i_k}]-[1])\bigr\}\biggr)\otimes ([\chi_S]-[1])\otimes \cdots 
\otimes([\chi_S]-[1])
\end{eqnarray*}
of $C_{n+1}(H(R'))$ at $c^{-1}$. For simplicity, we set 
\begin{equation}
  A_S=\sum_{k=1}^m([\chi_{i_k}]-[1])\ ,\quad B_S=
\sum_{p=1}^{m-1}\bigl\{([\chi_{i_p}]-[1])\otimes([\prod_{k>p}\chi_{i_k}]-[1])\bigr\}\ .
\end{equation}
Similarly, we can now see that 
the isomorphisms
\begin{equation}\label{7}
E_{(\chi_{i_k},\chi_S\ldots, \chi_{S})}\xrightarrow{\sim} 
 \bigotimes_{p=1}^mE_{(\chi_{i_k},\chi_{i_p},\chi_S,\ldots, \chi_{S})}
\end{equation}
which give the next step in unraveling the inverse of (\ref{2}) are described by evaluating at $c^{-1}$  the element
\begin{equation}
([\chi_{i_k}]-[1])\otimes B_S\otimes ([\chi_S]-[1])\otimes \cdots 
\otimes([\chi_S]-[1]) .
\end{equation}
The combined effect (for $k=1,\cdots, m$) of all of these on the tensor product 
of (\ref{6}) is given by evaluating at $c^{-1}$  the element
\begin{equation} A_S\otimes B_S\otimes (([\chi_S]-[1]))^{\otimes (n-2)}\ .
\end{equation}
The next step is unraveling the first remaining $\chi_S$ in 
$E_{(\chi_{i_k},\chi_{i_p},\chi_S,\ldots, \chi_{S})}$. As above, we can see that this is given by the elements 
\begin{equation}
([\chi_{i_k}]-[1])\otimes ([\chi_{i_p}]-[1])\otimes B_S\otimes ([\chi_S]-[1])^{\otimes (n-3)} 
\end{equation}
with combined effect 
\begin{equation}
A_S^{\otimes 2}\otimes B_S\otimes ([\chi_S]-[1])^{\otimes (n-3)}
\end{equation}
and so on.
Putting everything together we see that $(\chi_1\otimes\cdots\otimes\chi_n)(d_S)$ is given by
evaluating the element
\begin{equation}\label{psi}
\Psi_S(\chi_1,\cdots,\chi_n)=\sum_{j=0}^{n-1} A_S^{\otimes j}\otimes B_S\otimes 
([\chi_S]-[1])^{\otimes (n-j-1)} 
\end{equation}
of $C_{n+1}(H(R'))$ at $c$. Hence, the isomorphism (\ref{isom}) is given 
by an element $d\in R[G^n]^*$ which is such that
\begin{equation}
(\chi_1\otimes\cdots\otimes \chi_n)(d)=\biggl(\sum_{\substack{S\subset \{1,\ldots, n\}\\
S\neq\emptyset}}(-1)^{n-\#S}\Psi_S(\chi_1,\ldots,\chi_n)\biggr)(c)\ .
\end{equation}
For simplicity, set
\begin{equation}\label{phi}
\Phi(\chi_1,\ldots, \chi_n)=\sum_{ S\subset \{1,\ldots, n\}}(-1)^{n-\#S}\Psi_S(\chi_1,\ldots,\chi_n)\ 
\end{equation}
in $C_{n+1}(H(R'))$ (here by definition $\Psi_{\emptyset}=0$).
The proof of the proposition will follow if we show the following properties:
 
\begin{equation*}
\begin{split}
{\rm f}0)\qquad \ \ \ \ \Phi(1,\ldots, 1)=0,\qquad\qquad\qquad\qquad \qquad \qquad \qquad \qquad \qquad \ \ \ \ \ \ \ \ \\
{\rm f}1)\qquad \ \ \ \ \Phi(\chi_{\sigma(1)},\ldots, \chi_{\sigma(n)})=\Phi(\chi_1,\ldots, \chi_n),
\hbox{\rm\  for all $\sigma\in S_n$,}\ \ \ \ \ \ \ \ \ \ \ \ \ \ \ \ \\
{\rm f}2)\qquad \ \ \ \ \Phi(\chi_0\chi_1,\chi_2,\ldots,\chi_n)+\Phi(\chi_0,\chi_1,\chi_3,\ldots,\chi_n)=\ \ \ \ \ \ \ \ \ \ \ \ \ \ \ \ \ \ \ \ \ \ \ \\
=\Phi(\chi_0,\chi_1\chi_2,\chi_3,\ldots,\chi_n)+\Phi(\chi_1,\chi_2,\ldots,\chi_n).\ \ \ \ \ \ \ \ \ \
\end{split}
\end{equation*}
Property (f0) is obvious and it is enough to concentrate on (f1) and (f2). 
Let $F=\{\prod_{i=0}^nx_i^{k_i}\ |\ k_i\in\Z\}$ be the free abelian group generated by 
the symbols $x_0,x_1,\ldots, x_n$ and  let us consider 
$C_k(F):={\rm Sym}^k_{\Z[F]}I[F]$, for $k\geq 1$. Recall that we denote by $I[F]^k$
the $k$-th power of the augmentation ideal $I[F]$ of the group ring $\Z[F]$.

\begin{lemma} \label{ideal}
The multiplication morphism $a_1\otimes\cdots\otimes a_k\mapsto a_1\cdots a_k$
induces an isomorphism
$$
C_k(F)={\rm Sym}^k_{\Z[F]}I[F]\xrightarrow{\sim} I[F]^k\subset \Z[F] .
$$
\end{lemma}

\begin{Proof} 
In this case, $\Z[F]\simeq \Z[u_0,u_0^{-1},\ldots ,u_{n}, u_{n}^{-1}]$ (the
ring of Laurent polynomials in $n+1$ indeterminants) with $I[F]$ corresponding to the ideal 
$(u_0-1,\ldots, u_{n}-1)$. Consider the ideal $I=(v_0,\ldots,v_{n})$ in the polynomial
ring $\Z[\underline v]:=\Z[v_0,\ldots,v_n]$. Multiplication
${\rm Sym}^k_{\Z[\underline v]}I\xrightarrow {}I^k$ gives an isomorphism
and the desired statement follows from this fact by setting $v_i=u_i-1$ and localizing.\endproof
\end{Proof}
\medskip

Suppose that $y_i$, $1\leq i\leq n$,
are elements of $F$. The identities (\ref{psi}), (\ref{phi}) with $\chi_i$ replaced by
$y_i$ can be used to define elements $\Psi_S(y_1,\ldots, y_n)$,
$\Phi(y_1,\ldots, y_n)\in C_{n+1}(F)$. The group homomorphism $F\to H(R')$ given by $x_i\mapsto \chi_i$ induces 
a homomorphism 
$C_{n+1}(F)\to C_{n+1}(H(R'))$ which sends the elements $\Phi(x_1,\ldots, x_n)$, $\Phi(x_0x_1,\ldots, x_n)$,
 to $\Phi(\chi_1,\ldots, \chi_n)$, $\Phi(\chi_0\chi_1,\ldots, \chi_n)$ etc.
Hence, it is enough to show
\smallskip
\begin{equation*}
\begin{split}
{\rm g}1)\qquad \ \ \ \, \Phi(x_{\sigma(1)},\ldots, x_{\sigma(n)})=\Phi(x_1,\ldots, x_n),
\hbox{\rm\  for all $\sigma\in S_n$,}\ \ \ \ \ \ \ \ \ \ \ \ \ \ \\
{\rm g}2)\qquad \ \ \ \Phi(x_0x_1,x_2,\ldots,x_n)+\Phi(x_0,x_1,x_3,\ldots,x_n)=\ \ \ \ \ \ \ \ \ \ \ \ \ \ \ \ \ \ \ \ \ \ \\
=\Phi(x_0,x_1x_2,x_3,\ldots,x_n)+\Phi(x_1,x_2,\ldots,x_n).\ \ \ \ \ \ \ \ 
\end{split}
\end{equation*}

In what follows, we will identify $C_k(F)$ with $I[F]^k$ using the
multiplication morphism of Lemma \ref{ideal}. Furthermore, we will
eliminate the brackets from the notation of elements of the group ring $\Z[F]$.

As above, if $S=\{i_{1}<\cdots <i_m\}\subset \{1,\ldots ,n\}$, we set
\begin{equation}
\begin{split}
A_S=A_S(y_1,\ldots, y_n)&=\sum_{k=1}^my_{i_k}-m ,\\ 
B_S=B_S(y_1,\ldots, y_n)&=\sum_{p=1}^{m-1} \bigl\{(y_{i_p}-1)(\prod_{k>p}y_{i_k}-1)\bigr\} ,
\end{split}
\end{equation}
and also
\begin{equation}
P_S=P_S(y_1,\ldots, y_n)=\prod_{k=1}^my_{i_k}-1\ .
\end{equation}

\begin{lemma}\label{powers}
$\ \Psi_S(y_1,\ldots, y_n)=P_S^n-A_S^n\ .$
\end{lemma}

\begin{Proof} 
By the definition of $\Psi_S(y_1,\ldots, y_n)$ we have
\begin{equation}
\Psi_S(\chi_1,\cdots,\chi_n)=\sum_{j=0}^{n-1}  A_S^{ j} B_S
P_S^{n-j-1}=B_S\cdot(\sum_{j=0}^{n-1}  A_S^{j}
P_S^{n-j-1}).
\end{equation}
However, observe that by telescoping we find
\begin{equation}
\begin{split}
B_S(y_1,\ldots, y_n)=&\sum_{p=1}^{m-1} \bigl\{(y_{i_p}-1)(\prod_{k>p}y_{i_k}-1)\bigr\}=\\
=&(\prod_{k=1}^my_{i_k}-1)-(\sum_{k=1}^my_{i_k}-m)=P_S-A_S .
\end{split}
\end{equation}
The result now follows from the standard identity.\endproof
\end{Proof}
\medskip

Lemma \ref{powers} and the definition of $\Phi(y_1,\ldots, y_n)$ (cf. (\ref{phi})) now imply
\begin{equation}\label{phiform}
\Phi(y_1,\ldots, y_n)=\sum_{S\subset \{1,\ldots, n\}}(-1)^{n-\#S}((\prod_{i\in S}y_i-1)^n-(\sum_{i\in S}y_i-\#S)^n)\ .
\end{equation}
\smallskip

Notice that (g1) now follows immediately.
It remains to show (g2). To do that we will compare
terms between the two sides of the equation (g2). 
Let us consider the left hand side of the equation.
Using (\ref{phiform}) above we can see that it is a sum:
\begin{equation*}
Z+Y(x_0x_1,x_2,\ldots,x_n)+Y(x_0,x_1,x_3,\ldots,x_n),
\end{equation*}
 where we set
\begin{equation*}
Y(y_1,\ldots, y_n)=-\sum_{ S\subset \{1,\ldots, n\}
 }(-1)^{n-\#S}(\sum_{i\in S}(y_i-1))^n\ ,
\end{equation*}
and where $Z$ is a sum of terms which are either of the form $(-1)^{n-(\#T-1)}(\prod_{i\in T}x_i-1)^n$
or of the form $(-1)^{n-\#T}(\prod_{i\in T}x_i-1)^n$, for certain $T\subset\{0,1,\ldots, n\}$.
More specifically, we can write
$$
Z=2\Sigma_1+\Sigma_2+\Sigma_3.
$$
where

$$\Sigma_1=\sum_{  \{0,1,2\}\cap T=\emptyset}(-1)^{n-\#T}(\prod_{i\in T}x_i-1)^n,$$

$$\Sigma_2=\sum_{ \{0,1,2\}\subset T }(-1)^{n-(\#T-1)}(\prod_{i\in T}x_i-1)^n,$$  

$$\Sigma_3=\sum_{  \#(T\cap\{0,1,2\})=1}(-1)^{n-\#T}(\prod_{i\in T}x_i-1)^n.$$

\noi (In these sums,
 $T$ ranges over subsets of $\{0,1,\ldots, n\}$.) The rest of the terms cancel out since they appear twice but with different signs
in $\Phi(x_0x_1,x_2,\ldots,x_n)$ and $\Phi(x_0,x_1,\ldots,x_n)$.

Similarly, a careful look at the right hand side of (g2)
reveals that it is 
equal to the sum
\begin{equation*}
Z+Y(x_0,x_1x_2,x_3,\ldots,x_n)+Y(x_1,x_2,\ldots,x_n)\ .
\end{equation*}

Now observe that the identity
$$
\sum_{S\subset \{1,\ldots, n\}}(-1)^{n-\#S}(\sum_{i\in S}z_i)^n=n!z_1z_2\cdots z_n\ 
$$
gives
$$
Y(y_1,\ldots, y_n)=-n!\prod_{k=1}^n(y_i-1) \ .
$$
Hence, by the above, we can now conclude that (g2) is equivalent to the identity
\begin{multline}
(x_0x_1-1)(x_2-1)\cdots 
(x_n-1)+(x_0-1)(x_1-1)\cdots (x_n-1)=\\=(x_0-1)(x_1x_2-1)\cdots (x_n-1)+(x_1-1)(x_2-1)\cdots (x_n-1)\ ,
\end{multline}
which is easily seen to be true. This concludes the 
proof of the identity (g2) and of Proposition \ref{flat}.
\endproof
\bigskip

\section{Multiextensions and abelian sheaves} \label{sheaves}

In the next few paragraphs, we will not distinguish in our notation 
between a commutative group scheme $A$ over a scheme $T$ and the sheaf 
of abelian groups on the site $T_{\rm fppf}$ which is given by the sections
of $A$. We will use the derived category of the homotopy category of 
complexes of sheaves of abelian groups on $T_{\rm fppf}$
(recall that $T_{\rm fppf}$ is the site of $T$-schemes which are locally of finite presentation
with the fppf topology). Otherwise, we continue with the 
notations and general set-up of \S \ref{multiext}.
If $Y\to S$ is an object of $S_{\rm fppf}$,
we will denote by $\Z[Y]$ the abelian sheaf
on the fppf site of $S$ which is freely generated by the points of $Y$
(see [SGA4] IV 11).  If $Y$ and $Y'$ are two objects of $S_{\rm fppf}$ 
then there is a canonical
isomorphism: 
\begin{equation}\label{product}
\Z[Y]\Lder \Z[Y']=\Z[Y]\otimes
\Z[Y']\simeq \Z[Y\times_SY']. \end{equation}
The natural morphism of sheaves $Y\to \Z[Y]$ induces a canonical isomorphism ([SGA7I] VII 1.4):
\begin{equation} \label{exttorsor}
\Ext^1(\Z[Y], J)\xrightarrow{\sim} {\rm H}^1(Y, J_Y).
\end{equation} 
Now suppose that $H\to S$ is a commutative group scheme
which is finite locally free over $S$.
We will denote by $\ep_H: \Z[H]\to H$ the augmentation homomorphism.
Then, under \ref{exttorsor}, the homomorphism
$\Ext^1(\ep_H, J)$ is identified with the 
natural homomorphism
\begin{equation} \label{extforget}
\Ext^1(H, J)\to {\rm H}^1(H, J).
\end{equation}

Suppose $F$, $F'$ are abelian sheaves (on $S_{\rm fppf}$) and $E$
a complex of abelian sheaves which is bounded above. Then there is
a canonical spectral sequence
\begin{equation}\label{spectral}
\Ext^p(E, \underline{\Ext}^q(F, F')) \Rightarrow
\Ext^{p+q}(E, \underline{{\bf R}{\Hom}}(F, F'))
\end{equation}
which induces an exact sequence:
\begin{equation}\label{exactspectral}
 0\to \Ext^{1}(E, \underline{{\Hom}}(F, F'))
\to\Ext^{1}(E, \underline{{\bf R}{\Hom}}(F, F'))\to
\Ext^{1}(E, \underline{\Ext}^1(F, F')).
\end{equation}
There is also a canonical isomorphism:
\begin{equation} \label{flip}
\Ext^p(E, \underline{{\bf R}{\Hom}}(F, F'))\loniso \Ext^p(E\Lder F, F').
\end{equation}

\subsection{} By [SGA7I] VII (2.5.4.1) and 3.6.7
(see also loc. cit. 3.6.4, 3.6.5 and the remarks in VIII \S 0.2) there are canonical isomorphisms
\begin{equation}\label{nHom}
n{\hbox{\small -}}{\rm Ext}^0(H,J)\xrightarrow{\sim} 
{\rm Hom}(\underset{n} {\underbrace{H\otimes \cdots\otimes H}},J)\ ,
\end{equation}
\begin{equation}\label{nExt}
n{\hbox{\small -}}{\rm Ext}^1(H,J)\xrightarrow{\sim} {\rm Ext}^1(\underset{n} {\underbrace{H\Lotimes \cdots\Lotimes H}},J)\ .
\end{equation}
In the source of the second morphism $H\Lotimes \cdots \Lotimes H
=((H\Lotimes\cdots \Lotimes H)\Lotimes H)\Lotimes H$ is 
the complex, well-defined up to canonical isomorphism in the derived category,
obtained by applying successively 
the derived tensor product functor. 

When $J=\Gm$, the discussion in loc. cit. shows
that the diagram
\begin{equation}
\begin{matrix}
n{\hbox{\small -}}{\rm Ext}^1(H,\Gm)&\xrightarrow{\sim}& {\rm Ext}^1( H\Lotimes \cdots\Lotimes H,\Gm)\\
\downarrow{t}&&\downarrow\\
{\rm Pic}(H\times_S\cdots\times_SH)&\xrightarrow{\sim} & {\rm Ext}^1( \Z[H]\otimes \cdots\otimes \Z[H],\Gm)\\
\end{matrix}
\end{equation}
commutes. Here the second vertical isomorphism is $\Ext^1(\ep_H\Lotimes\cdots\Lotimes \ep_H,\Gm)$,
and the lower horizontal isomorphism is given by (\ref{exttorsor}) and (\ref{product}).

\subsection{}
We continue to assume that $H\to S$ is finite locally free.
Once again, we denote by $H^D=\underline{\rm Hom}(H, \Gm)$ the Cartier dual of $H$; 
let $\{\ ,\ \}: H^D\times H\to \Gm$ be the natural pairing. By [SGA7I] VIII
Prop. 3.3.1, $\underline{\rm Ext}^1(H^D,\Gm)=(0)$. Then, the exact
sequence (\ref{exactspectral}) gives an isomorphism
\begin{equation}\label{above}
\Ext^{1}(E, H^D)\xrightarrow{\sim}
\Ext^{1}(E, \underline{{\bf R}{\Hom}}(H, \Gm)) .
\end{equation}
By composing (\ref{nExt}) with  
(\ref{above}) and (\ref{flip}) we obtain 
a canonical isomorphism
\begin{equation}\label{nExt2}
\Ext^1(\underset{n-1} {\underbrace{H\otimes^{\rm L} \cdots\otimes^{\rm L}H}}, H^D)\xrightarrow{\sim}n{\hbox{\small -}}{\rm Ext}^1(H,\Gm),
\end{equation} 
hence also
\begin{equation}\label{nExt3}
{(n-1)}{\hbox{\small -}}{\rm Ext}^1(H,H^D)\xrightarrow{\sim} n{\hbox{\small -}}{\rm Ext}^1(H,\Gm).
\end{equation}
For $n=1$, (\ref{nExt2}) amounts to an isomorphism
\begin{equation}\label{nExt4}
{\rm H}^1(S, H^D)\simeq {\rm Ext}^1(\Z,H^D) \xrightarrow{\sim}{\rm Ext}^1(H,\Gm) .
\end{equation}
To describe this last isomorphism explicitly, suppose we start with
an $H^D$-torsor $Q\to S$ which, under (\ref{exttorsor}),
corresponds to the extension
\begin{equation*} 0\to H^D\to Q'\to \Z\to 0 .
\end{equation*}
Tensoring with $H$ gives
an extension
\begin{equation*}
0\to H^D\otimes H\to Q'\otimes H\to H\to 0
\end{equation*}
which we can push out by $H^D\otimes H\to \Gm; a\otimes h\mapsto \{a,h\}$, to obtain an extension
of $H$ by $\Gm$. We can see that this push-out extension is isomorphic to
\begin{equation}\label{water}
1\to \Gm\to  (Q\times_S H\times_S \Gm)/H^D\to H\to 0\, .
\end{equation}
Here the (representable) fppf sheaf in the middle is the quotient sheaf
for the action of $H^D$ on $Q\times_S H\times_S \Gm$
given by $ (q, h, u)\cdot a=( q\cdot a, h, \{a,h\}^{-1}u)$
and has group structure given via descent by $(q,h,u)\cdot (q',h', u')=(q,hh',\{q^{-1}q',h'\}uu') $.
The isomorphism (\ref{nExt4}) and the explicit extension
(\ref{water}) are discussed in detail in [Wa]; see Theorems 2 and 3.
When $H^D$ is constant, then 
loc. cit. Theorem 3 implies that the extension (\ref{water}) is the negative of the
extension (\ref{groth}) which was associated 
to the $H^D$-torsor $Q$ in \S \ref{torsors1}; cf. Remark \ref{LchiRemark} (b)
below.  

In fact, we can obtain a similar description for the 
map (\ref{nExt3}) (cf. [SGA7I] VIII (1.1.6) where 
the details of this construction for $n=2$ are left 
to the reader):   Suppose that $Q\to H^{n-1}$ is the $H^D_{H^{n-1}}$-torsor
supporting the structure of an $n-1$-extension of $H$ by $H^D$. The 
construction (\ref{water}) applied to $S=H^{n-1}$
provides us with an extension of $H_{H^{n-1}}$ by $\Gm_{H^{n-1}}$ .
The underlying $\Gm_{H^n}$-torsor over $H_{H^{n-1}}=H^n$ 
then supports a canonical structure of $n$-extension whose isomorphism class
is the image of the class of $Q$ under the map (\ref{nExt3}).

\subsection{} 
We continue to assume that $H\to S$ is finite locally free.
For future use we observe that the following diagram
is commutative:
\begin{equation}\label{bigcd}
\begin{matrix}
n{\hbox{\small -}}{\rm Ext}^1(H,\Gm)&\xrightarrow{t} &{\rm Pic}(H^n) &\xrightarrow{\ \Delta_n^*\ } &{\rm Pic}(H)\\
&&&&\\
(\ref{nExt3})\Big\uparrow\wr\ \ \ \ \ \ \ &&&&  \ \ \Big\uparrow\Delta^*_2\\
&&&&\\
(n-1){\hbox{\small -}}{\rm Ext}^1(H,H^D) &&&& {\rm Pic}(H\times H)\\
&&&&\\
(\ref{4.7})\Big\downarrow t\ \ \ \ &&&&  \Big\uparrow t\\
&&&&\\
{\rm H}^1(H^{n-1},H^D)&\xrightarrow{\Delta_{n-1}^*}&  {\rm H}^1(H,H^D)&
\xrightarrow{\underset{\sim}{(\ref{nExt4})}}&\Ext^1(H_H,\Gm_H). \\
\end{matrix}
\end{equation}
This follows from the description of the maps (\ref{nExt3}), (\ref{nExt4})
in the previous paragraph.

\subsection{} \label{Lsection} Suppose now that $S=\Spec(R)$ and $H=G^D_S=\Spec(R[G])$ is the Cartier dual
of the finite constant abelian group scheme $G_S$. Let $T\to S$ be an $S$-scheme; 
Suppose $ q: Q\to T$ is a $G$-torsor; the construction (\ref{water}) gives a corresponding extension
of $G^D_T$ by $\Gm_T$. Suppose that $S'=\Spec(R') \to S$ is another $S$-scheme and consider a character
$\chi: G\to {R'}^*$; this corresponds to a point $S'\to G^D_S$ which we will still denote by $\chi$.
 Now suppose that 
$T'$ is an $S'\times_ST$-scheme
and consider the morphism $f: T' \to S'\times_ST\xrightarrow{(\chi,{\rm id})}G^D_T=G^D_S\times_ST$.
By pulling back the $\Gm_T$-torsor underlying 
the extension (\ref{water}) along 
$f: T' \to G^D_T$ we obtain a $\Gm_{T'}$-torsor (i.e an invertible
sheaf) $\L^Q_{f}$ over $T'$.  
By definition, the class of $\L^Q_{f}$
is the image of the class of $Q$
under the composition 
\begin{equation}\label{Lchimap}
{\rm H}^1(T,G)
\xrightarrow{(\ref{nExt4})}\Ext^1(G^D_T,\Gm_T)
\xrightarrow{\ }\Pic(G^D_T)\xrightarrow{{f}^*}\Pic(T') .
\end{equation}

Now recall (\S \ref{2d}) that $q_*(\O_Q)$ is  actually
a coherent sheaf of $\O_T[G]=\O_T\otimes_R R[G]$-modules;
we may think of it as a coherent 
$\O_{G^D_T}$-module which by (\ref{27}) is invertible. 

\begin{lemma}\label{Lchi}
 Let $G$ act on $q_*(\O_Q) \otimes_{\O_{T}}\O_{T'}$ via
$g\cdot (b\otimes t')=g\cdot b\otimes \chi(g)t'$. 
The sheaf of invariants $(q_*(\O_Q)\otimes_{\O_{T}}\O_{T'})^G$ is an
invertible sheaf of $\O_{T'}$-modules  
and we have
$$
\L^Q_{f}\simeq  (q_*(\O_Q)\otimes_{\O_{T}}\O_{T'})^G.
$$
\end{lemma}

\begin{Proof} This is a special case of [Wa] Theorem 3.
It   also follows directly from the explicit description
of the middle sheaf in the extension \ref{water} as a quotient
and the fact that in this case of free $G$-action taking quotient
commutes with base change (see (\ref{basechange})). \endproof
\end{Proof} 

\begin{Remark}\label{LchiRemark}
{\rm a) Suppose that we take $T'=S'\times_S T$ and $f=(\chi,{\rm id})$.
For simplicity, we set $\L^Q_{\chi}:=\L^Q_{(\chi, {\rm id})}$.
Then, Lemma \ref{Lchi} implies that
\begin{equation}\label{compa1}
\L^Q_{\chi}\simeq \O_{Q,\chi^{-1}}
\end{equation}
with the right hand side defined in \S \ref{2d}.  
\medskip

b) Recall the description of commutative extensions given in \ref{extension}.
Consider the points $(\chi, {\rm id}): S'\to G^D_T$
given by $S'$-valued characters $\chi: G\to \O_{S'}^*$
as in (a) above.
By \ref{extension}, and the definition of $\L^Q_{\chi}$,
there are functorial isomorphisms
\begin{equation}\label{Lmulti}
\L^Q_{\chi}\otimes_{\O_{S'}} \L^Q_{\chi'}\xrightarrow{\sim} \L^Q_{\chi\chi'}
\end{equation}
for any pair of characters $\chi$, $\chi'$.
Theorem 3 of [Wa] implies that (\ref{compa1})
takes the isomorphisms (\ref{Lmulti}) above to the isomorphisms (\ref{torsorIso})
defined in \S \ref{2d} using the multiplication of $\O_Q$.
In other words, (\ref{compa1}) is actually obtained  from
an isomorphism between the extension (\ref{water})
and the negative of the extension (\ref{groth}); we are not going to need
this more precise statement.}
\end{Remark} 
\bigskip

\section{Multiextensions of finite multiplicative group schemes}\label{overZ}

Suppose that $S=\Spec(\Z)$ and $H=G^D$, the Cartier dual
of a finite abelian {\sl constant} group scheme $G$.
If $G\simeq \Z/n_1\Z\times\cdots\times\Z/n_r\Z$, then
\begin{equation}\label{directpr}
H\simeq \mu_{n_1}\times\cdots \times\mu_{n_r}
\end{equation} 
where $\mu_{k}=\Spec(\Z[x]/(x^{k}-1))$ denotes the group scheme 
of $k$-th roots of unity over $\Z$. Our goal in this section
is to understand the category of $n$-extensions of $H$ by $\Gm$.
The main result is Theorem \ref{N}.
 
\subsection{} Let us suppose that $n\geq 2$.
\begin{lemma}\label{discrete}
With the above notations and assumptions 
\begin{equation*}
n{\hbox{\small -}}{\rm Ext}^0(H,\Gm)=({\rm id})\ . 
\end{equation*}
\end{lemma} \begin{Proof} We have ($n\geq 2$)
\begin{equation*}
\Hom(\underset{n}{\underbrace {H\otimes\cdots \otimes H}},\Gm)\simeq 
\Hom(\underset{n-1}{\underbrace {H\otimes\cdots \otimes H}},H^D)\ .
\end{equation*} 
Each element of this last group is given by a morphism
$H^{n-1}=H\times_S\cdots _S\times H\to H^D$. 
Since $H^{n-1}$ is connected and $H^D\simeq G$ is constant
any such morphism factors through the identity section; hence this group is trivial. The result now follows from
(\ref{nHom}).  \endproof
\end{Proof}

\begin{Remark} \label{discRe}
{\rm a) Lemma \ref{discrete} shows that for $n\geq 2$ the Picard category 
of $n$-extensions of $H$ by $\Gm$ is ``discrete", i.e there is at most 
one isomorphism between any two objects.

b) As a consequence of (a), any two symmetric $n$-extensions
of $H$ by $\Gm$
which are isomorphic as $n$-extensions are also isomorphic
as {\sl symmetric} $n$-extensions. Indeed, by (a), the diagram 
(\ref{isosymm}) of isomorphisms of $n$-extensions automatically
commutes. In particular, if an $n$-extension of $H$ by $\Gm$
is trivial as an $n$-extension then it is also trivial as a 
symmetric $n$-extension.
}
\end{Remark}

\subsection{}

In what follows we will study the group $n{\hbox{\small -}}{\rm Ext}^1(H,\Gm)$
of isomorphism classes of $n$-extensions of $H=G^D$ by $\Gm$.
We begin by introducing some
notations.

If $C$ is an abelian group and $m\geq 1$ an integer, we will
denote by $C/m$, resp. $_mC$, the cokernel, resp. kernel,
of the map $C\to C$ given by multiplication by $m$.
Set $\zeta_{r}=e^{2\pi i/r}$ 
and for simplicity denote by $C(r)$
the ideal class group $\Cl(\Q(\zeta_r))$ 
of the cyclotomic field $\Q(\zeta_{r})$.
We will identify $(\Z/r\Z)^*$ with the Galois group $\Gal(\Q(\zeta_{r})/\Q)$ 
by sending $a\in (\Z/r\Z)^*$ 
to $\sigma_a$ defined by $\sigma_a(\zeta_{r})=\zeta_{r}^a$.
Now let $p$ be a prime number;
we will denote by $v_p$, resp. $|\ |_p$ the usual
$p$-adic valuation, resp. $p$-adic absolute value. 
Consider the Teichmuller character
$\omega: (\Z/p\Z)^*\ \lo\ \Z_p^*$ characterized by $a= \omega(a)\ {\rm mod}\ p\Z_p$. 
For simplicity, set $\Delta={\rm Gal}(\Q(\zeta_p)/\Q)$.
We will view $\Delta$ as a direct factor of ${\rm Gal}(\Q(\zeta_{p^k})/\Q)$
for any $k\geq 1$.
Suppose $D$ is  a $\Z[\Delta]$-module
which is annihilated by a power of $p$. For $i\in\Z$ we set 
\begin{equation*}
D^{(i)}=\{d\in D\ |\ \sigma_a(d)=\omega^i(a)d, \hbox{\rm \ for all\ }a\in (\Z/p\Z)^*\} .
\end{equation*}
We have
\begin{equation*}
D=\bigoplus_{0\leq i\leq p-2}D^{(i)}.
\end{equation*}
We will consider the 
groups $\Hom(C(p^k), p^{-k}\Z/\Z)$, $k\geq 1$; these are naturally 
${\rm Gal}(\Q(\zeta_{p^k})/\Q)$-modules
via
\begin{equation}\label{action}
(\sigma_a(\phi))(c)=\phi(\sigma^{-1}_a(c))\ \ \hbox{\rm for $\phi: C(p^k)\to p^{-k}\Z/\Z$}.
\end{equation}
Note that the norm $C(p^k)\to C(p^{k-1})$ for the extension $\Q(\zeta_{p^k})/\Q(\zeta_{p^{k-1}})$ 
induces   a homomorphism
\begin{equation*}
N_{k-1}: \Hom(C(p^{k-1}), p^{-(k-1)}\Z/\Z)\to \Hom(C(p^k), p^{-k}\Z/\Z)\ .
\end{equation*}
 
\begin{Definition}
For $n\geq 1$, $m\geq 1$, let $\,\CC(n;p^m)$ be the group of $m$-tuples 
$$
(f_k)_{1\leq k\leq m}\ ;\quad f_k \in \Hom(C(p^k), p^{-k}\Z/\Z)
$$
which satisfy
\begin{eqnarray*}
i)\ \ \ \ \ \ \ \ \ \  \ \sigma_a(f_k ) &=& a^{n-1}f_k,\hbox{\rm \ for all\ }a\in(\Z/p^k)^*,\\
 ii)\ \ \ \ N_{k-1}(f_{k-1} ) &=&p^{n-1}f_k\ .
\end{eqnarray*}
\end{Definition}
\smallskip

\begin{Remark}\label{idealclincl}
{\rm Property (i) implies that
\begin{multline*}
\ \ \ \ \CC(n;p^m)\subset \bigoplus_{1\leq k\leq m}\Hom(C(p^k), p^{-k}\Z/\Z)^{(n-1)}=\\
=\bigoplus_{1\leq k\leq m}\Hom((C(p^k)/p^k)^{(1-n)}, p^{-k}\Z/\Z).\ \ \ \ \ \ \ 
\end{multline*}
In particular, since  $(C(p^k)/p^k)^{(0)}=(0)$ we obtain $\CC(1;p^m)=(0)$.}
\end{Remark}

One of the main results in this section is:
\begin{prop}\label{idealcl}
There is a natural injective homomorphism 
\begin{equation*}
\psi_n: \nExt(\mu_{p^m},\Gm)\xrightarrow{}  \bigoplus_{1\leq k\leq m}\Hom((C(p^k)/p^k)^{(1-n)}, p^{-k}\Z/\Z) 
\end{equation*}
with image the subgroup $\,\CC(n;p^m)$.
\end{prop}
\medskip

Before we consider  the proof we will discuss some
consequences of this result.

\begin{cor}\label{idealcl1}
If $p$ is a regular prime, then 
$\nExt(\mu_{p^m},\Gm)=(0)$.
\end{cor}

Now observe that $\CC(n;p)\simeq\Hom(C(p), \Z/p)^{(n-1)}$.
Hence, for $m=1$ the result  amounts to:

\begin{cor}\label{idealcl2}
There are natural isomorphisms
\begin{equation*}
\nExt(\mu_{p},\Gm)\xrightarrow{\sim} \Hom(C(p), \Z/p)^{(n-1)} \simeq \Hom(
(C(p)/p)^{(1-n)}, \Z/p)\ .
\end{equation*}
\end{cor}
\smallskip

Recall that the $k$-th Bernoulli number $B_k$ is defined by
\begin{equation*}
\frac{t}{e^t-1}=\sum^{\infty}_{k=0}B_k\frac{t^k}{k!}\ .
\end{equation*}
Also, recall that by [Bo], the Quillen ${\rm K}$-groups ${\rm K}_{m}(\Z)$ are finite for
even integers $m\geq 2$. We set $h^+_p=\# {\rm Cl}(\Q(\zeta_p+\zeta_p^{-1}))$ and 
\begin{equation*}
e(n)=\begin{cases}
\ \ 1&,\ \text{if $n=1$},\\
\ \ {\rm numerator\,}({B_{n}}/{n})&,\ \text{if $n$ is even},\\
\displaystyle{\prod_{p,p|h^+_p}}\ord_p(\#{\rm K}_{2n-2}(\Z))&,\ \text{if $n>1$ is odd}.
\end{cases}
\end{equation*}
As we shall now see, Corollary \ref{idealcl2} 
can now be used to obtain:

\begin{thm} \label{N}
For every finite abelian group $G$, the group of isomorphism
classes of $n$-extensions
$\nExt(G^D, \Gm)$
is annihilated by 
$$
 \prod_{p|e(n)}\ord_p(\# G) .
$$
In particular, if $(\#G, e(n))=1$, then $\nExt(G^D, \Gm)=(0)$.
\end{thm}
 
\begin{Proof}
Using (\ref{nExt}) we can see that the group 
$\nExt(G^D, \Gm)$ is annihilated by the order $\#G$
and that it can be written as direct sum  
\begin{equation*}
\bigoplus_{p|\#G}\nExt(G_p^D,\Gm)
\end{equation*}
where $G_p$ is the $p$-Sylow subgroup of $G$. 
The desired result will now follow if we show that  
$p{\not|}e(n)$ implies 
$\nExt(G_p^D,\Gm)=(0)$. Using (\ref{nExt}) again and employing
the long exact cohomology sequence which is obtained by unraveling $G_p^D$
into its ``simple pieces" (each isomorphic to $\mu_p$) we see that it
is enough to show that $p{\not|}e(n)$ implies that 
$\nExt(\mu_p,\Gm)=(0)$. Corollary \ref{idealcl2} then implies that
it suffices to show that when $p{\not|}e(n)$, 
we have $(C(p)/p)^{(1-n)}=(0)$. This now follows
from well-known results on cyclotomic ideal class groups
([W], [Ku], [So]). For the convenience
of the reader we sketch the argument
(we can assume that $p$ is odd). First of all, when $n\geq 2$
is even the result follows directly from Herbrand's theorem 
([W] Theorem 6.17) and Kummer's congruences ([W] Cor. 5.14 and Cor. 5.15). 
To deal with the case that $n$ is odd, we will use the cohomology groups
${\rm H}^i(\Z[1/p], \Z_p(n)):=\underset{\leftarrow m}{\lim}\ 
{\rm H}^2_{\rm et}(\Z[1/p], \mu_{p^m}^{\otimes n})$.
By [Ku] Lemma 1.2 we have
\begin{equation*}
{\rm H}^2(\Z[1/p], \mu_{p}^{\otimes n}) \simeq (C(p)/p)^{1-n},
\end{equation*}
while since ${\rm H}^3(\Z[1/p], \Z_p(n))=(0)$ we can see that
\begin{equation*}
{\rm H}^2(\Z[1/p], \Z_p(n))\otimes_{\Z_p}{\Z/p}\simeq {\rm H}^2(\Z[1/p], \mu_{p}^{\otimes n}). 
\end{equation*}
(See [So], [Ku] for more details.) For $n\geq 2$, there is a surjective Chern character
([DF]; see [So])
\begin{equation*}
{\rm ch}: {\rm K}_{ 2n-2}(\Z)\to {\rm H}^2(\Z[1/p], \Z_p(n))\ .
\end{equation*}
Since $(C(p)/p)^{(0)}=(0)$ and
$p{\not|}h^+_p$ implies that $(C(p)/p)^{(1-n)}=(0)$
for $n>1$ odd, these facts
imply the result.\endproof
\end{Proof}
 
\begin{Remark}\label{low}
{\rm a) We have $B_2=1/6$, $B_4=-1/30$ and  ${\rm K}_4(\Z)$ is trivial 
(see [Ro]; in fact, for our purposes it suffices to know that ${\rm K}_4(\Z)$ has at most $6$-power torsion.
This is somewhat simpler and is shown in Soul\'e's addendum to [LS]).
Hence, we see that Theorem \ref{N} implies that for all finite abelian groups $G$
$$
\nExt(G^D,\Gm)=(0),\text{ for\ } n=1,2,3,4 .
$$

b) Note that [So] gives a doubly exponential
bound on the size of $e(n)$ for $n$ odd. 
However, according to the Kummer-Vandiver conjecture,
$p{\not|}h^+_p$. Assuming this we could replace in the statement of
Theorem \ref{N} $e(n)$ by $e'(n)$ given by $e'(n)=e(n)$ if $n$ is even, $e'(n)=1$
if $n$ is odd. Actually when the prime
divisors of $\#G$ satisfy the Kummer-Vandiver conjecture
(which is true for all primes $<12\cdot 10^6$ by the computations 
of [BCEM])
we have 
$$
\nExt(G^D,\Gm)=(0),\text{ for\ } 1\leq n\leq 11 .
$$
Indeed, the first $e'(n)$ which is not equal to $1$ is $e'(12)=691$.

c) Note that the Quillen-Lichtenbaum conjecture, coupled with the
argument in the proof of Theorem \ref{N} above, implies that for $n>1$ odd we have
$e(n)=2^{a}\cdot\#{\rm K}_{2n-2}(\Z)$ ($a\in \Z$).}
\end{Remark}

Before we continue, we observe that Theorem \ref{N}
together with the results of the previous section
imply:

\begin{thm}\label{Ncube}
Let $G$ be a finite abelian group and $n\geq 1$.
If $\L$ is an invertible sheaf on $H=G^D_{\Spec(\Z)}=\Spec(\Z[G])$ which supports an $n+1$-cubic structure
then 
$$
\L^{\otimes M_n}\simeq \O_H
$$
where
$$
M_n=M_n(G)=\prod_{k=2}^n\prod_{p,p|e(k)}\ord_p(\# G)\ .
$$
\end{thm}

\noindent{\sc Proof of Theorem \ref{Ncube}.} \ \ \
It follows from Theorem \ref{N}, Lemma \ref{61}, Lemma \ref{62} 
and Remark \ref{extra} (in view of Remark \ref{discRe} (b) and the fact that ${\rm Pic}(\Z)=(0)$).
\endproof \smallskip

This combined with Remark \ref{low} (a) gives

\begin{cor}\label{NcubeCor}
Let $G$ be a finite abelian group and $1\leq n\leq 4$.
If $\L$ is an invertible sheaf on $H=G^D_{\Spec(\Z)}=\Spec(\Z[G])$ which supports an $n+1$-cubic structure
then $\L\simeq \O_H$.
\end{cor}
\medskip

\noindent{\sc Proof of Proposition \ref{idealcl}.}  \ Recall that we set $S=\Spec(\Z)$. When $n=1$ 
the Proposition follows from
(\ref{nExt4}), Remark \ref{idealclincl} and the fact that ${\rm H}^1(S,\Z/p^m)=(0)$. 
Now assume that $n\geq 2$ and let $r\geq 1$. Consider the 
homomorphisms 
\begin{equation*}
\delta_i: {\rm H}^1(\mu_{r}^{n-1},\Z/r)\to {\rm H}^1(\mu_{r}^{n},\Z/r),\quad 1\leq i\leq n-1,
\end{equation*}
obtained as $m_i^*-p_i^*-q_i^*$ where 
\begin{eqnarray*}
m_i: \mu_{r}^{n}\to \mu_r^{n-1}\ ;&\quad (x_1,\ldots, x_{n})&\mapsto (x_1,\ldots, x_{i}x_{i+1},\ldots, x_n),\\
p_i: \mu_{r}^{n}\to \mu_r^{n-1}\ ; &\quad (x_1,\ldots, x_{n})&\mapsto (x_1,\ldots, x_{i},x_{i+2},\ldots, x_n),\\
q_i: \mu_{r}^{n}\to \mu_r^{n-1}\ ; &\quad (x_1,\ldots, x_{n})&\mapsto (x_1,\ldots, x_{i-1}, x_{i+1},\ldots, x_n).
\end{eqnarray*}

\begin{lemma}\label{exact}
There is an exact sequence
\begin{equation*}
0\to {(n-1)}{\hbox{\small -}}{\rm Ext}^1(\mu_r, \Z/r)\xrightarrow{t} 
{\rm H}^1(\mu_{r}^{n-1},\Z/r)\xrightarrow{\oplus_{i}\delta_i}\bigoplus_{1\leq i\leq n-1}{\rm H}^1(\mu_{r}^{n},\Z/r).
\end{equation*}
where $t$ is the forgetful map.
\end{lemma}

\begin{Proof}
First suppose that $P$ is a $\Z/r$-torsor over $\mu_r^{n-1}$ whose class $(P)\in {\rm H}^1(\mu_{r}^{n-1},\Z/r)$
is in the kernel of all the homomorphisms $\delta_i$, $1\leq i\leq n-1$. 
We claim that there is a  structure of $n-1$-extension on $P$: By our assumption, there 
are isomorphisms
\begin{equation}\label{ci}
c_i: p_i^*P\cdot q_i^*P\xrightarrow{ \sim} m_i^*P
\end{equation}
of $\Z/r$-torsors over $\mu_r^n$. (We use $\cdot$ to denote the composition of
$\Z/r$-torsors.) Let $e_k: S\to \mu_r^k$ be the identity section and set $P_0=e_{n-1}^*P$;
all $\Z/r$-torsors over $S=\Spec(\Z)$ are trivial, so we can fix a trivialization $\alpha:\Z/r\xrightarrow{\sim} P_0$.
There are natural isomorphisms
$e^*_{n}m_i^*P\simeq P_0$, $e^*_{n}p_i^*P\simeq P_0$, $e^*_nq_i^*P\simeq P_0$.
Hence we obtain an isomorphism of $\Z/r$-torsors over $S$
\begin{equation}\label{ci1}
e_n^*c_i: P_0\cdot P_0\xrightarrow{\sim} P_0.
\end{equation}
Using $\alpha$ we see that $e_n^*c_i$ becomes an isomorphism
\begin{equation*}\label{ci2}
\Z/r \cdot \Z/r \to \Z/r 
\end{equation*}
which is necessarily given by addition in $\Z/r$ followed by translation by an element 
$a_i\in {\rm H}^0(S,\Z/r)=\Z/r$. Modify
$c_i$ by composing with the translation on $m^*_iP$ given by $-a_i$. We can now observe 
that the modified isomorphisms $c_i$ satisfy the conditions of the partial
composition laws (i.e of (\ref{isomulti}))
in the definition of $n-1$-extension: Indeed, this amounts to comparing
certain isomorphisms of $\Z/r$-torsors over $\mu_r^k$, $k=n$, $n+1$, $n+2$.
Note that any two isomorphisms of $\Z/r$-torsors over $\mu_r^k$
differ by the addition of an element in ${\rm H}^0(\mu_n^k,\Z/r)$. Since $\mu_r^k$
is connected we have 
$$
{\rm H}^0(\mu_r^k,\Z/r)\xrightarrow{\overset{e^*_k}{\sim}}{\rm H}^0(S,\Z/r)=\Z/r\ .
$$
and so it is enough to compare the pull-backs of these isomorphisms
via $e_k$. We can now see using the above discussion
that these pull-backs agree.

It remains to show that the map 
${(n-1)}{\hbox{\small -}}{\rm Ext}^1(\mu_r, \Z/r)\to {\rm H}^1(\mu_{r}^{n-1},\Z/r)$ is injective: 
Suppose that a $\Z/r$-torsor $P$ over $\mu_r^{n-1}$
with an $n-1$-extension structure given by the  composition
laws $c_i$ as in (\ref{ci}) is trivial as a $\Z/r$-torsor. Pick a trivialization
$\beta: (\Z/r)_{\mu^{n-1}_r}\xrightarrow{\sim} P$. This induces 
trivializations of $m_i^*P$, $p^*_iP$, $q_i^*P$
under which $c_i$ is identified with 
$$
(\Z/r)_{\mu^{n}_r}\cdot (\Z/r)_{\mu^{n}_r}\to (\Z/r)_{\mu^{n}_r}
$$
given by $(x,y)\mapsto x+y+a'_i$ with $a'_i\in {\rm H}^0(\mu_r^n,\Z/r)=\Z/r$.
The compatibility condition (\ref{compa}) on the $c_i$'s implies that $a'_i$
is independent of the index $i$ (see [SGA7I] VII 2.2
for a similar argument in the case of biextensions); we will denote it by $a'$.
By composing $\beta$ with the translation by $-a'\in\Z/r={\rm H}^0(\mu_r^{n-1},\Z/r)$ 
we can see that $c_i$
becomes the map
\begin{equation*}
(\Z/r)_{\mu^{n}_r}\cdot (\Z/r)_{\mu^{n}_r}\to (\Z/r)_{\mu^{n}_r}
\end{equation*}
given by addition. Therefore, $P$ with the composition laws $c_i$
is isomorphic to the trivial $n-1$-extension.
\endproof
\end{Proof}
\smallskip

We now continue with the proof of Proposition \ref{idealcl}. Lemma \ref{exact} applied to $r=p^m$
and (\ref{nExt3}) implies that, for $n\geq 2$, it is enough to show there is
a natural isomorphism
\begin{equation*}
\CC(n;p^m)\xrightarrow{\sim} {\rm ker}({\rm H}^1(\mu_{p^m}^{n-1},\Z/p^m)
\xrightarrow{\oplus_{i}\delta_i}\bigoplus_{1\leq i\leq n-1}{\rm H}^1(\mu_{p^m}^{n},\Z/p^m)).
\end{equation*}
To identify the kernel above, we will follow a technique used by Mazur in [M] \S 2.
For the convenience of the reader we repeat some of Mazur's arguments.
Suppose that $X$ and $Y$ are any two schemes equipped with $\fp$-valued points
\begin{equation}\label{diagram}
X\leftarrow \Spec(\fp)\rightarrow Y.
\end{equation}
We will use the symbol $X\vee Y$ to refer to any scheme theoretic union 
of $X$ and $Y$ along a subscheme which is a nilpotent extension of $\Spec(\fp)$. 
For $Y=S=\Spec(\Z)$ we set
\begin{equation*}
\tilde {\rm H}^1(X,\Z/p^m)={\rm H}^1(X\vee S,\Z/p^m)
\end{equation*} (fppf or \'etale cohomology).
There is an exact sequence
\begin{equation}\label{MV}
0\to \tilde {\rm H}^1(X,\Z/p^m)\to {\rm H}^1(X,\Z/p^m)\to {\rm H}^1(\Spec(\fp),\Z/p^m)
\end{equation}
obtained using the Mayer-Vietoris exact sequence for \'etale cohomology,
the fact that $\Spec(\fp)$ is connected and that ${\rm H}^1(S,\Z/p^m)=(0)$.
Hence, $\tilde {\rm H}^1(X,\Z/p^m)$ is independent of the exact scheme theoretic union
of $X$ and $S$ used in the definition. A similar calculation, shows that for any
diagram as in (\ref{diagram}), we have
\begin{equation}\label{add}
\tilde {\rm H}^1(X\vee Y,\Z/p^m)=\tilde {\rm H}^1(X,\Z/p^m)\oplus \tilde {\rm H}^1(Y,\Z/p^m) .
\end{equation}

Now set $T_{p^k}=\Spec(\Z[\zeta_{p^k}])$, $1\leq k\leq m$;
this is a closed subscheme of $\mu_{p^k}$. Class-field theory
gives a natural isomorphism
\begin{equation*}
{\rm H}^1(T_{p^k}, \Z/p^m)\xrightarrow{\sim}\Hom(C(p^k),\Z/p^m).
\end{equation*}
Since the unique prime ideal of $\Z[\zeta_{p^k}]$ 
that lies above $(p)$ is principal, by the exact sequence (\ref{MV}) and the definition of the Artin map, we have
\begin{equation}\label{cft}
\ti{\rm H}^1(T_{p^k}, \Z/p^m)={\rm H}^1(T_{p^k}, \Z/p^m)\xrightarrow{\sim} \Hom(C(p^k),\Z/p^m).
\end{equation}
Now observe that we have canonical identifications
\begin{equation*}
(\Z/p^k)^*={\rm Aut}(\mu_{p^k})={\rm Aut}(T_{p^k})
\end{equation*}
where $a\in(\Z/p^k)^*$ acts via the operation ``raising to the $a$-th power"
on $\mu_{p^k}$. The isomorphism (\ref{cft}) is compatible with 
the action of $(\Z/p^k)^*$ by functoriality of cohomology on the
one side and by (\ref{action}) on the other. If $\pi_k: T_{p^k}\to T_{p^{k-1}}$
is the natural projection, there is a commutative diagram
\begin{equation}\label{norm2}
\begin{CD}
\ti{\rm H}^1(T_{p^{k-1}}, \Z/p^m)@={\rm H}^1(T_{p^{k-1}}, \Z/p^m)@>\sim>>\Hom(C(p^{k-1}),\Z/p^m)\\
@V\pi^*_kVV @V\pi^*_kVV @VVN_{k-1}V\\
\ti{\rm H}^1(T_{p^{k}}, \Z/p^m)@={\rm H}^1(T_{p^{k}}, \Z/p^m)@>\sim>>\Hom(C(p^{k}),\Z/p^m)\\
\end{CD}
\end{equation}
with $N_{k-1}$ induced by the norm.
Now notice that $\mu_{p^m}^s$
for any $s\geq 1$, can be obtained as a wedge (in the sense of $\vee$ defined above) 
of several copies of $T_{p^k}$, $1\leq k\leq m$, with $S$. More precisely, 
$\mu_{p^m}^s$ is the wedge of $S$ with
\begin{equation*}
\bigvee_{1\leq k\leq m}\, \big[\bigvee_{ (a_1;\cdots;a_s)\in {\bf P}^{s-1}(\Z/p^k)}T_{p^k}\big].
\end{equation*}
 Using (\ref{add}) we can deduce that
\begin{equation}\label{dirsum}
{\rm H}^1(\mu_{p^m}^s,\Z/p^m)=\bigoplus_{1\leq k\leq m}
\bigoplus_{(a_1;\cdots;a_s)\in {\bf P}^{s-1}(\Z/p^k)} \ti{\rm H}^1(T_{p^k},\Z/p^m).
\end{equation}
Notice that an element $(a_i)=(a_1,\ldots, a_s)\in (p^{-k}\Z/\Z)^s$
defines a group scheme homomorphism 
\begin{equation*}
(a_i): \mu_{p^k}\to \mu^s_{p^m}\ ;\quad x\mapsto (x^{p^ka_1},\ldots, x^{p^ka_t})
\end{equation*}
and a scheme morphism
\begin{equation*}
(a_i): T_{p^k}\subset \mu_{p^k}\to \mu^s_{p^m}.
\end{equation*}
Since $\pi_k(x)=x^p$ we have a commutative diagram
\begin{equation}\label{cd}
\begin{CD}
T_{p^k}@>(pa_i)>>\mu^s_{p^m}\\
@V\pi_kVV  @|\\
T_{p^{k-1}}@>(pa_i)>>\mu^s_{p^m},\\  
\end{CD}
\end{equation}
where in the first, resp. second line, $(pa_i)$ is considered
as an element of $(p^{-k}\Z/\Z)^s$, resp. of $(p^{-(k-1)}\Z/\Z)^s$.

Set 
$
U^s_{p^k}=(p^{-k}\Z/\Z)^s-(p^{-k+1}\Z/\Z)^s
$;
if $(a_i)$ is in $U^s_{p^k}$ then the corresponding morphism
is a closed immersion. Now consider the group of maps
\begin{equation*}
{\rm Maps}_{(\Z/p^k)^*}(U^s_{p^k}, \ti{\rm H}^1(T_{p^k},\Z/p^m))
\end{equation*}
which are compatible with the natural action of $(\Z/p^k)^*$
on domain and range. Note that ${\bf P}^{s-1}(\Z/p^k)\simeq (\Z/p^k)^*\backslash U^s_{p^k}$.
We can define a homomorphism
\begin{equation}\label{map}
{\rm H}^1(\mu_{p^m}^s,\Z/p^m)\to \bigoplus_{1\leq k\leq m}
{\rm Maps}_{(\Z/p^k)^*}(U^s_{p^k}, \ti{\rm H}^1(T_{p^k},\Z/p^m)) \end{equation}
by sending $h\in {\rm H}^1(\mu_{p^m}^s,\Z/p^m)$ to $(a_i)\mapsto (a_i)^*h$.
Using (\ref{dirsum}) we see that (\ref{map}) is an isomorphism.
We can also consider the map
\begin{equation}\label{map2}
{\rm H}^1(\mu_{p^m}^s,\Z/p^m)\to \bigoplus_{1\leq k\leq m}
{\rm Maps}_{(\Z/p^k)^*}((p^{-k}\Z/\Z)^s, \ti{\rm H}^1(T_{p^k},\Z/p^m))
\end{equation}
given by the same rule as the one above. The map (\ref{map2}) is injective and
using (\ref{cd}) we can see that its image
is the subgroup of all elements $(\phi_k)_{1\leq k\leq m}$ which
are such that $\phi_k((pa_i))=\pi_k^*\phi_{k-1}((pa_i))$.
Let us denote this subgroup by $F(s;p^m)$.
By applying the above to $s=n-1$, $n$ we can conclude that there are 
commutative diagrams
\begin{equation}\label{cd2}
\begin{CD}
{\rm H}^1(\mu_{p^m}^{n-1},\Z/p^m) @>\delta_i>>  {\rm H}^1(\mu_{p^m}^n,\Z/p^m)\\
@V(\ref{map2})_{n-1}VV  @VV(\ref{map2})_{n}V\\
F(n-1;p^m) @>\delta'_i>> F(n;p^m) \\
\end{CD}
\end{equation}
\smallskip

\noindent where the $k$-th component of $\delta_i'((\phi_k)_{1\leq k\leq m})$
is the map given by
\begin{equation*}
(\ldots,a_i,a_{i+1},\ldots)\mapsto \phi_k(\ldots,a_i+a_{i+1},\ldots)-
\phi_k(\ldots,a_i,\ldots)-
\phi_k(\ldots,a_{i+1},\ldots).
\end{equation*}

It now follows that the kernel of $\oplus_i\delta_i$ is isomorphic 
to the group of $m$-tuples $(\phi_k)_{1\leq k\leq m}$ 
of  multilinear maps 

\begin{equation*}
\phi_k: (p^{-k}\Z/\Z)^{n-1}\to \ti{\rm H}^1(T_{p^k},\Z/p^m)=\Hom(C(p^k),\Z/p^m)
\end{equation*}

which satisfy
\begin{eqnarray*}
i)\ \ \ \ \phi_k(ax_1,\ldots, ax_{n-1})&=& \sigma_a(\phi_k(x_1,\ldots, x_{n-1})),\hbox{\rm \ for all\ }a\in(\Z/p^k)^*,\\
\ \ \ ii)\ \ \ \ \phi_k(px_1,\ldots, px_{n-1})&=&\pi_k^*(\phi_{k-1}(px_1,\ldots, px_{n-1}))\ .
\end{eqnarray*}

Note  that a multilinear map $\phi_k: (p^{-k}\Z/\Z)^{n-1}\to  \Hom(C(p^k),\Z/p^m)$
has image which is contained in ${}_{p^k}\Hom(C(p^k),\Z/p^m)\simeq\Hom(C(p^k),p^{-k}\Z/\Z)$; such a map
is uniquely determined by $f_k:=\phi_k(p^{-k},\ldots, p^{-k})\in \Hom(C(p^k),p^{-k}\Z/\Z)$.
Using (\ref{norm2}) and the multilinearity we see that conditions (i) and (ii) above translate to
\begin{eqnarray*}
i)'\ \ \ \ \ \ \ \ \ \  \ \sigma_a(f_k ) &=& a^{n-1}f_k,\hbox{\rm \ for all\ }a\in(\Z/p^k)^*,\\
 ii)'\ \ \ \ N_{k-1}(f_{k-1} ) &=&p^{n-1}f_k\ .
\end{eqnarray*}

The proof of Proposition \ref{idealcl} now follows.\endproof
\medskip

\begin{Remark}\label{diagonal}
{\rm We can see from the proof that the injective homomorphism 
\begin{equation}
\psi_n: {n}{\hbox{\small -}}{\rm Ext}^1(\mu_{p^m}, \Gm)\hookrightarrow  \bigoplus_{1\leq k\leq m}\Hom((C(p^k)/p^k)^{(1-n)}, p^{-k}\Z/\Z) 
\end{equation}
is obtained as follows: Consider the homomorphism
\begin{equation}\label{compo1}
\psi'_n: {n}{\hbox{\small -}}{\rm Ext}^1(\mu_{p^m}, \Gm)\to {\rm H}^1(\mu_{p^m},\Z/p^m),
\end{equation}
defined as the composition
\begin{multline*}
\ \ \ \ \ {n}{\hbox{\small -}}{\rm Ext}^1(\mu_{p^m}, \Gm)\to {(n-1)}{\hbox{\small -}}{\rm Ext}^1(\mu_{p^m}, \Z/p^m)\xrightarrow{t}\\
\xrightarrow{t}{\rm H}^1(\mu_{p^m}^{n-1},\Z/p^m)\xrightarrow{\Delta_{n-1}^*}
{\rm H}^1(\mu_{p^m},\Z/p^m)\ \ \ \ \ \ \ \ \ \ \ \ 
\end{multline*}
where the first arrow is  the inverse of (\ref{nExt3}), $t$ is the forgetful
map and $\Delta_{n-1}^*$ is the pull-back along the diagonal
$\Delta_{n-1}:\mu_{p^m}\to \mu_{p^m}^{n-1}$ . Then $\psi_n$ is given by 
the composition of $\psi'_n$ with the isomorphism
\begin{equation*}
{\rm H}^1(\mu_{p^m},\Z/p^m)\xrightarrow{\sim} \bigoplus_{1\leq k\leq m}\Hom(C(p^k) , p^{-m}\Z/\Z)
\end{equation*} 
obtained by (\ref{cft}) and (\ref{dirsum}). Indeed, the maps $\phi_k$
in the proof of the Proposition are determined by their image on the ``diagonal" elements
$(p^{-k},\ldots, p^{-k})$.
}
\end{Remark}
\medskip

\section{Reflection homomorphisms} \label{reflsect}

In the next few paragraphs, we elaborate on the constructions of the previous section.
We continue with the same assumptions and notations.  
In particular, we again write $T_{p^k}=\Spec(\Z[\zeta_{p^k}])$ which we think of as
a closed subscheme of $\mu_{p^k}$.
We will denote by $\widetilde {\mu_{p^m}}=\bigsqcup_{0\leq k\leq m}T_{p^k}$ the normalization 
of the scheme $\mu_{p^m}$ and by $\nu: \widetilde {\mu_{p^m}}\to \mu_{p^m}$ the natural
projection map. Our main goal is to express the composition
\begin{equation}\label{diaghomo}
n{\hbox{\small -}}{\rm Ext}^1(\mu_{p^m},\Gm)  \xrightarrow{t} 
{\rm Pic}(\mu_{p^m}^n)  \xrightarrow{ \Delta_n^* } {\rm Pic}(\mu_{p^m})\xrightarrow{ \nu^* }
 {\rm Pic}(\widetilde{\mu_{p^m}})
\end{equation}
in terms of the classical ``reflection homomorphisms" (see below).
We do this in Corollary \ref{Reflcor}. We can then deduce some 
additional results on the pull-back $\nu^*\L$ of an 
invertible sheaf $\L$  on $\mu_{p^m}$ with hypercubic structure.
\medskip

\subsection{ }
Consider the homomorphism (\ref{Lchimap}) described in \ref{Lsection}
\begin{equation*}
R_k: {\rm H}^1(T_{p^k},\Z/p^k)\to   {}_{p^k}\Pic(T_{p^k})\ ;\quad Q\mapsto \L^Q_{f(\chi_0)} 
\end{equation*}
for $G=\Z/p^k$, $S'=T'=T=T_{p^k}$, $\chi_0: T\hookrightarrow \mu_{p^k}$
the natural closed immersion which corresponds to the character 
$\chi_0: \Z/p^k\to \Z[\zeta_{p^k}]^*$,
$\chi_0(1)= \zeta_{p^k}$, and $f(\chi_0): T'=T \to {\mu_{p^k}}_{T}$ the morphism
$$
f(\chi_0):T \xrightarrow{\Delta} T \times_S T \xrightarrow{(\chi_0, {\rm id})} {\mu_{p^k}}_{T}=\mu_{p^k}\times T.
$$
Using (\ref{cft}) (class field theory) and $C(p^k)=\Pic(T_{p^k})$ 
we see that this amounts to a homomorphism
\begin{equation*}
R_k: \Hom(C(p^k),p^{-k}\Z/\Z)\to {}_{p^k}C(p^k).
\end{equation*}
If $Q\to T_{p^k}$ is a $\Z/p^k$-torsor
there is an unramified Galois extension $N$ of $\Q(\zeta_{p^k})$ with Galois group $\Z/p^k$
and ring of integers $\O_N$ such that the $\Z/p^k$-torsor $Q$ is $Q=\Spec(\O_N)$.
Lemma \ref{Lchi} implies that $\L^Q_{f(\chi_0)}$
is isomorphic to the invertible sheaf which corresponds to the locally free rank $1$ 
$\Z[\zeta_{p^k}]$-module 
\begin{equation}
L^Q_{\chi_0}:=\{\xi\in \O_N\,|\,\sigma_a(\xi)=\chi^{-1}_0(a)\xi=\zeta_{p^k}^{-a}\xi,   \text{ for all }a\in \Z/p^k\} .
\end{equation}
Notice that $N/\Q(\zeta_{p^k})$ is a Kummer extension. Therefore,
it can be obtained by adjoining the $p^k$-th root of 
an element $b\in \Q(\zeta_{p^k})^*$: $N=\Q(\zeta_{p^k})(\sqrt[p^k]{b})$.
We can arrange so that the element $\sqrt[p^k]{b}$ gives a generic section of $\L^Q_{f(\chi_0)}$; 
the corresponding
divisor of $\L^Q_{f(\chi_0)}$ is given by a fractional ideal $I$ 
of $\Q(\zeta_{p^k})$ such that
$I^{p^k}=(b)$. The class of $\L^Q_{f(\chi_0)}$ corresponds to the class
$(I)$ under the isomorphism $\Pic(T_{p^k})\simeq C(p^k)$.
Using this we can see that $R_k$ coincides with the
classical {\sl ``reflection homomorphism"}
\begin{equation*}
\Hom(C(p^k),p^{-k}\Z/\Z)\to {}_{p^k}C(p^k)
\end{equation*}
(see for example [W] \S 10.2 for the case $k=1$; actually the reflection homomorphism defined there is  the negative of the one above):

Now observe that the definition of $\L^Q_{f(\chi_0)}$ implies
\begin{equation*}
a^*(\L^Q_{f(\chi_0)})\simeq \L^{a^*Q}_{f(\chi_0^a)} ,\quad a\in (\Z/p^k)^*,
\end{equation*}
where $a^*$ denotes the pull-back by the Galois automorphism 
$a:T_{p^k}\to T_{p^k}$. Using (\ref{Lmulti}) we see that this gives
\begin{equation}\label{Galoisreflect}
a^*(\L^{Q}_{f(\chi_0)})\simeq  \L^{a^*Q}_{f(\chi_0^a)}\simeq (\L^{a^*Q}_{f(\chi_0)})^{\otimes a}.
\end{equation}

The isomorphism (\ref{Galoisreflect}) now implies that 
$R_k$ ``reflects" between odd and even eigenspaces,
in fact it decomposes into a direct sum of
\begin{multline}
R_k^{(n)}: \Hom((C(p^k)/p^k)^{(1-n)},p^{-k}\Z/\Z)=\\=\Hom(C(p^k),p^{-k}\Z/\Z)^{(n-1)}\to ({}_{p^k}C(p^k))^{(n)}.\ \ \ \ \ \ \ \ \ 
\end{multline}
for $0\leq n\leq p-2$ (cf. [W] \S 10.2).

\subsection{ } For notational simplicity, we set $r=p^m$.
Recall $\nu: \widetilde{\mu_r}\to \mu_r$ is the normalization morphism.
Let us consider the homomorphism
\begin{equation}\label{compo2}
{\rm H}^1(\mu_{r}, \Z/r)\to \Pic(\mu_{r})\xrightarrow{\nu^*} \Pic(\widetilde{\mu_r})=\bigoplus_{1\leq k\leq m}C(p^k)
\end{equation}
where the first arrow is the composition 
\begin{equation}\label{compo3} {\rm H}^1(\mu_{r}, \Z/r)\xrightarrow{(\ref{nExt4})}\Ext^1({\mu_{r}}_{\mu_{r}},\Gm_{\mu_r})\xrightarrow{t} 
\Pic(\mu_{r}\times \mu_{r})\xrightarrow{\Delta^*_2}\Pic(\mu_{r}).
\end{equation}

Recall that  (\ref{cft}) and (\ref{dirsum}) give an isomorphism
\begin{equation}\label{compos4}
{\rm H}^1(\mu_{r}, \Z/r)\xrightarrow{\sim} \bigoplus_{1\leq k\leq m}\Hom(C(p^k),p^{-m}\Z/\Z).
\end{equation}
Now let us restrict the map (\ref{compo2}) to the subgroup of 
${\rm H}^1(\mu_{r}, \Z/r)$ that corresponds to 
$$
\bigoplus_{1\leq k\leq m}\Hom(C(p^k),p^{-k}\Z/\Z)
$$
under (\ref{compos4}). We obtain a homomorphism
\begin{equation}\label{restr}
R: \bigoplus_{1\leq k\leq m}\Hom(C(p^k),p^{-k}\Z/\Z)\xrightarrow{\ } \bigoplus_{1\leq k\leq m}C(p^k).
\end{equation}  
By unraveling the definition of $R$ we see
that the description
of ``reflection homomorphisms" in the above paragraph implies that
\begin{prop}\label{Refldirect}
The homomorphism $R$ is a direct sum $R=\oplus_{1\leq k\leq m}R_k$
with
\begin{equation*}
R_k: \Hom(C(p^k),p^{-k}\Z/\Z)\to {}_{p^k}C(p^k)
\end{equation*} the ``reflection homomorphism" as defined above.\endproof
\end{prop}
\smallskip

We now obtain:

\begin{cor} \label{Reflcor}
There is a commutative diagram
\begin{equation*}
\begin{matrix}
n{\hbox{\small -}}{\rm Ext}^1(\mu_{p^m},\Gm)\ \ \  \ \xrightarrow{t}\ \  
{\rm Pic}(\mu_{p^m}^n)&  \ \xrightarrow{ \Delta_n^*\ }\  \ {\rm Pic}(\mu_{p^m})&\xrightarrow{ \nu^*\ }
\ {\rm Pic}(\widetilde{\mu_{p^m}})\\
&&\\
\psi_n\Big\downarrow\ \ \ \ \ \ \ \ \ \ \ \ \ \ \ \ &&\ \ \cup\\
&&\\
\displaystyle{\bigoplus_{1\leq k\leq m}\Hom((C(p^k)/p^k)^{(1-n)},p^{-k}\Z/\Z)}& \xrightarrow{\ \displaystyle{\oplus_{1\leq k\leq m}R^{(n)}_k}\ }\ &
\displaystyle{\bigoplus_{1\leq k\leq m}({}_{p^k}C(p^k)})^{(n)}.\\
\end{matrix}
\end{equation*}
\end{cor}

\begin{Proof}
Recall that by Remark \ref{diagonal}, the homomorphism $\psi_n$ is given by the composition of 
$\psi_n': n{\hbox{\small -}}{\rm Ext}^1(\mu_{p^m},\Gm)\to {\rm H}^1(\mu_{p^m}, \Z/p^m)$
with the isomorphism (\ref{compos4}). The result follows now from the definitions
of the homomorphisms $R$ and $\psi_n'$, Proposition \ref{Refldirect},
and the commutative diagram (\ref{bigcd}) for $H=\mu_{p^m}$:
Indeed, we can observe that the composite homomorphism (\ref{compo3}) essentially gives
a half of the commutative diagram (\ref{bigcd})
for $H=\mu_{p^m}$.  
\endproof
\end{Proof}

\begin{Remark} \label{Refltriv}
{\rm If the prime $p$ satisfies the Kummer-Vandiver conjecture, then 
the reflection maps $R^{(n)}_k$ are all trivial; indeed, either $n$ or $1-n$ is even
and so either $({}_{p^k}C(p^k))^{(n)}$ or $(C(p^k)/p^k)^{(1-n)}$ is trivial
([W] Cor. 10.6).
Then the composition $\nu^*\cdot \Delta^*_n\cdot t$ along the first row
of the above diagram is also the trivial homomorphism.}
\end{Remark}
\medskip

\subsection{ }
We now combine the above to obtain
an additional result about invertible sheaves with hypercubic structures
over $H=\Spec(\Z[G])$, for $G$ any finite abelian group. For an integer $u\geq 1$ set 
\begin{equation*}
e'(u)=\begin{cases}
\  1 &, \ \text{if $u$ is odd},\\
\  {\rm numerator\,}({B_{u}}/{u})&, \ \text{if $u$ is even}.
\end{cases} 
\end{equation*}
(cf. Remark \ref{low} (b)). Set
\begin{equation*}
M'_n(G)=\prod_{ u=1}^n\prod_{p|e'(u)}\ord_p(\# G).
\end{equation*}

\begin{thm}\label{KVtrivial}
Assume that all the prime divisors of $\#G$ satisfy the
Kummer-Vandiver conjecture. We denote
by $\nu: \widetilde {H}\to H$ 
the normalization morphism. Suppose that $\L$ is an invertible sheaf on 
$H$ which supports an $n+1$-cubic structure $\xi$ and set
$C=\text{\rm GCD}(M'_n(G), n!!)$.
Then 
$
\nu^*\L^{\otimes C}\simeq\O_{\widetilde H}.
$
In particular, if in addition all the prime 
divisors of $\#G$ are $\geq n+1$, then $\nu^*\L\simeq\O_{\widetilde H}$.
\end{thm}

\begin{Proof} Suppose that $G=G_{p_1}\times\cdots \times G_{p_k}$ is the decomposition
of $G$ into its $p$-Sylow subgroups. Set $H_{p_j}=(G_{p_j})^D=\Spec(\Z[G_{p_j}])$ and let ${\rm pr}_j: H\to H_{p_j}$
be the natural projection.
If all the prime divisors $p_j$, $1\leq j\leq k$, of $\#G$ satisfy the Kummer-Vandiver conjecture
we have $\L^{\otimes M'_{n}(G)}\simeq \O_{H}$
by Theorem \ref{Ncube} (cf. Remark \ref{low} (b)). On the other hand, Corollary \ref{flat2}
gives the ``Taylor expansion"
\begin{equation*}
\nu^*\L^{\otimes n!!}\simeq \bigotimes_{i=0}^{n-1}
(\nu^*\Delta^*_{n-i}E(\L^{(i)},\xi^{(i)}))^{\otimes (-1)^i(n-i-1)!!}\otimes \nu^*0^*\L^{\otimes n!!}.
\end{equation*}
Here, the invertible sheaf $E(\L^{(i)},\xi^{(i)})$ carries the structure of an $(n-i)$-extension
of $H$ by $\Gm$. Notice that, since ${\rm Pic}(\Z)=(0)$, $0^*\L$ 
is trivial. Our goal is to show that the invertible sheaves 
$\nu^*\Delta^*_{n-i}E(\L^{(i)},\xi^{(i)})$, $0\leq i\leq n-1$, are also trivial.
This would imply that we have $\nu^*\L^{\otimes n!!}\simeq\O_{\widetilde H}$
from which, given the above discussion, the result follows.
Observe that (\ref{nExt}) implies that there is an isomorphism
of multiextensions
\begin{equation}\label{decoSylow}
E(\L^{(i)},\xi^{(i)})\simeq \bigotimes_{j=1}^k ({\rm pr}_j\times\cdots\times{\rm pr}_j)^* (E^i_j)\ ,
\end{equation}
where $E^i_j$ is an $(n-i)$-extension of $H_{p_j}$ by $\Gm$. Using (\ref{decoSylow})
we see that it is enough to prove that the invertible sheaves $\nu^*\Delta^*_{n-i}(E^i_j)$ 
are trivial, where now $\nu$ and $\Delta_{n-i}$ are the normalization and diagonal 
morphisms for the group scheme $H_{p_j}$. In fact, since the normalization $\widetilde H_{p_j}$
is the disjoint union of components corresponding to characters of $G_{p_j}$ and these 
factor through prime power order cyclic quotients we can see that 
we can reduce to the case of a prime power order cyclic group.
More precisely, it is enough to show the following statement:
If $E$ is an $(n-i)$-extension of $\mu_{p^m}$ by $\Gm$, then the invertible sheaf
$\nu^*\Delta^*_{n-i}(E)$ is trivial.
Corollary \ref{Reflcor} applied to $n-i$ implies that the composition
\begin{equation}\label{above89}
(n-i){\hbox{\small -}}{\rm Ext}^1(\mu_{p^m},\Gm)\xrightarrow{t} {\rm Pic}(\mu_{p^m}^{n-i}) 
\xrightarrow {\nu^*\Delta^*_{n-i}}
{\rm Pic}(\widetilde{\mu_{p^m}})
\end{equation} 
factors through the reflection homomorphisms. Hence, if $p$ satisfies the Kummer-Vandiver conjecture
then the homomorphism (\ref{above89}) is trivial
(Remark \ref{Refltriv}).
Therefore,  the invertible sheaves 
$\nu^*\Delta^*_{n-i}(E)$  are  trivial.
The result now follows.\endproof
\end{Proof}
\medskip

\subsection{} \label{8d} Suppose that $G=\Z/p$ and that $\L$ is an invertible sheaf on 
$\mu_{p}$ which supports an $d+2$-cubic structure $\xi$ with $p> d+1$. 
As above, Corollary \ref{flat2} gives
\begin{equation*}
\nu^*\L^{\otimes (d+1)!!}\simeq \bigotimes_{i=0}^{d}\nu^*(\Delta^*_{d+1-i}E(\L^{(i)},\xi^{(i)}))^{\otimes (-1)^i(d-i)!!}.
\end{equation*}
Since, by Theorem \ref{Ncube}, the invertible sheaf $\L$ is $p$-power torsion and ${\rm GCD}(p, (d+1)!!)=1$,
we can write
\begin{equation*}
\nu^*\L \simeq \bigotimes_{i=0}^{d}\nu^* \Delta^*_{d+1-i}E'_{d+1-i},
\end{equation*}
where $E'_{d+1-i}$ is an invertible sheaf on $\mu_p^{d+1-i}$ with a $(d+1-i)$-extension structure.
Let us denote by $t_{j}(\L,\xi)$ the image of $E'_j$ under the homomorphism
\begin{equation*}
\psi_j: j{\hbox{\small -}}{\rm Ext}^1(\mu_{p},\Gm)\to \Hom((C(p)/p)^{(1-j)},p^{-1}\Z/\Z)
\end{equation*}
of Proposition \ref{idealcl} (In this case, $\psi_j$ is an isomorphism; cf. Corollary \ref{idealcl2}).
By [Ri], the pull back $\nu^*: \Pic(\mu_p)\to \Pic(\tilde \mu_p)={\rm Cl}(\Q(\zeta_p))$ is an isomorphism
and we can use it to identify these class groups.
Therefore, Corollary \ref{Reflcor} and the above equality now implies that we can write
\begin{equation}\label{84}
[\L]=\sum_{j=1}^{d+1} R^{(j)}( t_{j}(\L,\xi))
\end{equation}
in $\Pic(\mu_p)={\rm Cl}(\Q(\zeta_p))$ with $R^{(j)}=R^{(j)}_1: \Hom((C(p)/p)^{(1-j)},p^{-1}\Z/\Z)\to ({}_pC(p))^{(j)}$
the reflection homomorphism.
\bigskip

\section{The hypercubic structure on the determinant of cohomology}\label{hypercubedet}

In this section, we explain some of the results of F. Ducrot [Du] and
we show how the work in [Du] can be used to deduce the main result in this 
section, Theorem \ref{detC}.

For every non-empty finite set $I$ we denote by 
$C(I)$ the set of all subsets of $I$. Suppose that $\PP$ is a s.c.  Picard category
(\S \ref{picard}).
By definition, an {\sl $I$-cube} in $\PP$ is a family $K=(K_v)_{v\in C(I)}$ 
of objects of $\PP$, indexed by the set $C(I)$ of all subsets  of $I$. 
We denote by $I\mi Cube(\PP)$ the category of $I$-cubes
in $\PP$. 

If $n\geq 1$, then an {\sl $n$-cube} in $\PP$ is by definition 
an $I$-cube in $\PP$ for some $I$ with $\# I=n$. A morphism between
two $n$-cubes $K=(K_v)_{v\in C(I)}$, $K'=(K'_{w})_{w\in C(I')}$
is a pair $\Phi=(\varphi,\phi)$ where $\varphi: I\xrightarrow{\sim} I'$
is bijective and $\phi: K\xrightarrow {\sim}\varphi^*(K')$ is an (iso)morphism
of $I$-cubes. The category of all $n$-cubes 
in $\PP$ will be denoted by $n\mi Cube(\PP)$. 

This ``cubic" terminology is motivated by 
the fact that if $I=\{1,2,\ldots, n\}$ we may think of $C(I)$ as the 
vertices of the standard $n$-dimensional cube $C_n$ in ${\bf R}^n$
by sending a subset $v\subset I$ to the point $(v_i)_{1\leq i\leq n}$
with $v_i=1$ if  $i\in v$, $v_i=0$ if $i\not\in v$. Hence, after choosing 
an order for $I$, we can visualize the objects $K_v$ of an $I$-cube $K$
as being placed on the
vertices of the $n$-dimensional cube $C_n$. 

For $v\in C(I)$, set $s(v)=(-1)^{\#I-\# v}$.
If $K$ is an $n$-cube in $\P$ we now set
\begin{equation}\label{sigma} \Sigma(K)=\sum_{v\in C(I)}(-1)^{s(v)}K_v\ 
\end{equation}  
(cf. Lemma \ref{construct}). We can see that an isomorphism $\Phi: K\to K'$ between $n$-cubes induces
an isomorphism
\begin{equation}\label{sigmaphi}
\Sigma(\Phi): \Sigma(K)\xrightarrow{\sim} \Sigma(K')
\end{equation}
in $\PP$. In fact, we obtain a functor $\Sigma: n\mi Cube(\PP)\to \PP$.

Suppose that $K$ is an $I$-cube and let $J_0$, $J$ be disjoint subsets of $I$ with $\#J=q$.
We call the $J$-cube $(K_{J_0\cup J'})_{J'\subset J}$ a $q$-face of $K$.
If $i\in I$, we will denote by ${\rm Front}_i(K)$, ${\rm Back}_i(K)$ the
$I-\{i\}$-cubes given by $J_0=\{i\}$, $J_0=\emptyset$. We may think of 
${\rm Front}_i(K)$, ${\rm Back}_i(K)$ as the two ``opposite" $n-1$-faces of $K$ obtained by 
restricting the $i$-th coordinate to be equal to $1$, resp. $0$.
There is a canonical isomorphism (cf. Lemma \ref{construct})
\begin{equation}\label{sigmadec}
\Sigma({\rm Front}_i(K))+(-\Sigma({\rm Back}_i(K)))\xrightarrow{\sim} \Sigma(K).
\end{equation}  
Conversely, if $I=J\cup \{i\}$ and $A$, $B$ are two $J$-cubes, we will
denote by  
$A\buildrel{i}\over{\hbox{\rm ---}}B$ the $I$-cube whose $i$-th back face and 
$i$-th front face are respectively $A$ and $B$. We obtain a canonical isomorphism
\begin{equation}\label{sigmadec2}
\Sigma(B)+(-\Sigma(A))\xrightarrow{\sim} \Sigma(A\buildrel{i}\over{\hbox{\rm ---}}B).
\end{equation} 

\begin{Definition}
{\rm By definition, a {\sl decorated $n$-cube} in $\P$ is 
an $n$-cube $K$ together with isomorphisms 
(``$2$-face trivializations") 
\begin{equation}\label{face}
m_F: \uno\xrightarrow{\sim} \Sigma(F),
\end{equation}
for each $2$-face $F$ of $K$, which satisfy the following condition:
Suppose that $L$ is a $3$-face of $K$ which corresponds to the subsets $J_0$ and $J$ as above.
For each pair of indices $i, j\in J$ the $2$-face trivializations for
${\rm Front}_i(L)$, ${\rm Back}_i(L)$, ${\rm Front}_{j}(L)$, ${\rm Back}_{j}(L)$
should be compatible with the canonical isomorphism given by (\ref{sigmadec})
\begin{equation}
\Sigma({\rm Front}_i(L))+(-\Sigma({\rm Back}_i(L)))\xrightarrow{\sim} \Sigma(L)\xrightarrow{\sim} 
\Sigma({\rm Front}_{j}(L))+(-\Sigma({\rm Back}_{j}(L))).
\end{equation}}
\end{Definition}

\subsection{}
 
Now let $\Sfr$ be a site and suppose that $a: \P\to {\Sfr}$ is a s.c. Picard $\Sfr$-stack.
For our applications, $\Sfr=S_{\rm fppf}$.
A (decorated) $n$-cube $K$ in $\P$ is a (decorated) $n$-cube 
in the fiber category $\P_T$ for some $T\in {\rm Ob}(\Sfr)$;
in particular, all the objects $K_v$, $v\in C(I)$, satisfy $a(K_v)=T$. 

The $n$-cubes in $\P$ form naturally an $\Sfr$-stack $n\mi Cube(\P)$.
There is also an $\Sfr$-stack $n\mi Cube_d(\P)$ whose objects are decorated $n$-cubes in $\P$,
and morphisms are given by  morphisms $\Phi=(\varphi,\phi): K\to K'$ 
in $n\mi Cube(\P)$ which are compatible with the isomorphisms (\ref{face}) in the following sense:
If $F'\subset K'$ is a $2$-face, then let $F$ be the corresponding (via $\varphi$) $2$-face of $K$.
Let $a(\Phi)=\tau: T\to T'$.
We ask that
the following diagram commutes:
\begin{equation}
\begin{CD}
\uno_{T}@>m_F>> \Sigma(F)\\
@VVV @VV  \Sigma(\Phi_{|F})V\\
\uno_{T'} @>m_{F'}>> \Sigma(F')\ ;\\
\end{CD}
\end{equation}
here the left vertical arrow is the canonical morphism lifting $\tau: T\to T'$
and we denote by $\Phi_{|F}: F\to F'$ the isomorphism of $2$-cubes obtained 
from $\Phi$ ``by restriction".
 
Now suppose that $K$, $K'$ are two $I$-cubes in $\PP$ and let  
$u: {\rm Front}_i(K)\xrightarrow{\sim} {\rm Back}_i(K')$ be an isomorphism of
$I-\{i\}$-cubes for some $i\in I$. Then we can define the ``glueing" $K\ast_iK'$ 
of $K$ and $K'$ along $u$; this
is an  $I$-cube with
${\rm Back}_i(K\ast_iK')={\rm Back}_i(K)$, ${\rm Front}_i(K\ast_iK')={\rm Front}_i(K')$. 
By combining the isomorphisms given by (\ref{sigmadec}) with $\Sigma(u)$ we obtain a natural isomorphism
\begin{equation}\label{sigmaglue}
\Sigma(K\ast_iK')\xrightarrow{\sim} \Sigma(K)+\Sigma(K').
\end{equation} 
If both $K$ and $K'$
are decorated and the isomorphism $u$ is compatible with the $2$-face trivializations 
in the sense above (i.e if it is an isomorphism of decorated cubes), 
then $K\ast_iK'$ becomes naturally a decorated $I$-cube with 
$2$-face trivializations induced by those of $K$ and $K'$ ([Du] \S 1.5). 
Also note that for $i\in I$, $J=I-\{i\}$, if $A$ is a decorated $J$-cube, then the 
$I$-cube $A\buildrel{i}\over{\hbox{\rm ---}}A$ is also naturally decorated.

\subsection{} Let $\PP$, $\QQ$ be two s.c. Picard $\Sfr$-stacks and
suppose that $\delta: \PP\to \QQ$ is an  $\Sfr$-functor. We will also denote by $\delta$
the $\Sfr$-functor
\begin{equation*}
\delta:  n\mi Cube(\PP)\to n\mi Cube(\QQ)\ 
\end{equation*}
which sends the  $I$-cube $K=(K_v)_{v\in C(I)}$
to $(\delta(K_v))_{v\in C(I)}$. 

Compose with 
$\Sigma: n\mi Cube(\QQ)\to\QQ$ to obtain a functor:
\begin{equation*}
\delta^n:=\Sigma\cdot \delta: n\mi Cube(\PP)\to\QQ\ ;\quad 
\delta^n(K)=\sum_{v\in C(I)}(-1)^{s(v)}\delta(K_v)\ .
\end{equation*}
We can observe that if $K$, $K'$ are  $I$-cubes with
an isomorphism 
$u: {\rm Front}_i(K)\xrightarrow{\sim}{\rm Back}_i(K')$, $i\in I$,
then by using (\ref{sigmaglue}) we can obtain
a (natural) isomorphism
\begin{equation}\label{tri2}
\delta^n(K\ast_iK')\xrightarrow{\sim} \delta^n(K)+\delta^n(K').
\end{equation}
Also notice that for $i\in I$, $J=I-\{i\}$, and $A$ a decorated $J$-cube,
(\ref{sigmadec2}) provides us with an isomorphism
\begin{equation}\label{tri3}
\uno \xrightarrow{\sim} \delta^n(A\buildrel{i}\over{\hbox{\rm ---}}A).
\end{equation}

We can now restrict $\delta^n$ to decorated $n$-cubes to obtain a functor
\begin{equation*}
\delta^n: n\mi Cube_d(\PP)\to\QQ\,. 
\end{equation*}

Let us also consider the trivial functor 
\begin{equation*}
e: n\mi Cube_d(\P)\to \QQ\ ;\ \ \ e(K)=\uno,\ e(\Phi)={\rm Id}_{\uno}.
\end{equation*}

 \begin{Definition}\label{cubeDudef}
{\rm {(Ducrot [Du])} An $n$-cubic structure on $\delta$ is an $\Sfr$-functor isomorphism 
\begin{equation*}
\Xi: e\xrightarrow {\sim}\delta^n
\end{equation*}
(between functors on decorated $n$-cubes)
which satisfies the following ``glueing" condition:
If $K$, $K'$ are decorated $I$-cubes, $\#I=n$,
with a (compatible) isomorphism 
$u: {\rm Front}_i(K)\xrightarrow{\sim}{\rm Back}_i(K')$, for some $i\in I$,  
then the diagram
\begin{equation*}
\begin{CD}
\uno+\uno@>\Xi(K)+\Xi(K')>> \delta^n(K)+\delta^n(K')\\
@A\wr AA @AA(\ref{tri2})A\\
\uno@>\Xi(K\ast_iK')>>\delta^n(K\ast_iK')\\
\end{CD}
\end{equation*}
is commutative. 
}
\end{Definition}

\begin{Remark} \label{defrem}{\rm a) In [Du], one actually 
finds a definition of the notion of $n$-cubic structure 
on a functor $\delta:\PP\to\QQ$ where $\QQ$ is a Picard stack which is not necessarily
strictly commutative. This definition involves a functor $\Sigma:n\mi Cube(\QQ)\to \QQ$ ([Du] 1.3.2)
which generalizes the functor $\Sigma$ given in the strictly commutative case above.
(Roughly speaking, to define $\Sigma$, Ducrot adjusts  our ``naive" definition
by using an appropriate combination of the signs $\varepsilon(K_v)$ that give the obstruction
to strict commutativity; when $\QQ$ is s.c. this coincides with the definition that
we have used here). This definition of $n$-cubic structure in [Du] also includes an additional condition:

(``Normalization")\  For $i\in I$, $A$ a decorated $I-\{i\}$-cube and 
$K=(A\buildrel{i}\over{\hbox{\rm ---}}A)$ the corresponding decorated $I$-cube,
the isomorphism $\Xi(K):\uno\xrightarrow{\sim} \delta^n(K)$ is equal to
the isomorphism (\ref{tri3}). 

We can see that when $\QQ$ is s.c this normalization condition is implied
by the glueing condition above. Hence,  in this case, the definition above is equivalent
to the definition of [Du]. Indeed, apply the glueing condition to $K=K'=(A\buildrel{i}\over{\hbox{\rm ---}}A)$
with $u$ the identity.
We have $K\ast_i K=K$ and so we obtain a commutative diagram
\begin{equation*}
\begin{CD}
\uno+\uno@>\Xi(K)+\Xi(K)>> \delta^n(K)+\delta^n(K)\\
@A\wr AA @AA(\ref{tri2})A\\
\uno@>\Xi(K)>>\delta^n(K).\\ 
\end{CD}
\end{equation*}
Now observe that in the s.c. Picard stack $\QQ$, there is a similar 
commutative diagram but with $\Xi(K)$ replaced by 
the isomorphism (\ref{tri3}). This implies that $\Xi(K)$ coincides with $(\ref{tri3})$
since the pair $(\uno, \epsilon)$
of an identity object of $\QQ$ together with $\epsilon: \uno+\uno \xrightarrow {\sim}\uno$
is unique up to {\sl unique} isomorphism. 

b) Suppose that $K$ is a decorated $I$-cube for $I=\{1,\ldots, n\}$ and let $\sigma\in S_n$. The permutation
$\sigma$ induces a morphism of decorated $n$-cubes $\sigma: \sigma^*K\to K$. This 
provides us with an isomorphism 
\begin{equation}\label{tri4}
\delta^n(\sigma): \delta^n(\sigma^*K)\xrightarrow {\sim}\delta^n(K).
\end{equation}
Since $\Xi$ is supposed to be a functor isomorphism we must have 
\begin{equation}
\delta^n(\sigma)\cdot \Xi(\sigma^*K)=\Xi(K).
\end{equation}}
\end{Remark}
\smallskip

\subsection{} It follows directly from the definitions that an additive $\Sfr$-functor $F: \P\to \QQ$
induces an $\Sfr$-functor between decorated $n$-cubes 
\begin{equation}
 F: n\mi Cube_d(\P)\to n\mi Cube_d(\QQ)\ .
\end{equation}

\begin{lemma}\label{addcube}
If $\delta: \PP\to\QQ$ supports an $n$-cubic structure and $F: \PP'\to \PP$
is an additive functor, then the composite $\delta\cdot F: \PP'\to \QQ$
also supports an $n$-cubic structure.
\end{lemma}

\begin{Proof}
Since by definition $(\delta\cdot F)^n= \delta^n\cdot F$,
a trivialization $\Xi$ of $\delta^n$ defines a trivialization of $(\delta\cdot F)^n$.
Now if $K$, $K'$ are decorated $I$-cubes  with a (compatible) isomorphism
$u: {\rm Front}_i(K)\xrightarrow{\sim} {\rm Back}_i(K')$, then $F(u)$ gives a compatible
isomorphism between the $i$-th front and $i$-th back faces of $F(K)$
and $F(K')$ respectively and we have
\begin{equation}\label{addglue}
 F(K)\ast_i  F(K')=  F(K\ast_i K')
\end{equation}
(as decorated $I$-cubes).
We can now see using (\ref{addglue}) that the above 
trivialization satisfies the glueing property
of the definition.  
\endproof
\end{Proof}
 
\subsection{} Recall that we denote by 
${\mathcal {PIC}}(T)$ the s.c. Picard $S_{\rm fppf}$-stack 
given by invertible $\O_{T\times_SS'}$-sheaves on the schemes
$T\times_SS'$, $S'\to S$ a fppf morphism. Let $H\to S$ be a fppf abelian
group scheme. By Yoneda equivalence, there is natural bijection between $S_{\rm fppf}$-functors 
$\delta: H\to {\mathcal {PIC}}(S)$ (where $H$ also denotes the  
s.c. Picard stack represented by the group scheme) and
invertible $\O_H$-sheaves on the scheme $H$ given by $\delta\mapsto \L_\delta:=\delta(H\xrightarrow{{\rm id}} H)$.

\begin{lemma}\label{twocubes}
Suppose that the functor $\delta: H\to {\mathcal {PIC}}(S)$ 
supports an $n$-cubic structure
 in the sense of Definition 
\ref{cubeDudef} above. Then the invertible sheaf $\L_\delta$ over $H$
is equipped with an $n$-cubic structure 
in the sense of Definition \ref{cubicdef}. 
\end{lemma}
\begin{Proof} For simplicity, set $\L=\L_\delta$. Let $I=\{1,\ldots, n\}$.
We can see that, since the Picard stack $H$ is discrete, decorated $I$-cubes
in $H$ correspond bijectively to ordered $n+1$-tuples $(a_i)_i$,
$i=0,\ldots ,n$, of $S'$-valued points $a_i\in H(S')$:
\begin{equation}
(a_i)_i\mapsto K_{a_0}(a_1,\ldots, a_n):={\text{the cube given by\ }} K_v=a_0+\sum_{i, v_i=1}a_i.
\end{equation}
For simplicity, if $a_0=0$, we will denote the above cube by
$K(a_1,\ldots, a_n)$. Now notice that, by the definition of $\delta^n$,
there is a canonical isomorphism
\begin{equation}\label{twocubes1}
\delta^n(K(a_1,\ldots, a_n))\simeq (a_1,\ldots, a_n)^*\Theta_n(\L).
\end{equation}
First suppose that $\delta$ supports an $n$-cubic structure $\Xi$; then 
by (\ref{twocubes1}), the functor isomorphism $\Xi$ evaluated at the $n$-cube $K(a_1,\ldots, a_n)$
with $S'=H^n$ and $a_i={{\rm pr}_i}: H^n\to H$, gives a trivialization $\xi$ of 
the invertible sheaf $\Theta_n(\L)$ over $H^n$. Now take all $a_i=0$ (as $S$-points) and set for simplicity 
$K(\underline{0})=K(0,\ldots, 0)$. The $n$-cubic structure $\Xi$ provides us
with an isomorphism
\begin{equation*}
\Xi(K(\underline{0})): \O_S\xrightarrow{\sim} \delta^n(K(\underline{0}))
\end{equation*}
which by Remark \ref{defrem} (a) agrees with the natural $\O_S\xrightarrow{\sim} \delta^n(K(\underline{0}))$
given by the contraction isomorphisms. In view of 
(\ref{twocubes1}) this implies that $\xi$ satisfies
Property (c0) of Definition \ref{cubicdef}. To examine Property (c1)
take again $S'=H^n$ and $a_i={{\rm pr}_i}: H^n\to H$.
We have $\sigma^*K(a_1,\ldots, a_n) =K(a_{\sigma(1)},\ldots, a_{\sigma(n)})$
and under (\ref{twocubes1}) the isomorphism (\ref{tri4}) corresponds to
${\mathfrak P}_\sigma$ (see (\ref{permu})). Hence,
we can see that Remark \ref{defrem} (b) implies that 
$\xi$ satisfies (c1) of Definition \ref{cubicdef}.
It remains to discuss Property (c2). For this, we take $S'=H^{n+1}$ and 
$a_i={{\rm pr}_i}: H^{n+1}\to H$, $0\leq i\leq n$.
Notice that we have
\begin{equation}\label{twocubes3}
K(a_0,a_1,a_3,\ldots, a_{n})\ast_2 K_{a_1}(a_0,a_2,\ldots, a_n)=K(a_0,a_1+a_2,a_3,\ldots, a_n),
\end{equation}
\begin{equation}\label{twocubes4}
K(a_1,a_2,\ldots, a_{n})\ast_1 K_{a_1}(a_0,a_2,\ldots, a_n)=K(a_0+a_1,a_2,\ldots, a_n).
\end{equation}
For example, when $n=2$ these are explained by the following diagram of three ``glued" squares
\begin{equation}
\label{associativite2}
\xymatrix{
a_2 \ar@{-}[d]\ar@{-}[r] 
&a_1+a_2 \ar@{-}[d] \ar@{-}[r] 
&a_0+a_1+a_2 \ar@{-}[d]   \\
0 \ar@{-}[r]
&a_1 \ar@{-}[d] \ar@{-}[r] \ar@{-}[d]
&a_0+a_1 \ar@{-}[d]   \\
&0 \ar@{-}[r] &a_0
}
\end{equation}
(the general case is just harder to draw: Here the vertical
glueing on the right   corresponds to (\ref{twocubes3}) while the horizontal
to (\ref{twocubes4})). The relations (\ref{twocubes3}), (\ref{twocubes4})
together with (\ref{tri2}) give isomorphisms
\begin{multline*}
\delta((K(a_0,a_1,a_3,\ldots, a_n))\otimes\delta((K_{a_1}(a_0,a_2,\ldots, a_n))\simeq 
\delta((K(a_0,a_1+a_2,a_3,\ldots, a_n))\, ,\\ 
\delta((K(a_1,a_2,\ldots, a_n))\otimes\delta((K_{a_1}(a_0,a_2,\ldots, a_n))\simeq 
\delta((K(a_0+a_1,a_2,\ldots, a_n))\, . \ \ \ \ \
\end{multline*}
Combining these gives:
\begin{multline*}
\delta^n(K(a_0+a_1,a_2,\ldots, a_n))\otimes 
\delta^n(K(a_0,a_1,a_3,\ldots, a_{n}))\simeq \\
\simeq \delta^n( K(a_0,a_1+a_2,a_3,\ldots, a_n))\otimes
\delta^n( K(a_1,a_2,\ldots, a_n))\, . 
\end{multline*}
We can see, by  using (\ref{twocubes1}), that this isomorphism corresponds to $\mathfrak Q$
of (\ref{Q}). It now follows that the glueing condition on $\Xi$ implies that 
the trivialization $\xi$ respects the isomorphism $\mathfrak Q$.
In other words, Property (c2) of Definition \ref{cubicdef} is true for $\xi$.
\endproof
\end{Proof}

\begin{Remarknumb} {\rm 
In fact, it is true that there is a bijective correspondence between $n$-cubic structures
on the functor $\delta: H\to {\mathcal {PIC}}(S)$ in the sense of Definition 
\ref{cubeDudef} above and $n$-cubic structures on $\L_\delta$
in the sense of Definition \ref{cubicdef}: Assume that $\L=\L_\delta$ 
supports an $n$-cubic structure $\xi$ in the sense of 
\ref{cubicdef}. In view of 
(\ref{twocubes1}) the isomorphism $\xi$ gives a functorial
trivialization of $\delta^n(K(b_1,\ldots, b_n))$
for all points $b_1,\ldots ,b_n$ of $H$. The relation
(\ref{twocubes4}) together with (\ref{tri2}) give a canonical
isomorphism
\begin{equation}\label{twocubes5}
\delta^n(K_{b_0}(b_1,\ldots, b_n))\simeq 
\delta^n(K(b_0+b_1,\ldots, b_n))\otimes \delta^n(K(b_0, b_2,\ldots, b_n))^{-1}.
\end{equation}
We can now use (\ref{twocubes5}) and the aforementioned 
trivialization of $\delta^n(K(b_1,\ldots, b_n))$ to define a functorial trivialization of 
$\delta^n(K_{b_0}(b_1,\ldots, b_n))$. One can verify that Properties
(c0), (c1) and (c2) of Definition \ref{cubicdef} imply that this
trivialization of $\delta^n(K_{b_0}(b_1,\ldots, b_n))$
defines an $n$-cubic structure on $\delta$. Since we are not going to use this
statement we choose to omit the  details. }
\end{Remarknumb}

\subsection{} \label{det}
Let $h: Y\to S$ be a projective and flat morphism of
relative dimension $d$ over the
locally Noetherian scheme $S$. If $\Hh$ is a locally free coherent sheaf
of $\O_Y$-modules on $Y$ then the total derived image ${\bf R}h_*(\Hh)$ 
in the derived category of complexes of sheaves of $\O_S$-modules which are bounded below
is ``perfect" (i.e it is locally on $S$ quasi-isomorphic to a 
bounded complex of finitely generated free $\O_S$-modules, see [SGA6] III). Hence, by 
[KM], we can associate to ${\bf R}h_*(\Hh)$ a graded invertible 
sheaf 
\begin{equation*}
{\det}_* {\bf R}h_* =(\det {\bf R}h_*(\Hh), {\rm rk}({\bf R}h_*(\Hh))
\end{equation*}
on $S$ (the ``determinant of cohomology").
By restricting to $\Hh$ which are invertible, we obtain a functor, which we will denote 
again
by
$
\det_* {\bf R}h_* 
$,
from the Picard category of invertible sheaves
on $Y$ to the Picard category of graded invertible sheaves on $S$.
In fact, by loc. cit. the formation of the determinant of cohomology 
commutes with arbitrary base changes $S'\to S$; hence we can define 
a corresponding $S_{\rm fppf}$-functor 
\begin{equation*}
{\det}_* {\bf R}h_*: {\mathcal {PIC}}(Y)\to {\mathcal {PIC}}_*(S).
\end{equation*}
By the main result of [Du] (Theorem 4.2) this functor 
supports a canonical $d+2$-cubic structure (in the sense of [Du] Definition 1.6.1).
Let us remark here that since ${\mathcal {PIC}}_*(S)$ is not strictly commutative 
we do have to refer to [Du] for the definition 
of cubic structure. (See Remark \ref{defrem} (a).) Now observe
that the substack ${\mathcal {PIC}}_{\rm ev}(S)$ of ${\mathcal {PIC}}_*(S)$ given by graded line bundles $(\M, e)$
with $e$ always {\sl even} is a s.c. Picard stack. The forgetful functor 
$f: {\mathcal {PIC}}_{\rm ev}(S)\to {\mathcal {PIC}}(S)$ is an additive functor.
(Note that this is not true for the forgetful functor ${\mathcal {PIC}}_{*}(S)\to {\mathcal {PIC}}(S)$.)
This fact together with [Du] Theorem 4.2 and Remark \ref{defrem} (a) implies the following:
\medskip

i) The functor ${\mathcal {PIC}}(Y)\to  {\mathcal {PIC}}(S)$
given by $\Hh\mapsto \det {\bf R}h_*(\Hh)^{\otimes 2}$ 
supports a canonical $d+2$-cubic structure in the sense 
of Definition \ref{cubeDudef}.
\medskip

ii) Let ${\mathcal {PIC}}^{\rm ev}(Y)$ be the substack
of ${\mathcal {PIC}}(Y)$ given by invertible  sheaves
for which the function ${\rm rk}({\bf R}h_*(\Hh))$
is always even. Then the functor ${\mathcal {PIC}}^{\rm ev}(Y)\to  {\mathcal {PIC}}(S)$
given by $\Hh\mapsto \det {\bf R}h_*(\Hh)$ 
supports a canonical $d+2$-cubic structure in the sense 
of Definition \ref{cubeDudef}.

\begin{Remark} \label{easy}
{\rm 
In fact, one can give a somewhat more direct proof of (i) and (ii) by following the general strategy of the proof of
the main theorem of [Du]. The argument is considerably less involved 
since we do not have to deal with the thorny problem
of signs that complicates the proof of [Du] Theorem 4.2.}
\end{Remark}

\subsection{} Let $h: Y\to S$ be as in the 
previous paragraph and assume in addition that $S$ is the spectrum of a Dedekind ring $R$
with field of fractions $K$. 
Suppose that $\pi: X\to Y$ is a $G$-torsor 
with $G$ a finite abelian group. The construction of 
\S\ref{2d} gives an additive $S_{\rm fppf}$-functor
\begin{equation}
F: G^D_S\to  {\mathcal {PIC}}(Y)\ ;\quad F(\chi)=\O_{X,\chi} .
\end{equation}
Let  $h_{G^D}: Y\times_S G^D_S\to G^D_S$ be the base change of $h: Y\to S$;
we can view $\pi_*(\O_X)$ as an invertible sheaf on $Y\times_S G^D$  
which is isomorphic to the value $\O_{X,\chi_0}$ of the functor $F$ on the 
``universal" point $\chi_0={\rm id}: G^D_S\to G^D_S$
(see (\ref{27})). Now notice that the invertible $\O_{Y\times_SG^D}$-sheaf $\pi_*(\O_X)$
is $\#G$-torsion (this follows  
from (\ref{27}) and (\ref{torsorIso})). Therefore,  
we have an equality of Euler characteristics (ranks)
\begin{equation*}
{\rm rk}({\bf R}{h_{G^D}}_*(\pi_*(\O_X)))={\rm rk}({\bf R}{h_{G^D}}_*(\O_{Y\times_SG^D_S})).
\end{equation*}
(For example, this equality can be deduced as follows: In the Grothendieck group
${\rm K}_0(Y\times_S G^D)$ of locally free coherent modules on $Y\times_SG^D$
the rank zero element $z=[\pi_*(\O_X)]-[\O_{Y\times_S G^D}]$ is nilpotent.  
Now since $\pi_*(\O_X)$ is a $\#G$-torsion invertible sheaf, it follows that a power of $\#G$ annihilates 
$z$. Hence, $z$ has  Euler characteristic equal to zero and this implies the equality.) 
Using the projection formula and the flatness of $Y\to S$, we see that
the   locally constant function $G^D_S\to \Z$ given by ${\rm rk}({\bf R}{h_{G^D}}_*(\O_{Y\times_SG^D_S}))$
always takes the value $\chi(Y,\O_Y)=\chi(Y_K, \O_{Y_K}):=\sum_{i}(-1)^i{\rm rk}_K{\rm H}^i(Y_K,\O_{Y_K})$.
We can now conclude using base change that for any $S'$-valued point $\chi$ of $G^D_S$ 
the function ${\rm rk}({\bf R}{h_{S'}}_*(\O_{X,\chi}))$ is constant with value $\chi(Y_K, \O_{Y_K})$.
This value is even if and only if the 
arithmetic genus $g(Y_{K})=(-1)^d(\chi(Y_{K}, \O_{Y_K})-1)$  of the generic fiber is odd. 
 
Now let $\kappa={\rm gcd}(2, g(Y_{K}))$
and consider the composition 
\begin{equation}
\delta:=\det {\bf R}h_*^{\otimes \kappa}\cdot F: G^D_S\to {\mathcal {PIC}}(S)
\end{equation}
given by $\chi\mapsto \det {\bf R}h_*(\O_{X,\chi})^{\otimes \kappa}$. By 
Lemma \ref{addcube} and (i), (ii) of \S \ref{det} above (applied when $g(Y_K)$ is even, odd respectively) 
we conclude that the functor $\delta$ has a canonical $d+2$-cubic 
structure in the sense of Definition \ref{cubeDudef}.
Using Lemma \ref{twocubes} and $\pi_*(\O_X)=\O_{X,\chi_0}$, we see that
the $d+2$-cubic structure on $\delta $ 
equips the invertible sheaf $\L_{\delta }=\det {\bf R}{h_{G^D}}_*(\pi_*(\O_X))^{\otimes \kappa}$  over $G^D_S$
with a corresponding $d+2$-cubic structure in the sense of Definition \ref{cubicdef}. We summarize:

\begin{thm}\label{detC}
Let $h: Y\to S$ be a projective and flat morphism of  relative dimension $d$ over the spectrum of a Dedekind ring
with field of fractions $K$; set $\kappa={\rm gcd}(2, g(Y_{K}))$ with $g(Y_{K})=(-1)^d(\chi(Y_{K}, \O_{Y_K})-1)$ the 
arithmetic genus of the generic fiber of $Y\to S$. Suppose that 
$\pi: X\to Y$ is a $G$-torsor for the
finite abelian group $G$. Then the invertible sheaf 
$\det {\bf R}{h_{G^D}}_*(\pi_*(\O_X))^{\otimes \kappa}$  over $G^D_S$
supports a $d+2$-cubic structure.\endproof
\end{thm}

\begin{Remark} \label{RRexpl}
{\rm  The Taylor expansion formula (\ref{flat2}) for $\L=\det {\bf R}{h_{G^D}}_*(\pi_*(\O_X))^{\otimes \kappa}$
(with its $d+2$-cubic structure provided by Theorem \ref{detC}) closely resembles a Riemann-Roch formula
``without denominators". To explain this, assume for simplicity that the morphism 
$h:Y\to S$ is smooth and that $\kappa=1$. Let $\chi: G\to {R'}^*$
be a character of $G$ with values in the Dedekind ring $R'$ 
and consider the invertible sheaf $\O_{X,\chi}$ on $Y'=Y\times_RR'$. The Grothendieck-Riemann-Roch
theorem for the morphism $h':Y'\to S'=\Spec(R')$ gives
\begin{equation}\label{RRfor}
[\det {\bf R}{h'_*}(\O_{X,\chi})]=[\det {\bf R}{h'_*}(\O_{Y'})]+
\sum_{i=0}^{d} h'_*\biggl[\frac{\ \ c_1(\O_{X,\chi})^{d+1-i}}{(d+1-i)!\ \ }
\cap {\rm Td}_{i}((\Omega^1_{Y/R})^*)\biggr]
\end{equation}
in ${\rm Pic}(S')_\Q={\rm CH}^1(S')_\Q$ (see for example [Fa] Theorem 1.7 or [Fu] \S 15).
Of course, this result is of little use since in our case of interest 
${\rm Pic}(S')_\Q={\rm Pic}(R')\otimes_\Z\Q=(0)$. Our point is that the Taylor expansion (\ref{flat2})
provides a ``without denominators" version of this equality. 
In fact, we can conjecture a
precise relationship between the terms of (\ref{flat2}) and those of 
(\ref{RRfor}): Let $\chi: S'\to G^D_S$
be the morphism given by the character $\chi$. First of all, we have
$
\chi^*(\L)\simeq\det {\bf R}{h'_*}(\O_{Y,\chi})$,
$\chi^*(0^*\L)\simeq \det {\bf R}{h'_*}(\O_{Y'})$.
Let us consider the rest of the terms. The standard formula for the Todd class gives a (universal) expression 
$
{\rm Td}_i((\Omega^1_{Y'/R'})^*)=N_i\cdot T_i
$
with $N_i\in \Q$, $T_i$ an integral linear combination 
of Chern classes of $\Omega^1_{Y'/R'}$ in ${\rm CH}^i(Y')$
($N_0=1$, $N_1=1/2$, $N_2=1/12$, $N_3=1/24$, $N_4=1/720$; in general $N_i$ is
given using Bernoulli numbers).
We conjecture that for all $i=0,\ldots, d$
\begin{equation}
[\chi^*(\delta(\L^{(i)},\xi^{(i)})]=
(-1)^{i}\frac{(d+1)!!N_i}{(d-i)!!(d+1-i)!}h'_*\biggl[c_1(\O_{X,\chi})^{d+1-i}\cap T_i\biggr]
\end{equation}
in ${\rm Pic}(S')$. (It follows from the theorem of von Staudt-Clausen  -[Wa] Theorem 5.10-
that this coefficient is integral.) In the case of arithmetic surfaces ($d=1$),
this conjecture follows from the Deligne-Riemann-Roch theorem ([De] (7.5.1)). In fact, in this case 
we can see that $
\delta(\L^{(0)},\xi^{(0)})\simeq <\pi_*(\O_X),\pi_*(\O_X)>$, 
$
\delta(\L^{(1)},\xi^{(1)})\simeq <\pi_*(\O_X),\omega_{Y\times_S G^D_S/G^D_S}>
$
(notations as in loc. cit.) and that the Taylor expansion (\ref{flat2})
amounts to the Deligne-Riemann-Roch formula. We will leave the details 
for another occasion.
}\end{Remark}
\bigskip

\section{The multiextension associated to a $\Z/p$-torsor.}\label{multidet}

In this section, we assume that $G=\Z/p$, $p$ a prime, and that $\pi: X\to Y$ is a $G$-torsor 
with $Y\to S=\Spec(\Z)$ projective and flat of relative dimension $d$.  
In addition, we will
assume that the generic fiber $Y_\Q$ is normal. By Theorem \ref{detC} the line bundle
$\L_{\delta}=\det {\bf R}{h_{G^D}}_*(\pi_*(\O_X))^{\otimes \kappa}$ on $G^D$ is 
equipped with a canonical $d+2$-cubic structure $\xi$. 
We will now give a rather explicit description of the 
isomorphism class of the corresponding (via the construction of \S \ref{cubicext}) 
$d+1$-extension $E(\L_{\delta },\xi)$ on 
$\Theta_{d+1}(\L_{\delta })$. 
Using Remark \ref{diagonal} we see that it is enough to describe the
image of the class of $E(\L_{\delta },\xi)$ under the
homomorphism
\begin{multline}\label{compo1new}
\ \ \ \ {d+1}{\hbox{\small -}}{\rm Ext}^1(\mu_{p}, \Gm)\xrightarrow{\psi_{d+1}'}  
{\rm H}^1(\mu_{p},\Z/p)\hookrightarrow  {\rm H}^1(\Spec(\Q(\zeta_{p})),\Z/p),\ \ \ \ 
\end{multline}
where  the second arrow is
given by restriction along the generic point $\Spec(\Q(\zeta_{p}))\to\mu_p$.
 Now set $U=X _{\Q(\zeta_p)}$, 
$V=Y _{\Q(\zeta_p)}$. Then $U\to V$ is a $G$-torsor of 
normal projective schemes of dimension $d$ over $\Q(\zeta_p)$. 
Consider  the character $\chi_0: \Z/p\Z\to \Q(\zeta_p)^*$ given by $\chi_0(1)=\zeta_p$. 
By \S \ref{2d}, the $\chi_0$-isotypic part 
$\O_{U,\chi_0}$ of the $\O_V$-sheaf $\pi_*(\O_U)$ is an invertible sheaf on $V$
which is $p$-torsion. Let $A$ be a generic section of this invertible sheaf; $A$
is then an element of the function field $k(U)$ and
$a\in \Z/p$ acts on it by multiplication by $\chi_0(a)=\zeta_p^a$.
Hence,
 $F=A^p$ belongs to the function field $k(V)$. 
Using the section $A$ we can associate
to the invertible sheaf $\O_{U,\chi_0}$ a Cartier divisor $D$ on $V$.
We can see that
 $p\cdot D$ is equal to the principal  Cartier divisor $(F)$. 

Consider the
group $Z_0(V)$ of $0$-cycles on $V$ and the 
Chow group ${A}_0(V)$ of $0$-cycles on $V$ modulo 
rational equivalence (see [Fu], 2.4). Suppose that $P$ is a closed point of $V$
with residue field $k(P)$ (a finite extension of $\Q(\zeta_p)$). 
Specializing the cover $U\to V$ at $P_i$ gives a
$G$-torsor over $\Spec(k(P))$ which, by Kummer theory (Remark \ref{rem23}), 
corresponds to $f\in k(P)^*/(k(P)^*)^p$. Set $\lambda(P)=
{\rm Norm}_{k(P)/\Q(\zeta_p)}(f)  \in \Q(\zeta_p)^*/(\Q(\zeta_p)^*)^{p}.
$
By extending linearly we obtain a group homomorphism 
$$
\lambda: Z_0(V)\to \Q(\zeta_p)^*/(\Q(\zeta_p)^*)^{p} .
$$

\begin{lemma}\label{chow}
The homomorphism $a$ factors through rational equivalence and
provides us with a group homomorphism
$$
\lambda: {A}_0(V)\to \Q(\zeta_p)^*/(\Q(\zeta_p)^*)^{p}\ .
$$
\end{lemma}

\begin{Proof}
Consider two $0$-cycles $Z$, $Z'$ which are rationally equivalent.
By definition, there is then a projective reduced curve $C\subset V$
passing through all points in the supports of $Z$ and $Z'$
and a function $\phi$ on the normalization $q: \tilde C\to C$
such that $q_*((\phi))=Z-Z'$ where the push-forward $q_*$ of the divisor
$(\phi)$ by $q$ is defined as in [Fu] 1.4. A standard moving argument,
shows that we can choose a generic section $A$ of $\O_{U,\chi_0}$
as above, such that the divisor of $F=A^p$ is disjoint
from all the points in the supports of $Z$ and $Z'$
and from the generic points of $C$. Specializing $F$ on $C$
gives an element $f$ in the function algebra (the product of the
function fields of the irreducible components) $k(C)=k(\tilde C)$.
Now if $Q$ is a point of $ C$ (resp. of $\tilde C$) which is not in
the support of $q_*((f))$ (resp. of $(f)$) then we denote by $f(Q)$ the Norm from 
$k(Q)$ to $\Q(\zeta_p)$ of the ``evaluation" of $f$ at $Q$. Extending by 
linearity, we can make sense of $f(Z)$ and $f(Z')$.
By our choice of $f$, we have
$$
\lambda(Z)=f(Z)\ {\rm mod\ } (\Q(\zeta_p)^*)^p,\qquad 
\lambda(Z')=f(Z')\ {\rm mod\ } (\Q(\zeta_p)^*)^p,
$$
(see Remark \ref{rem23} (a)) and therefore 
$\lambda(Z)\cdot \lambda(Z')^{-1}=f(q_*((\phi)))\ {\rm mod\ } (\Q(\zeta_p)^*)^p$.
By the definitions of $q_*$ and $f(-)$ we now obtain that
$f(q_*((\phi)))=f((\phi))$ where in the 
second expression $f$ (resp. $(\phi)$) is regarded as
a function (resp. a divisor) on the smooth projective curve 
$\tilde C$. By Weil reciprocity on $\tilde C$,  $f((\phi))=\phi((f))$.
Now observe that the divisor $(f)$  on $\tilde C$ is the $p$-th multiple of the 
pull back of the Cartier divisor $D$ via $\tilde C\to V$. Therefore, $f(q_*((\phi)))=\phi((f))$ is in 
$(\Q(\zeta_p)^*)^p$ and this completes the proof.\endproof
\end{Proof}
\smallskip
\begin{Remark} \label{CoTh}
{\rm The homomorphism $\lambda$ is closely connected to the 
cohomological Abel-Jacobi map described in [C-TS]. For simplicity, set $K=\Q(\zeta_p)$.
Following Bloch,  Colliot-Th\'el\`ene and Sansuc construct a ``characteristic homomorphism"
\begin{equation}
\Phi: A_0(V)\to 
{\rm Ext}^1_{{\rm Gal}(\bar K/K)}({\rm Pic}(V_{\bar K}), {\bar K}^*)
\end{equation}
(extensions of discrete ${\rm Gal}(\bar K/K)$-modules). One can see from their construction
that pulling back along $\Z/p\to {\rm Pic}(V_{\bar K})$
given by $1\mapsto D$, and composing with the isomorphism  ${\rm Ext}^1_{{\rm Gal}(\bar K/K)}(\Z/p, {\bar K}^*)\simeq K^*/(K^*)^p$ given by Hilbert's Theorem 90 gives the homomorphism $\lambda$ above.}
\end{Remark}

Since $D$ is a Cartier divisor, the $d$-th self intersection 
$D^d:=D\cap\, \cdots\, \cap D$ makes sense
as an element of the Chow group 
${A}_0(V)$ (see [Fu], 2.4) and we can set
\begin{equation*}
\lambda(X/Y)= \lambda (D^d)\ .
\end{equation*}
Recall $\kappa={\rm gcd}(2,g(Y_\Q))$.

\begin{prop}\label{explicit}
Set $L=\Q(\zeta_p)[y]/(y^p-\lambda(X/Y)^{-\kappa})$; then $T(X/Y)=\Spec(L)$ is a $\Z/p$-torsor over 
$\Spec(\Q(\zeta_p))$ whose class in ${\rm H}^1(\Spec(\Q(\zeta_{p})),\Z/p)$
is the image of the isomorphism class of the $d+1$-extension $E(\L_{\delta},\xi)$
under the injective map (\ref{compo1new}).
\end{prop}

\begin{Proof}  For simplicity, let us denote by $E(X/Y)$ the multiextension $E(\L_{\delta },\xi)$
associated to the torsor $X\to Y$.
Consider 
the base change  of the multiextension $E(X/Y)$
by the morphism $\Spec(\Q)\to S=\Spec(\Z)$. This is isomorphic to the multiextension $E(X_\Q/Y_\Q)$
associated to the torsor $X_\Q\to Y_\Q$.
Let $\tau$ be 
the image of the isomorphism class of $E(X_\Q/Y_\Q)$ under the map
\begin{multline}\label{compo1new1}
\ \ \ \ {d+1}{\hbox{\small -}}{\rm Ext}^1({\mu_{p}}_{/_\Q}, \Gm_{/_\Q})\xrightarrow{\psi_{d+1}'}  
{\rm H}^1 ({\mu_{p}}_{/_\Q},\Z/p)\to  {\rm H}^1(\Spec(\Q(\zeta_{p})),\Z/p).\ \ \ \ 
\end{multline}
(Here ${\mu_{p}}_{/_\Q}$, $ \Gm_{/_\Q}$ is short hand notation for the group schemes
$\mu_p\times_S\Spec(\Q)$, $\Gm\times_S\Spec(\Q)$ over $\Spec(\Q)$. The
first arrow is defined by the construction of (\ref{compo1}) performed 
over the base $\Spec(\Q)$.)
We can see that it is enough to show that  the class of $T(X/Y)$
is equal to the image $\tau$. By unraveling the definition of $\psi'_{d+1}$
(Remark \ref{diagonal})
we can now see that the torsor class $\tau$ corresponds 
(under the ``Waterhouse isomorphism" (\ref{nExt4})) to the 
$\Gm\times_S{\Spec(\Q(\zeta_p))}$-extension over $\mu_p\times_{S}{\Spec(\Q(\zeta_p))}$
given by the additive functor on characters of $\Z/p$
\begin{equation}\label{93}
\chi\mapsto E(X_{\Q}/Y_{\Q})_{(\chi_0, \ldots, \chi_0, \chi)}
\end{equation}
(see \S \ref{extension}). According to the definition of $E(X_\Q/Y_\Q)$,
\begin{equation*}
E(X_{\Q}/Y_{\Q})_{(\chi_0, \ldots, \chi_0, \chi)}
=\bigotimes ^{d+1}_{k=0}\Biggl(\bigotimes_{i_1<\cdots <i_k}\bigl(
{\rm det}{\bf R}h'_*(\O_{U,\chi[i_1]}\otimes\cdots \otimes \O_{U, \chi[i_k]})\bigr)^{\otimes\kappa}\Biggr)^{(-1)^{d+1-k}}
\end{equation*}
with $\chi[i]=\chi_0$ if $i\neq d+1$, $\chi[d+1]=\chi$, and $h':V\to \Spec(\Q(\zeta_p))$ the structure morphism. 
In fact, we can view the right hand side of the above equality
as a value of the ``intersection sheaf" ([De] or [Du] \S 5): If $\L_i$, $i=1,\ldots, d+1$, are 
invertible sheaves on $V$, then by definition
\begin{equation}
I(\L_1,\ldots,\L_{d+1}):=\bigotimes ^{d+1}_{k=0}\Biggl(\bigotimes_{i_1<\cdots <i_k}
{\rm det}{\bf R}h'_*(\L_{i_1}\otimes\cdots \otimes \L_{i_k})\Biggr)^{(-1)^{d+1-k}},
\end{equation}
and so
\begin{equation}\label{94}
E(X_{\Q}/Y_{\Q})_{(\chi_0, \ldots, \chi_0, \chi)}=I(\O_{U,\chi_0},\ldots, \O_{U,\chi_0},\O_{U,\chi})^{\otimes \kappa}.
\end{equation}

Now let us return to the notations before the statement of
the Proposition. A standard argument shows that we can write $\O_{U,\chi_0}=\M_1\otimes \M_2^{-1}$
with $\M_i$ very ample invertible sheaves on $V$
that have sections $s^j_1$ and $s^j_2$, $j=1,\ldots, d$, respectively, with the following property:
All sequences $ (\underline {s})$ of the form $(s^1_{i_1}, s^2_{i_2},\ldots, s^d_{i_d})$,
for $i_1,\ldots, i_d\in \{1,2\}$ are regular.
(By definition, this means that
 the corresponding sequences of elements of the local rings $\O_{V,x}$, $x\in V$, which are obtained
from the $s^j_i$ using local trivializations of $\M_i$ are regular.)
Let us denote by $w(\underline {s})$ the number of indices $j=1,\ldots, d$
in the sequence $ (\underline {s})$ where the subscript is equal to $2$.
For each such sequence $ (\underline {s})$, the scheme theoretic intersection $Z_{(\underline {s})}={\rm div}(s^1_{i_1})
\cap {\rm div}(s^2_{i_2})\cap \cdots \cap {\rm div}(s^d_{i_d})$ is a
union of reduced closed points of $V$ and we have 
\begin{equation}\label{choweq}
\sum_{(\underline {s})}(-1)^{w(\underline {s})}[Z_{(\underline {s})}]=D^d
\end{equation}
in the Chow group ${A}_0(V)$. Using the additivity and restriction properties 
of the intersection sheaves ([Du] \S 5,  [De]) and (\ref{94})
we can see that the functor (\ref{93}) is isomorphic (as an additive functor) to 
the   functor 
\begin{equation}\label{96}
\chi\mapsto \bigotimes_{(\underline {s})}
{\rm Norm}_{Z_{(\underline {s})}/\Spec(\Q(\zeta_p))}
({\O_{U,\chi}}|_{Z_{(\underline {s})}})^{\otimes(-1)^{w(\underline {s})}\cdot\kappa}\ .
\end{equation}
We have a natural isomorphism ${\O_{U,\chi_0}}|_{Z_{(\underline {s})}}\simeq \O_{\pi^{-1}(Z_{(\underline {s})}),\chi_0}$
where $\pi^{-1}(Z_{(\underline {s})})\to Z_{(\underline {s})}$ is the $\Z/p$-torsor obtained
by base changing $\pi_\Q: X_\Q\to Y_\Q$ along $Z_{(\underline {s})}\to Y_\Q$. 
We can now see that, under the description of Remark \ref{rem23} (a), the 
$\Z/p$-torsor $\tau$ is given by the invertible sheaf
$$
\bigotimes_{(\underline {s})}{\rm Norm}_{Z_{(\underline {s})}/\Spec(\Q(\zeta_p))}
(\O_{\pi^{-1}(Z_{(\underline {s})}),\chi_0^{-1}})^{\otimes(-1)^{w(\underline {s})}\cdot\kappa}\
$$
together with the natural trivialization of its $p$-th power induced by
the isomorphism
$\O_{U,\chi^{-1}_0}^{\otimes p} \simeq \O_V$ and the additivity of the norm.
(Recall that the $\Z/p$-torsor $\tau$ corresponds --under (\ref{nExt4})--
  to the extensions given by the additive functors
(\ref{93}) and (\ref{96}). The inverse $\chi_0^{-1}$ in the expression 
above occurs, instead of 
$\chi_0$, because of Remark \ref{LchiRemark} (a)).
It now follows from Remark \ref{rem23} (b) that a Kummer element 
for this $\Z/p$-torsor over $\Spec(\Q(\zeta_p))$ is given by taking the inverse of the 
product over all
$(\underline {s})$ 
of the Norms 
from $Z_{(\underline {s})}$ to $\Spec(\Q(\zeta_p))$ of Kummer elements 
for the $\Z/p$-torsors $\pi^{-1}(Z_{(\underline {s})})\to Z_{(\underline {s})}$.
Hence, the result now follows from (\ref{choweq}) and Lemma \ref{chow}.\endproof
\end{Proof}

\begin{Remark} \label{explicitRe}

{\rm Assume that $p>d+1$.

a)  Proposition  \ref{explicit} gives an explicit description of the
$\Z/p$-extension of $\Q(\zeta_p)$ that corresponds via class field theory 
to the element $t_{d+1}(\L,\xi)\in \Hom(C(p), \Z/p)$
(defined in the  section \S \ref{8d}) when $(\L,\xi)$ is given by 
$\det {\bf R}{h_{G^D}}_*(\pi_*(\O_X))^{\otimes \kappa}$ with its canonical cubic structure. 
Indeed, we can see that 
the class of the $\Z/p$-extension
that corresponds to $t_{d+1}(\L,\xi)$ is  equal to 
$\displaystyle{\frac{1}{(d+1)!}(T(X/Y))}$ in ${\rm H}^1(\Spec(\Q(\zeta_p),\Z/p)$. 
\smallskip

b)   We conjecture that the elements $\{t_{d+1}(\L,\xi)\}$, $X\to Y$ 
ranging over all $\Z/p\Z$-torsors
with $Y\to \Spec(\Z)$ projective and flat of relative dimension $d$,
generate the target group $\Hom((C/p)^{(-d)}, \Z/p\Z)$. 
To this moment we have no direct evidence 
to support  this  conjecture. Any calculation proves to be quite difficult:
Indeed,   the lowest dimension for which
the group $\Hom((C/p)^{(-d)}, \Z/p\Z)$ is non-trivial for $p<12\cdot 10^6$
is $d=11$; then we have a non-trivial group for $p=691$. On the other hand, it is easier
to examine analogous statements over bases more general than $\Spec(\Z)$; some (indirect)
evidence can be collected this way. We intend to return to this topic
in a subsequent paper.

c) Suppose $d\geq 2$, the generic fibers $X_\Q$, $Y_\Q$ are smooth and the $p$-torsion 
divisor $D$ on $V=Y_{\Q(\zeta_p)}$ is algebraically equivalent to zero, i.e corresponds to
a $p$-torsion point in the Picard variety ${\rm Pic}^0(V)$. Then $ t_{d+1}(\L,\xi) =0$.
(Therefore examples of non-trivial elements $ t_{d+1}(\L,\xi) $ 
can only come from $p$-torsion in the Neron-Severi group of $V$.) Let us quickly sketch 
a proof of this fact. (Here we will give a straightforward 
argument. Later, in the course of the proof of Theorem \ref{main3},
we will show a stronger result; see Remark \ref{albtri}). Notice that we have a commutative diagram 
\begin{equation}
\begin{CD}
(A_0(V))_{{\rm degree}=0} @>{\Phi}>> {\rm Ext}^1_{{\rm Gal}(\bar K/K)}({\rm Pic}(V_{\bar K}), {\bar K}^*)\\
@V  VV @VVV\\
{\rm Alb}(V)(K) @>{\Phi}>> {\rm Ext}^1_{{\rm Gal}(\bar K/K)}({\rm Pic}^0(V_{\bar K}), {\bar K}^*).\\
\end{CD}
\end{equation}
Here $K=\Q(\zeta_p)$, ${\rm Alb}(V)$ is the Albanese variety 
of $V$ and $\Phi$ is the homomorphism of [C-TS] as in Remark \ref{CoTh}.
Using that remark and the above, we can see that the homomorphism 
$\lambda: (A_0(V))_{{\rm degree}=0}\to K^*/(K^*)^p$ factors through   
$(A_0(V))_{{\rm degree}=0}\to {\rm Alb}(V)(K)\to K^*/(K^*)^p$. Using Proposition \ref{explicit},
we now see that to show  $ t_{d+1}(\L,\xi) =0$ is enough to show that the self-intersection
$D^d\in (A_0(V))_{{\rm degree}=0}$ has trivial image in ${\rm Alb}(V)(K)$. It is enough to show
that the corresponding $p$-torsion point of ${\rm Alb}(V)(K)$ is trivial in 
${\rm Alb}(V)(k({\mathfrak P}))\subset {\rm Alb}(V)(\overline{k({\mathfrak P})})$
(that is after  reduction modulo some prime $\mathfrak P$ of $K$ that does not divide $p$ and where $V$ has good reduction;
her $k({\mathfrak P})$ is the residue field of $\mathfrak P$). In this case intersections commute with specialization,
and so this reduction is obtained as the image of the $d$-th self-intersection $\bar D^d\in 
A_0(V_{\overline{k({\mathfrak P})}})$ of $\bar D\in  {\rm Pic}^0(V_{\overline{k({\mathfrak P})}})$.
However $\bar D$ is a torsion element in a divisible group and the 
intersection pairing is bilinear. It follows that $\bar D^d=0$ in $A_0(V_{\overline{k({\mathfrak P})}})$.

}
\end{Remark}

\begin{Remark} \label{Poincare}
{\rm 
Our construction in this paragraph actually gives a map
\begin{equation*}
{\rm H}^1(Y_\Q,\Z/p)\xrightarrow{\ \ \ } \Q(\zeta_p)^*/(\Q(\zeta_p)^*)^p={\rm H}^1(\Q(\zeta_p), \Z/p).
\end{equation*}
It is not hard to see that this map is a group homomorphism.
When $h_\Q: Y_\Q\to \Spec(\Q)$ is smooth it can also be derived 
from 
\'etale duality for the smooth projective morphism $h_\Q$. We will leave the details
to the interested reader.
}
\end{Remark}

\bigskip
\section{Galois module structure}\label{Galoismod}

In this section, we complete the proofs of the
main results stated in the introduction.
We begin with some preliminaries.

\subsection{}  \label{10a}
Let $B$ be an associative ring
with unit, which we assume is left Noetherian.
Most of the time we will take $B$ to be the group
ring $R[G]$ of a finite group $G$ with $R$ commutative Noetherian.
We denote by $\Gr_0(B)$ (resp. $\Kr_0(B)$)
the Grothendieck group of finitely generated (resp. finitely generated projective)
left $B$-modules, and by $\Gr^{\rm red}_0(B)$, resp. $\Kr_0^{\rm red}(B)$,
the quotient of $\Gr_0(B)$, resp. $\Kr_0(B)$, by the subgroup generated by the
class $[B]$ of the free $B$-module $B$. We will denote by $a^{\rm red}$ the image
of the element $a$ of $\Gr_0(B)$, resp. $\Kr_0(B)$,
in $\Gr_0^{\rm red}(B)$, resp. $\Kr_0^{\rm red}(B)$.
Denote by $c: \Kr_0(B)\to \Gr_0(B)$, 
$c^{\rm red}: \Kr^{\rm red}_0(B)\to \Gr^{\rm red}_0(B)$,
the natural forgetful  homomorphisms.
If $B$ is commutative, we will denote by $\Pic(B)$
the Picard group of $B$.  
If $B$ is commutative, a finitely generated $B$-module is projective
if and only if it is
locally free. Taking highest exterior powers of
locally free $B$-modules defines in this case a
group homomorphism
$$
i: \Kr_0^{\rm red}(B)\to \Pic(B).
$$

We will denote by $D^+(B)$ the derived category
of the homotopy category of complexes of left $B$-modules
which are bounded below. Recall that a complex $C^\bullet$
in $D^+(B)$ is called ``perfect", if it is isomorphic in $D^+(B)$
to a bounded complex $P^\bullet$ of
finitely generated projective
left $B$-modules. Then
the element
\begin{equation}\label{chiPerf}
\chi(C^\bullet)=\sum_i(-1)^i[P^i]\in \Kr_0(B)
\end{equation}
depends only on the isomorphism class of $C^\bullet$ in $D^+(B)$.

In what follows, all the schemes will be
separated and Noetherian.
For a scheme $Y$, we will denote by
$\Gr_0(Y)$, resp. $\Kr_0(Y)$, the Grothendieck group of
coherent, resp. coherent locally free,
sheaves of $\O_Y$-modules, and by $\Pic(Y)$ the Picard group of $Y$.
If $Y$ is regular and has an ample invertible sheaf, then we have
$\Kr_0(Y)\simeq \Gr_0(Y)$ ([SGA6] IV).
If $Y=\Spec(B)$ is affine, we will identify quasi-coherent sheaves of $\O_Y$-modules
on $Y$ with $B$-modules. This gives natural
identifications $\Gr_0(B)=\Gr_0(\Spec(B))$, $\Kr_0(B)=\Kr_0(\Spec(B))$ and
$\Pic(B)=\Pic(\Spec(B))$.

\subsection{} \label{Euler} Now suppose that $S=\Spec(R)$, $R$ regular Noetherian
and that $\pi: X\to Y$ is a (right) torsor for a finite group $G$ with $h:Y\to S$ projective.  
Then $f:=\pi\cdot h: X\to S$ is also projective with a free (right) $G$-action. Denote by $\Gr_0(G,X)$,
resp. $\Kr_0(G,X)$, the Grothendieck group of $G$-equivariant coherent, resp. locally free coherent,
sheaves on $X$. By descent (see \S \ref{torsors2}) pulling back along the finite \'etale morphism $\pi$
gives isomorphisms
\begin{equation}
\pi^*: \Gr_0(Y)\xrightarrow{\sim} \Gr_0(G,X),\quad \pi^*: \Kr_0(Y)\xrightarrow{\sim} \Kr_0(G,X).
\end{equation}
Recall also (cf. \S\ref{torsors2}) that if $\F$ is a $G$-equivariant coherent sheaf on $X$
we can view
$\pi_*(\F)$ as a coherent sheaf of $\O_Y[G]$-modules 
on $Y$: If $V=\Spec(C)$ is an open affine subscheme of $Y$
and $U=\pi^{-1}(V)$, then the sections $(\pi_*(\F))(V)=\F(U)$ are naturally
a left $C[G]$-module.  We can then consider the right derived image
$\R \Gamma(Y,\pi_*(\F))$; this is a complex in $D^+(R[G])$ 
which computes the cohomology of $\pi_*(\F)$
([SGA6] III \S 2, IV \S 2).
We may give a complex isomorphic to 
$\R \Gamma(Y,\pi_*(\F))$ by taking the (bounded) \v{C}ech complex obtained by considering the
sections of $\pi_*(\F)$ on the intersections
of the sets in a finite affine cover $\{V_i\}_i$ of $Y$.

\begin{thm} \label{perfect} \ ([CEPT1] Theorem 8.3; see also [C], [CE].) \ \
Assume as above that  $\pi: X\to Y$ is a $G$-torsor, $h:Y\to S$ projective and
$S=\Spec(R)$ is
regular Noetherian. Let $\F$ be a $G$-equivariant coherent
sheaf on $X$. Then the complex
$\R \Gamma(Y,\pi_*(\F))$ in $D^+(R[G])$ is perfect. The elements
\begin{equation*}
\chi_f^P(\F):=\chi(\R\Gamma(Y,\pi_*(\F)))\in \Kr_0(R[G])
\end{equation*}
define a group homomorphism
(the equivariant projective Euler characteristic)
\begin{equation*}
\chi_f^P:\Gr_0(G, X)\to \Kr_0(R[G]).
\end{equation*}
\end{thm}

We also set $\bar\chi_f^P(\F):=(\chi_f^P(\F))^{\rm red}\in \Kr_0^{\rm red}(R[G])$.  
Note that $\bar\chi_f^P(\F) \in \Kr_0^{\rm red}(R[G])$ is the obstruction
for the complex $\R \Gamma(Y,\pi_*(\F))$ to be isomorphic in $D^+(R[G])$ to a bounded
complex of finitely generated {\sl free} $R[G]$-modules.
When  $f:X\to S$ is fixed,
we will usually write $\chi^P$, $\bar\chi^P$ instead of $\chi^P_f$, $\bar\chi_f^P$.

Suppose now that $R$ is the ring of integers of a number field. 
We will say that a $R[G]$-module $M$ is locally free if for
each prime ideal $\PP$ of $R$, $M\otimes_RR_\PP$ is a free $R_\PP[G]$-module.
By a theorem of Swan the notions of projective and locally free
coincide for finitely generated $R[G]$-modules ([Sw]). 
Therefore, we may identify $\Kr_0^{\rm red}(R[G])$ with the 
class group of finitely generated locally free  $R[G]$-modules and write ${\rm Cl}(R[G])$ instead of
$\Kr_0^{\rm red}(R[G])$.  Also, let us mention that in this case,
$\chi_f^P$ coincides with the homomorphism $f^{CT}_* :\Gr_0(G, X)\to {\rm CT}(R[G])=\Kr_0(R[G])$ of [C] and [CE]
(here ${\rm CT}(R[G])$ denotes the Grothendieck group of finitely generated
$R[G]$-modules which are  cohomologically trivial as $G$-modules; see [C] and [CEPT1] Lemma 8.5).

\subsection{ } \label{relationS}

In what follows, we continue with the general assumptions and notations of
the previous paragraph. In addition, we will assume that $G$ is commutative
and that the morphisms $f: X\to S$ and $h: Y=X/G\to S$ are flat.
We denote by 
$
 h_{G^D_S}:  Y\times_SG^D_S\to G^D_S
$
the base change of $h$. In what follows, we will omit the subscript 
$S$ from $G^D_S$. The morphism 
$h_{G^D}$ is also flat and projective.
For any coherent sheaf $\Hh$ of 
$\O_{Y\times_SG^D}$-modules on $Y\times_SG^D$, we can consider the 
complex ${\bf R}{h_{G^D}}_*(\Hh)$ in the derived category $D^+(\O_{G^D})$ of the homotopy category of
complexes of sheaves of $\O_{G^D}$-modules on $G^D$ which are bounded below. 
In fact, if $\Hh$ is locally free, then ${\bf R}{h_{G^D} }_*(\Hh)$ is perfect ([SGA6] III)
and we can consider the determinant of cohomology  $\det{\bf R}{h_{G^D}}_*(\Hh)$
(see \S \ref{det})
and the Euler characteristic
$\chi({\bf R}{h_{G^D}}_*(\Hh))\in \Kr_0(R[G])$. (Notice that, by [KM] Prop. 4,
over the affine scheme  $G^D$, a complex is perfect
in the sense of [SGA6] if and only if it is {\sl globally} quasi-isomorphic to
a bounded complex of sheaves associated to finitely generated locally free $R[G]$-modules;
this allows us to define $\chi({\bf R}{h_{G^D}}_*(\Hh))\in \Kr_0(R[G])$ by a formula similar
to (\ref{chiPerf}).)

As in \S \ref{2d} we may think
of the locally free coherent sheaf $\pi_*(\O_X)$ of $ \O_Y[G]$-modules on $Y$ 
as an invertible sheaf of $\O_{Y\times_SG^D}$-modules on $Y\times_SG^D$.
We can construct a bounded complex isomorphic to ${\bf R}{h_{G^D}}_*(\pi_*(\O_X))$
by using the \v{C}ech construction associated
to $\pi_*(\O_X)$ and the finite affine cover 
$\{V_i\times_SG^D\}_i$ of $Y\times_S G^D$, where $\{V_i\}_i$ is a finite affine cover of $Y$. 
It follows that the complexes ${\bf R}{h_{G^D}}_*(\pi_*(\O_X))$ and 
$\R \Gamma(Y,\pi_*(\O_X))$
(see \S \ref{Euler}) are  isomorphic
(here again we identify $R[G]$-modules with the corresponding $\O_{G^D}$-sheaves). Both complexes
$\R \Gamma(Y,\pi_*(\O_X))$ and ${\bf R}{h_{G^D}}_*(\pi_*
(\O_X))$ are perfect (by Theorem \ref{perfect} and [SGA6] III). We obtain
\begin{equation*}
i(\chi(\R \Gamma(Y,\pi_*(\O_X)))^{\rm red})
=[\det{\bf R}{h_{G^D}}_*(\pi_*(\O_X))]\in {\rm Pic}(G^D_S)={\rm Pic}(R[G])
\end{equation*}
and therefore
\begin{equation} \label{i}
i(\bar\chi_f^P(\O_X))=[\det{\bf R}{h_{G^D}}_*(\pi_*(\O_X))].
\end{equation}

If $R$ is of Krull dimension $1$, which is the main case we are interested in,
the same is true for $R[G]$. Then by [BM] Cor. 3.5, $i: \Kr^{\rm red}_0(R[G])
\to \Pic(R[G])=\Pic(G^D_S)$ is an isomorphism. In this case, we
will use $i$ to identify the two groups; this should not cause any confusion.
We can then write
\begin{equation} \label{i2}
\bar\chi_f^P(\O_X)=[\det{\bf R}{h_{G^D}}_*(\pi_*(\O_X))].
\end{equation}

\subsection{} \label{paralb} We continue with the general assumptions of the
previous paragraph. In particular, $G$ is a commutative group
and  $\pi: X\to Y$ is a $G$-torsor with $h: Y\to S=\Spec(R)$ projective and flat
($R$ is regular and Noetherian). Recall that in this case, $\pi_*(\O_X)$ is an invertible sheaf
over $Y\times_SG^D$ (\S \ref{2d}). The corresponding $\Gm$-torsor supports
a commutative extension 
\begin{equation}\label{alb1}
1\to \Gm_Y\to E\to G^D_Y\to 1
\end{equation}
of group schemes over $Y$ (see Remark \ref{rem22}).

\begin{Definition} \label{alba}
{ We will say that the $G$-torsor $\pi: X\to Y$
is of Albanese type over $S$ if there is a smooth  
commutative group scheme of finite type $A\to S$ with {\sl connected} fibers, a commutative extension
\begin{equation}\label{alb2}
1\to \Gm_Y\to {\mathcal E}\to A_Y\to 1
\end{equation}
of group schemes over $Y$ and a group scheme homomorphism $\phi: G^D_S\to A$ over $S$ 
such that $E\simeq ({\rm id}_Y\times_S \phi)^*({\mathcal E})$ as group 
scheme extensions. }
\end{Definition}

\begin{Remarknumb}\label{remalb}
{\rm  a) If $\PP$ is the invertible sheaf over $A_Y=Y\times_SA$ which corresponds to
the $\Gm_Y$-torsor ${\mathcal E}$, then our condition implies that 
$\pi_*(\O_X)\simeq ({\rm id}_Y\times_S \phi)^*({\PP})$ as invertible shaves on $G^D_Y$.

b) Let $\pi: X\to Y$ be a $G$-torsor of Albanese type
and suppose that $i:Z\to Y$ is a projective morphism such that $h\cdot i: Z\to S$
is flat. Then the $G$-torsor $\pi_Z: X\times_YZ\to Z$ obtained by base change is also of Albanese type.

c) Suppose 
$S=\Spec(k)$ with $k$ a field of characteristic prime to the order of $G$
and that $Y$ is smooth 
and projective over $\Spec(k)$. Assume that $Y$ contains a $k$-rational point.

The following construction motivates our terminology:
Take $A={\rm Pic}^0(Y)$
(the Picard abelian variety of $Y$) and suppose we have a group scheme immersion
$\phi: G^D\to A$. The $k$-rational point on $Y$
provides us with a morphism $Y\to {\rm Alb}(Y)$ to the Albanese abelian variety of
$Y$.  Set $B=A/G^D$ for the quotient abelian variety and consider
\begin{equation*}
0\to G\to B^{\rm dual}\to A^{\rm dual}\to 0
\end{equation*}
given by the  isogeny dual to $A\to B$. 
The canonical duality $A^{\rm dual}\simeq {\rm Alb}(Y)$ 
between the Albanese and  Picard varieties of $Y$
allows us to view this exact sequence as a $G$-torsor over ${\rm Alb}(Y)$.
We can restrict this along $Y\to {\rm Alb}(Y)$ to obtain a $G$-torsor $X\to Y$.
Such torsors are called of Albanese type by Lang [La]. 
One can now see (using for example [Mu] \S 15) that $X\to Y$ is 
also of Albanese type 
according to our definition above. The required extension 
\begin{equation}\label{alb2}
1\to \Gm_Y\to {\mathcal E}\to A_Y\to 1
\end{equation} 
is given as the pull-back along $Y\times_S A\to {\rm Alb}(Y)\times_S A\simeq A^{\rm dual}\times_S A$
of the extension over $A_{A^{\rm dual}}=A^{\rm dual}\times_S A $ given by the Poincare invertible sheaf
with its biextension structure.  }
\end{Remarknumb}

\begin{thm} \label{albcond} 
{  Suppose $S=\Spec(k)$ with $k$ a field of characteristic prime to the order of $G$
and that $Y$ is smooth and projective over $\Spec(k)$. Assume that 
$Y$ contains a $k$-rational point and let $X\to Y$ be a  $G$-torsor 
such that $X$ is geometrically connected. Then $X\to Y$ is 
of Albanese type if and only if all the torsion invertible sheaves $\O_{X,\chi}$,
with $\chi$ running over all characters of $G$ (see \S \ref{2d}), 
correspond to divisors on $Y\times_k k(\chi)$ which are algebraically trivial,
or equivalently if and only if it is obtained by pulling back a $G$-isogeny
of the Albanese abelian variety ${\rm Alb}(Y)$ along $Y\to {\rm Alb}(Y)$.}
\end{thm}

\begin{Proof}
This follows from the construction and properties
of the Picard and Albanese abelian varieties of $Y$ as in 
Remark \ref{remalb} (c) above.\endproof
\end{Proof}

\begin{Remarknumb}\label{remalb2} 
{\rm  Assume that $X\to Y$ is an abelian $G$-torsor over the ring of integers $\O_K$ of the number field $K$ such that the 
generic fiber $X_K\to Y_K$ is of Albanese type. Assume that $Y_K$ is smooth
and that $Y\to \Spec(\O_K)$ has a section. Let ${\mathcal A}^0$ be the connected component
of the Neron model ${\mathcal A}$ of ${\rm Pic}^0(Y_K)$. It would be interesting to find conditions 
under which $X\to Y$ is of Albanese type over $\Spec(\O_K)$ with $A={\mathcal A}^0$; this is always the case when $Y_K$
is a curve of genus $\geq 2$ and $Y$ is regular by [Ra] \S 8. 
}
\end{Remarknumb}

Let \hbox{$\pi: X\to Y$} be an (abelian) $G$-torsor
of Albanese type over $S=\Spec(R)$. Suppose now that the residue fields of all the 
generic points of $S$
are perfect.  Denote by $h_A: Y\times_SA\to A$ the base change of $h: Y\to S$
and let $\M_A=\det({\bf R}{h_A}_*(\PP))$, $\M= \det{\bf R}{h_{G^D}}_*(\pi_*(\O_X))$
(invertible sheaves on $A$ and $G^D_S$ respectively). 
Since $\pi_*(\O_X)\simeq ({\rm id}_Y\times_S \phi)^*({\PP})$  we
then have $\M\simeq \phi^*(\M_A) $. By [Br] Proposition 2.4 and the discussion before
it (see also [SGA7I] VIII), the invertible sheaf $\L_A$ over $A$ supports a 
canonical cubic structure (as in Definition \ref{cubicdef} with $n=3$). 
This statement is of course an extension  to this situation of the classical theorem of the cube
for line bundles over abelian varieties.
(Notice that Breen uses a slightly different definition of cubic structure;
see Remark \ref{32re} (a) or [Br] 2.8. However, it is not hard to see 
that this does not affect the truth of our statement.) Using the functoriality of cubic structures, 
we can conclude that $\M\simeq \phi^*(\M_A)$ also supports a cubic structure.  

\subsection{} We can now complete the proof of our main results 
(Theorems \ref{main1}, \ref{main2} and \ref{main3}).

As in these statements we assume that $G$ is a finite group
and $\pi: X\to Y$ is
a $G$-torsor with $h: Y\to \Spec(\Z)$ projective and flat of 
relative dimension $d$. We set $f=h\cdot \pi$ and
$S=\Spec(\Z)$.

A standard argument using Noetherian induction shows:

\begin{lemma}\label{generate}
The Grothendieck group $\Gr_0(Y)$ is generated by the classes
$[i_*\O_Z]$ with $i:Z\hookrightarrow Y$ an integral subscheme of $Y$.
There are two possibilities for such a $Z$: Either the morphism $Z\to \Spec(\Z)$ 
is  flat of relative dimension $d'\leq d$, or it factors through 
$\Spec(\fp)\to \Spec(\Z)$ for some prime $p$.\endproof
\end{lemma}

Recall the descent
isomorphism $\pi^*: G_0(Y)\xrightarrow{\sim} G_0(G,X)$. Let $i: Z\hookrightarrow Y$
be as in Lemma \ref{generate} and 
consider the $G$-torsor $\pi_Z: X\times_YZ\to Z$ obtained by pulling back
$\pi$ along the closed immersion $i$. Denote by $f': X\times_YZ\to \Spec(\Z)$
the structure morphism. Then it follows from the construction
of the projective Euler characteristic that
\begin{equation}\label{gen1}
\chi^P_{f}(\pi^*(i_*\O_Z))=\chi^P_{f'}(\O_{X\times_YZ}).
\end{equation}
If $Z\to \Spec(\Z)$ factors through $\Spec(\fp)\to \Spec(\Z)$
then it follows from the main theorem of   [Na]  that $\bar\chi^P_{f'}(\O_{X\times_YZ})=0$
(see also [P] 4.b Remark 3).
On the other hand, we notice that if $d'\leq d$ we have $M_{d'+1}(G)|M_{d+1}(G)$,  
$M'_{d'+1}(G)|M'_{d+1}(G)$ and also $C_{d'+1}(G)|C_{d+1}(G)$.
Using Lemma \ref{generate}, the additivity of the Euler characteristic
$\bar\chi^P_f$ (Theorem \ref{perfect}) and (\ref{gen1}), we can now see that
the proofs of  Theorems \ref{main1}, \ref{main2} are reduced to the case that 
$\F$ is the structure sheaf $\O_X=\pi^*\O_Y$ and $Y$ is integral. A similar argument 
using $(d'+1)!!|(d+1)!!$ for $d'\leq d$ together with Remark \ref{remalb} (b) 
shows that the proof of  Theorem \ref{main3} is also reduced to the case
that $\F$ is the structure sheaf of $X$ and $Y$ is integral.

Let us first discuss the 
completion of the proofs of Theorems \ref{main1} and \ref{main2}.
This borrows heavily from the arguments in [P]:
Noetherian induction
and an argument as in the proof of [P] Proposition 4.4 (a) shows that $(\#G)^{d+1}\cdot \bar\chi^P(\O_X)=0$.
Suppose first that $G$ is abelian. In this case, by Theorem \ref{detC},
the invertible sheaf $\det{\bf R}{h_{G^D}}_*(\pi_*(\O_X))^{\otimes \kappa}$ 
over $G^D_S$ supports a $d+2$-cubic structure. Hence, it follows from  \ref{i2} and Theorem \ref{Ncube}
that $\kappa\cdot \bar\chi^P(\O_X)$ is annihilated by $M_{d+1}(G)$. Since 
$(\#G)^{d+1}\cdot \bar\chi^P(\O_X)=0$ this
gives Theorem \ref{main1}  (b) for $\F=\O_X$
when $G$ is abelian. In fact, we can see that when $G$ is abelian
and $\F=\O_X$, the conclusion of Theorem \ref{main1} (b) 
with $\ep(G)$ replaced by $\ep(G,Y)$ holds true.

The fact that $(\#G)^{d+1}\cdot \bar\chi^P(\O_X)=0$
together with the ``localization" argument in the proof of [P] Proposition 4.5 shows that
the proofs of Theorems \ref{main1} and \ref{main2} in general 
reduce to the case that $G$ is an $l$-group, $l$ prime. 

Theorem \ref{main1} (b) now follows
from the abelian case just explained above. To show Theorem \ref{main1} (a) observe
that the argument in [P] p. 215-216 allows us to reduce the case of an $l$-group to that
of the case of a ``basic" $l$-group with $l$ odd, i.e to the case of a cyclic group
of odd prime order $l$. Then part (a) follows once again by the abelian case.
This completes the proof of Theorem \ref{main1}. 

The same argument from [P] also shows that the proof of Theorem \ref{main2} 
can be reduced to the case of a ``basic" $l$-group with $l$ odd, i.e to the case of a cyclic group
of odd prime order $l$. Then Theorem \ref{main2} follows from 
Theorem \ref{detC} and Theorem \ref{KVtrivial}.
\endproof
\medskip

Let us now discuss the proof of Theorem \ref{main3}. Recall that it is enough to deal
with the case that $\F=\O_X$ and $Y$ is integral. It is also enough to assume
that $d\geq 1$. Let us first discuss part (a):
Since $\pi: X\to Y$ is of Albanese type, by \S\ref{paralb}, the invertible sheaf
$\M:=\det{\bf R}{h_{G^D}}_*(\pi_*(\O_X))$ supports a $3$-cubic structure. Hence, by 
Corollary \ref{NcubeCor}, $\M$ is trivial. It follows from  (\ref{i2}) that  
$\bar\chi^P(\O_X)=0$ in ${\rm Cl}(\Z[G])$. This completes the proof of part (a).

It remains to prove Theorem \ref{main3} (b):

For simplicity, we set $T=X_\Q$, $U=Y_\Q$. We are 
assuming that $T\to U$ is of Albanese type
and so there is a commutative group scheme $A$ and an  
extension ${\mathcal E}$ as in
Definition \ref{alba}. We denote by $\PP$ the invertible 
sheaf over $A_U=U\times_\Q A$ given by
${\mathcal E}$. By Theorem \ref{detC}, the invertible sheaf 
$\M^{\otimes 2}=\det{\bf R}{h_{G^D}}_*(\pi_*(\O_X))^{\otimes 2}$ 
over $G^D_S$ supports a $d+2$-cubic structure $\xi$. We can
apply the same argument to the invertible sheaf $\M_A=\det({\bf R}{h_A}_*(\PP))$
over $A$: Here we start from the extension
\begin{equation}
1\to \Gm_U\to {\mathcal E}\to A_U\to 1
\end{equation}
and the corresponding additive functor  
$A\to {\mathcal {PIC}}(U)$; in this case, we have to work using the
fpqc topology. We obtain that $\M_A^{\otimes 2}$ has a canonical
$d+2$-cubic structure $\xi_A$ over $A$. Recall that the extension ${\mathcal E}$ pulls back
to the extension 
\begin{equation}
1\to \Gm_U\to E\to G^D_U\to 1
\end{equation}
using $\phi: G^D_\Q=G^D_{\Spec(\Q)}\to A$. This implies that
the $d+2$-cubic structure $\xi_\Q$ on the generic fiber $\M^{\otimes 2}_\Q$ 
(an invertible sheaf over $G^D_\Q=\Spec(\Q[G])$) is obtained by
pulling back $\xi_A$ along $\phi$; so there is an isomorphism 
of $d+2$-cubic structures $\xi_\Q\simeq \phi^*(\xi_A)$. 
However, as a corollary of the theorem of the cube on $A$ (see 
the argument at the end of the previous 
paragraph), $\M_A$ and hence also $\M_A^{\otimes 2}$ has a canonical
$3$-cubic structure. This induces a $d+2$-cubic structure $\xi_A'$
on $\M_A^{\otimes 2}$ using the construction of Lemma \ref{61}. 
We claim that this agrees with the $d+2$-cubic structure $\xi_A$:
To see this observe that the composition $\xi_A^{-1}\cdot \xi'_A$ is given by an invertible
regular function $c$ on $A^{d+2}$ which gives a $d+2$-cubic structure on the trivial invertible
sheaf. By \S \ref{various} $(c0')$  we have $c(a_1,\ldots, a_{d+2})=1$ if one of the $a_i$
is zero. By a classical lemma of Rosenlicht ([SGA7] VIII 4.1) this implies that $c$
is equal to $1$. Hence $\xi'_A=\xi_A$. Since $\xi_\Q\simeq \phi^*(\xi_A)$ and $\xi'_A=\xi_A$
comes from a $3$-cubic structure (using Lemma \ref{61}) we conclude 
that the same is true for $\xi_\Q$.
Now we will consider the Taylor expansion of Corollary \ref{flat2} for $(\M^{\otimes 2}, \xi)$
and $(\M^{\otimes 2}_\Q, \xi_\Q)$ (we take $n=d+1$). For simplicity,
set $\L=\M^{\otimes 2}$, $\L_\Q=\M^{\otimes 2}_\Q$.
By repeated application of Lemma \ref{62} and Proposition \ref{flat}
we see that since $(\L_\Q, \xi_\Q)$ is a $d+2$-cubic structure which comes from a 
$3$-cubic structure the $(d+1-i)$-extensions 
$E(\L^{(i)}_\Q, \xi^{(i)}_\Q)$ of $G^D_\Q$ by $\Gm$ are trivial when $d+1-i\geq 3$.
Since the construction of $(\L^{(i)}, \xi^{(i)})$ commutes with base change from $\Z$ to $\Q$,
this now implies that the $(d+1-i)$-extensions $E(\L^{(i)}, \xi^{(i)})$ of $G^D$ by $\Gm$
(over $S=\Spec(\Z)$) become trivial after base change to the generic fiber $\Spec(\Q)$
when $d+1-i\geq 3$. We claim that this implies that all these multiextensions are trivial. 
To see this let us consider the commutative diagram
\begin{equation}
\begin{CD}
n{\hbox{\small -}}{\rm Ext}^1(G^D,\Gm)
  @>{\sim}>> (n-1){\hbox{\small -}}{\rm Ext}^1(G^D,G) @> t>> {\rm H}^1((G^D)^{n-1},G)\\
@VVV @VVV @VVV\\
n{\hbox{\small -}}{\rm Ext}^1(G^D_\Q,\Gm) @>{\sim}>> (n-1){\hbox{\small -}}{\rm Ext}^1 (G^D_\Q,G)@> t_\Q>> {\rm H}^1  ((G^D_\Q)^{n-1},G).\\
\end{CD}
\end{equation}
Here the isomorphisms on the left side are given by (\ref{nExt3}). The homomorphism
$t$ is the forgetful map; the bottom row is obtained by the constructions 
that give the first row but performed over $\Spec(\Q)$. The vertical arrows are the base change 
homomorphisms. By Lemma \ref{exact} $t$ is injective (the argument immediately extends to the case
that $G$ is not cyclic).
The right vertical arrow is also injective. We conclude that
  the left vertical arrow is injective; this  implies the desired result:
the $(d+1-i)$-extensions $E(\L^{(i)}, \xi^{(i)})$ of $G^D$ by $\Gm$
(over $S=\Spec(\Z)$) are trivial when $d+1-i\geq 3$. In addition, by  Remark \ref{low} (a)
the $(d+1-i)$-extensions $E(\L^{(i)}, \xi^{(i)})$ are trivial when $d+1-i<3$.
This, together with ${\rm Pic}(S)=(0)$, implies that all the terms in the Taylor
expansion of $\L^{\otimes (d+1)!!}$ are trivial. Hence, by (\ref{i2}), we have $2((d+1)!!)\cdot \bar\chi^P(\O_X)=0$
in ${\rm Cl}(\Z[G])$. This concludes the proof of Theorem \ref{main3}.\endproof
\smallskip

\begin{Remarknumb} \label{albresext}
{\rm Suppose that $Y$ is normal of dimension $\geq 2$ and denote by 
\begin{equation}
Y\xrightarrow{h'}S'=\Spec(\O_K)\xrightarrow{\phi} S=\Spec(\Z)
\end{equation}
the ``Stein factorization" of $h$ where $\O_K$ is the ring of integers of the number field $K$. 
Assume that the abelian $G$-torsor $X\to Y$ is of Albanese type {\sl as a torsor over $S'$}.
Then the proof of Theorem \ref{main3} (a) above can be extended to show that $\bar\chi^P(\F)=0$
in $\Cl(\Z[G])$. Similarly, if  $X_K\to Y_K$ is of Albanese type 
over $\Spec(K)$ we can obtain that the conclusion of Theorem \ref{main3} (b) 
holds i.e that $2((d+1)!!)\cdot \bar\chi^P(\F)=0$ for all $G$-equivariant $\F$ on $X$. To see this it is enough to combine
the arguments for  $h':Y\to S'$ above with the following observations:

i) There is a canonical isomorphism
$$
\det{\bf R}{h_{G^D}}_*(\pi_*(\O_X))\xrightarrow{\sim}\det(\phi_*(\det{\bf R}{h'_{G^D}}_*(\pi_*(\O_X)))).
$$

ii) For any invertible sheaf $\M'$ on $G^D_{S'}$ and $n\geq 2$ 
the additivity of the Norm (\ref{rem23} (b)) gives  a canonical isomorphism
$$
{\mathfrak N}: {\rm Norm}_{S'/S}(\Theta_{n}(\M'))\xrightarrow{\sim} \Theta_{n}(\det(\phi_*(\M'))).
$$

iii) If $(\M', \xi')$ is an invertible sheaf with an $n$-cubic structure
over $G^D_{S'}$ then the isomorphism $\xi:={\mathfrak N}\cdot {\rm Norm}_{S'/S}(\xi')$ 
defines an $n$-cubic structure on $\M:=\det(\phi_*(\M'))$ over $G^D_S$.

We will leave the details to the reader.
}
\end{Remarknumb}

\medskip

Let us finish by explaining how we can deduce the other
results stated in the introduction. Corollary \ref{1.2}
follows directly from Theorem \ref{main1} (b), Theorem \ref{main3} (a) and the definition of the projective Euler
characteristic (Theorem \ref{perfect}). Notice that the improvement ($\ep(G,Y)$ in place
of $\ep(G)$) in Theorem \ref{main1} (b) and Corollary \ref{1.2} 
in the case that $\F=\O_X$ which was stated in the introduction
follows now directly from the argument in the proof of Theorem \ref{main1} (b) above. Corollary \ref{1.3}
follows from Theorem \ref{main1} (b) and Theorem \ref{main3} (a) in exactly the 
same way as explained in the proof of [P] Cor. 5.3. Finally,
if $G=\Z/p\Z$ with $p$ odd and $p>d+1$ then the equality 
\begin{equation*}
\bar\chi^P(\O_X)=\sum_{i=1}^{d+1}R^{(i)}(t_i(X/Y)) 
\end{equation*}
of the introduction follows from (\ref{i2}) and (\ref{84}) after setting
$\L=\det{\bf R}{h_{G^D}}_*(\pi_*(\O_X))^{\otimes \kappa}$ and
$t_i(X/Y)=\kappa^{-1}t_i(\L, \xi)\in \Hom((C(p)/p)^{(1-i)}, \Z/p)$. 
This combined with Remark \ref{explicitRe} (a)
gives 
\begin{equation}
t_{d+1}(X/Y)=\frac{1}{\kappa \cdot (d+1)!} [T(X/Y)]
\end{equation}
where $[T(X/Y)]$ is the element  
that corresponds to the unramified $\Z/p$-extension 
of $\Q(\zeta_p)$ given in Proposition \ref{explicit}.
\endproof

\smallskip
\begin{Remarknumb}\label{albtri}
{\rm It follows from the proof of Theorem \ref{main3} and 
the construction of the elements $t_i(X/Y)=\kappa^{-1}t_i(\L, \xi)$
that these are trivial for $i\geq 3$ when $Y$ is integral and $X_\Q\to Y_\Q$
is of Albanese type. }
\end{Remarknumb}
\bigskip
\medskip
\vfill\eject

\centerline{\sc APPENDIX}
\bigskip

In this appendix, we show how a well-known argument 
due to Godeaux and Serre ([Se] \S 20) combined with
an ``arithmetic" version of Bertini's theorem 
(based on the theorem of Rumely on the existence 
of integral points [Ru], [MB2]) allows us to construct ``geometric" $G$-torsors $\pi: X\to Y$
with $Y$ regular and $Y\to \Spec(\Z)$ projective and flat of relative dimension $d$
for any finite group $G$ and any integer $d\geq 1$.
In fact, for the $G$-torsors that we construct $Y\to \Spec(\Z)$ factors through 
a smooth morphism $Y\to \Spec(\O_K)$ where $\O_K$ is the ring of integers
of a number field $K$.

To explain this set $L=\Z[G]$ and for an integer $r\geq 1$ we denote by 
$S(r)=\oplus_{m\geq 0}S(r)_m={\rm Sym}_{\Z}(L^{\oplus r})$ 
the corresponding (graded) symmetric algebra with $G$-action.
We set $X(r)={\bf Proj}(S(r))$; this is then a projective space of (relative) 
dimension $s=r\cdot (\# G)-1$ over $\Spec(\Z)$
that supports a linear action of $G$. The quotient $Y(r)=X(r)/G$ is also
a projective scheme: We can see that, for sufficiently large integers $k$,
$$
Y(r)\simeq {\bf Proj}\bigl(\bigoplus_{m\geq 0}(S(r)_{mk})^G\bigl) ,
$$
and the graded algebra $\bigoplus_{m\geq 0}(S(r)_{mk})^G$ is generated by 
the free $\Z$-module $M:=S(r)_{k}^G$.  Set ${\bf P}(M)={\bf Proj}({\rm Sym}(M))$;
then
$Y(r)$ is a closed subscheme of ${\bf P}(M)$.
(The reader can consult [Se] \S 20 for the details of the argument in the corresponding situation 
over an algebraically closed field; the same argument readily applies to our case.)
Denote by $\pi(r): X(r)\to Y(r)$ the quotient morphism. Let $B(r)$ be the 
closed subscheme of $X(r)$ consisting of points with
non-trivial inertia subgroups and set $b(r)$ for the (reduced) image $\pi(r)(B(r))$. 
The group $G$ acts freely on the open subscheme $U(r):=X(r)-B(r)$,
the morphism $\pi(r): U(r)\to V(r):=Y(r)-b(r)$
is a $G$-torsor and $V(r)\to \Spec(\Z)$ is smooth of relative 
dimension $s$. One can now observe ([Se]) that 
for all sufficiently large $r$,
each fiber of $b(r)\to \Spec(\Z)$ has codimension $>d$ in the corresponding fiber 
of $Y(r)\to \Spec(\Z)$. 

Now let us
consider the dual projective space ${\bf P}(M^\vee)$
parametrizing hyperplane sections of ${\bf P}(M)$. 
Denote by ${\bf H}\subset {\bf P}(M)\times {\bf P}(M^\vee)$ 
the universal hyperplane section. 
Let us set $Q=({\bf P}(M^\vee))^{s-d}$. For $\phi: T\to Q$, $i=1,\ldots, s-d$,
let ${\bf H}^i_{\phi}$  be the  hyperplane in 
${\bf P}(M)\times T$ that corresponds to ${\rm pr}_i\cdot \phi: T\to {\bf P}(M^\vee)$
(this is the Cartesian 
product
\begin{equation*}
\begin{CD}
{\bf H}^i_{\phi}@>>> {\bf H}\\
@VVV @VVV\\
{\bf P}(M)\times T@>{\rm id}\times {\rm pr}_i\cdot \phi>> 
{\bf P}(M)\times {\bf P}(M^\vee).)\\
\end{CD}
\end{equation*}
Set $Y(r)_{\phi: T\to Q}$ for the scheme theoretic intersection
$$
(Y(r)\times T)\cap ({\bf H}^1_\phi\cap \cdots \cap {\bf H}^{s-d}_\phi)
$$
in ${\bf P}(M)\times T$; this is a scheme over $T$. Let $V_1$ be the subset of $x\in Q$ for which 
$$
Y(r)_{\Spec(k(x))\to Q}\subset Y(r)\times \Spec(k(x))
$$ 
{\sl does not} intersects $b(r)\times \Spec(k(x))$. Let $V_2$ 
be the subset of $x\in Q$ for which the projection $Y(r)_{{\rm id}: Q\to Q}\to Q$
is flat at all points that lie over $x$. Finally, 
let $V_3$ be the subset of $x\in Q$ for which
$Y(r)_{\Spec(k(x))\to Q}$ is smooth  over $\Spec(k(x))$
of dimension $d$. Set $V=V_1\cup V_2\cup V_3$;
this is a constructible subset of $Q$.
The proof of the usual Bertini theorem (over fields) applies to show that $V$ contains 
the generic point of each 
fiber $Q\to \Spec(\Z)$. Therefore, the complement $Q-V$ is contained
in a closed subscheme $Z$ which is such that $Q-Z\to \Spec(\Z)$ is surjective.
By [Ru] or [MB2] there is a
number field $K$ with integer ring $\O_K$ and an
integral point $\phi:\Spec(\O_K)\to Q-Z\subset Q$.
Since $Q-Z\subset V$, the pull-back of $X(r)\to Y(r)$ via
$Y(r)_\phi\to Y(r)\times \Spec(\O_K)$ now 
gives a $G$-torsor $X\to Y:=Y(r)_\phi$ with $Y\to \Spec(\O_K)$
projective and smooth of relative dimension $d$. The argument in [Se] 
now shows that the generic fiber $X_{K}$ is geometrically connected.

\vfill\eject

\end{document}